\documentclass[10pt,letter]{article}
\usepackage{amsmath}
\usepackage{amssymb}
\usepackage{amsfonts}
\usepackage{mathrsfs}
\usepackage{amsthm}
\usepackage{epsfig,epic}
\usepackage[english]{babel}
\usepackage{graphics}
\usepackage{graphicx}
\usepackage{eucal}
\usepackage{amscd}
\usepackage{amssymb}
\usepackage{fancyhdr}
\usepackage[hang]{subfigure}
\usepackage{array}
\usepackage[nottoc]{tocbibind}

\newtheorem{lemma}{Lemma}
\newtheorem{prop}{Proposition}
\newtheorem{defn}{Definition}
\newtheorem{thm}{Theorem}
\newtheorem{cor}{Corollary}
\newtheorem{rem}{Remark}

\newcommand{\be}{\begin{equation}}
\newcommand{\ee}{\end{equation}}
\newcommand{\bes}{\begin{equation*}}
\newcommand{\ees}{\end{equation*}}

\newcommand{\Sym}[1]{\mathcal{S} _{ #1} }
\newcommand{\ud}{\mathrm{d}}
\newcommand{\BS}[2]{S_{#2}(#1)}

\newcommand{\RL}[1]{A_{#1}}

\newcommand{\RLs}{A}

\newcommand{\R}{\mathcal{R}}
\newcommand{\V}{\mathcal{V}}
\newcommand{\Z}{\mathcal{Z}}
\newcommand{\RLp}[1]{A^{(#1)}}

\newcommand{\ZLp}[1]{A^{(#1)}}
\newcommand{\Orb}[2]{\mathcal{O}_{#2}(#1)}

\newcommand{\I}{\mathcal{I}}
\newcommand{\M}{\mathcal{M}}

\newenvironment{proofof}[2]{\begin{proof}[Proof of #1 \ref{#2}.]}{\end{proof}}

\title{Mixing of asymmetric logarithmic suspension flows over interval exchange transformations.}

\author{Corinna Ulcigrai\footnote{Mathematics Department, Princeton University.
{\it Email: \,}
ulcigrai@math.princeton.edu}}

\date{}

\begin{document}

\maketitle

\begin{abstract}
We consider suspension flows built over interval exchange transformations with the help of roof functions having an asymmetric logarithmic singularity. We prove that such flows are strongly mixing for a full measure set of interval exchange transformations.
\end{abstract}

\section{Introduction.}
\subsection{Motivation and main references.}
Hamiltonian systems with  multi-valued Hamiltonian functions on two dimentional tori 
 give rise 
 to area preserving flows 
 which can be decomposed into a finite number of components filled with periodic trajectories 
and one ergodic component (see \cite{Ar:top}). 
The flow on this ergodic component is isomorphic to a suspension flow built over a rotation of the circle with the help of a roof function which has asymmetric logaritmic singularities (see also Section \ref{defns} for precise definitions). 

The question about mixing of such flows, risen in the same paper \cite{Ar:top}, was answered by Sinai and Khanin in \cite{SK:mix}, where it was proved that, under a generic diophantine condition on the rotation angle, suspension flows with asymmetric singularities over a rotation are strongly mixing (see also \cite{Ka:mix}). 
The diophantine condition of \cite{SK:mix} was weakened by Kochergin 
in a series of works (\cite{Ko:nonI, Ko:nonII, Ko:som, Ko:wel}).

Mixing in these flows is produced by different deceleration rates near the singular points. Neighbouring points on a Poincar\'e transversal have different return times and this causes a phenomenon sometimes called stretching of the Birkhoff sums (the idea of how this stretching leads to mixing is explained in Section \ref{geommech}). A similar stretching of Birkhoff sums was also used by Fayad in \cite{Fa:ana} to construct mixing reparametrization of flows on $\mathbb{T}^3$.  

Mixing does not arise in suspension flows over rotations 
in the case of bounded variation roof functions \cite{Ko:abs}. The presence of a symmetric logarithmic singularity is also not enough, as it was shown by Kochergin in \cite{Ko:non}. 
Lema{\'{n}}czyk \cite{Le:sur} proved the absence of mixing if the Fourier coefficients of the roof function are of order $O(1/|n|)$ and showed with Fr{\c a}czek that these flows are disjoint in the sense of Furstenberg from all mixing flows \cite{FL:acl}. This condition is essentially sharp, see \cite{Ko:ami}. 

Consider, instead of $\mathbb{T}^2$, a compact orientable surface $M_g$ of higher genus ($g\geq 2$). A closed Morse one form $\omega$  
generates a Hamiltonian flow determined by the multivalued Hamiltonian $H$ locally defined by $\ud H = \omega$. 
The corresponding area preserving flow on $M_g$ can be decomposed into 
components filled by periodic orbits and 
components on which the flow is metrically isomorphic to a suspension flow over an interval exchange transformation (IET) (see e.g. \cite{Zo:how}). Interval exchange transformations are piecewise orientation-preserving isometries of an interval which appear naturally as first return maps of such flows on a transversal, as rotations do in the case of $\mathbb{T}^2$. 

It was proved by Katok in \cite{Ka:int} that suspension flows over IETs under roof functions of bounded variation are never mixing and (see   \cite{FL:ond}) are disjoint from mixing flows. 
On the other hand, Kochergin (see \cite{Ko:mix}) proves mixing for a class of roof functions over IETs which includes power-like singularities, which arise when the fixed points on the corresponding surface flow are degenerate. The presence of non-degenerate fixed points give rise to logarithmic singularities.
Fr{\c a}czek and Lema{\'{n}}czyk prove in \cite{FL:ond} that in the case of symmetric logarithmic singularities and typical IETs of $2$ or $3$ intervals 
 the suspension flows are also disjoint from mixing flows. 

In this paper we consider suspension flows over IETs of an arbitrary number of intervals with roof functions having a single asymmetric logarithmic singularity. We prove that for typical IETs such flows are strongly mixing. The case of several asymmetric singularities will be treated in another paper. 

As it was mentioned above the main mechanism of mixing is the stretching of Birkhoff sums. The proof of stretching in our case uses as a tool the Rauzy-Veech renormalization algorithm for IETs (see Section \ref{algorithmssec}). The condition on the IET which guarantees mixing is typical in view of a recent result in \cite{AGY:exp}. 

\subsection{Definitions and Main Result.}\label{defns}
\paragraph{Interval exchange transformations.}
Let $I^{(0)}=[0,1)$ and let $T:I^{(0)}\rightarrow I^{(0)}$ be an \emph{interval exchange transformation} (IET) of $d$ subintervals, i.e. a piecewise orientation preserving isometry of $I^{(0)}$ defined in the following way. Assign a permutation $\pi\in S_{d}$ and a partition of $I^{(0)}$ into $d$ subintervals, $I^{(0)}_1$, $I^{(0)}_2$, $\dots$, $I^{(0)}_{d}$, defined by a lengths vector $\underline{\lambda}=(\lambda_1, \lambda_2, \dots, 
\lambda_{d})$, $\lambda_i >0$, $\sum_{i=1}^{d}\lambda_i=1$, such that $\lambda_i=|I^{(0)}_i|$. Then $T$ permutes the subintervals according to $\pi$ so that under the action the transformation $I^{(0)}_i$ becomes the $\pi(i)$$^{th}$ interval, i.e. the order of the subintervals after applying $T$ is $I^{(0)}_{\pi^{-1}(1)}$,$I^{(0)}_{\pi^{-1}(2)}$ $\dots $ $I^{(0)}_{\pi^{-1}(d)}$. More precisely 
\begin{eqnarray}
I^{(0)}_{j} & \doteqdot &[ \sum_{i=1}^{j-1} \lambda_i , \sum_{i=1}^{j} \lambda_i [ \qquad j=1, \dots , d; \nonumber \\
T(x) & = & x  - \sum_{i=1}^{j-1} \lambda_i + \sum_{i=1}^{j-1} \lambda_{\pi^{-1}i} \quad \mathrm{for} \quad x \in I^{(0)}_{j}, \qquad j=1, \dots , d \nonumber . 
\end{eqnarray}
We shall often use the notation $T=(\lambda, \pi)$. 

\paragraph{Suspension flows.}
Let $f\in L^1 (I^{(0)}, dx)$ be a strictly positive function $f\geq m_f >0$ and assume $\int_{I^{(0)}} f(x) dx =1$. 
Further assumptions on $f$ will be formulated in Section \ref{1asymmetriclog}.
The \emph{phase space $X_f$} of the suspension flow is defined as
\bes
X_f \doteqdot \{ (x,y) |\quad  x \in I^{(0)}, \, 0\leq y < f(x) \}
\ees
and can be depicted as the set of points below the graph of the roof function $f$.
Introduce the normalized measure $\mu$ which is the restriction to $X_f$ of the Lebesgue measure $\ud x\,\ud y$.

The \emph{suspension flow built over $T$ with the help of the roof function f} is a one-parameter group $\{ \varphi_t \}_{t\in \mathbb{R}}$ of $\mu$-measure preserving transformations of $X_f$ whose action is generated by the following two relations:
\begin{equation} \label{suspflow}
 \left\{  \begin{array}{lll}
\varphi_t(x,y) &=& (x, y+t), \qquad \mathrm{if} \, 0\leq y+t < f(x); \\
\varphi_{f(x)}(x,0) &=& ( Tx,0).\\
		     \end{array} \right.
\end{equation}
Under the action of the flow a point of $(x,y) \in X_f$  moves with unit velocity along the vertical line up to the point $(x,f(x))$, then jumps instantly to the point $\left( T(x),0 \right)$, according to the base transformation. Afterwards it continues its motion along the vertical line and so on (see e.g. \cite{CFS:erg}).

We will denote 
by\footnote{The dependence on $T$ is omitted when there is no ambiguity.}
\bes
\BS{f,T}{r}( x) = \BS{f}{r}( x)  \doteqdot \sum_{i=0}^{r-1} f(T^i(x)), \qquad x\in I^{(0)}. 
\ees
the $r^{th}$ non-renormalized \emph{Birkhoff sum} of $f$ along the trajectory of $x$ under $T$. 

Given $x\in I^{(0)}$ denote by $r(x,t)$ the integer uniquely defined by
\be \label{defr}
r(x,t)\doteqdot \max \{ r\in \mathbb{N} \, | \quad  \BS{f}{r}(x) \leq t \},
\ee
which describes the number of \emph{discrete iterations} of the IET which the point $(x,0)$ undergoes before time $t$.
According to this notation the flow $\varphi_t$ defined by (\ref{suspflow}) acts as
\be \label{flowdef}
\varphi_t(x,0) = \left( T^{r(x,t)}(x), t- \BS{f}{r(x,t)}(x)\right).
\ee
For $t<0$, the action of the flow is defined as the inverse map. 


\paragraph{Single asymmetric logarithmic singularity.} \label{1asymmetriclog}
Assume that $f\in \mathscr{C}^2\left((0,1)\right)$ and there exist two positive constants $C^+>0$, $C^- >0$, 
such that
\be  \label{logsingforf''}
\lim_{x\rightarrow 0^+} \frac{f''(x)}{\frac{1}{x^2}}=C^+; \quad \qquad
 \lim_{x\rightarrow 1^-} \frac{f''(x)}{\left(\frac{1}{1-x}\right)^2}=C^- .\ee
It is easy to see that it implies that
\be\label{log}
\lim_{x\rightarrow 0^+} \frac{f(x)}{|\log x|}=C^+; \qquad \quad  \lim_{x\rightarrow 1^-} \frac{f(x)}{|\log (1-x) |}=C^- .
\ee
Hence we say in this case that $f$ has a \emph{logarithmic singularity} at the origin. The singularity is called \emph{asymmetric} if $C^+ \neq C^-$. 

\paragraph{Mixing.}
Recall that a flow $\{ \varphi_t\}_{t\in \mathbb{R}}$ preserving the measure $\mu$ is said to be \emph{mixing} if for each pair of measurable sets $A$, $B$, one has
\be \label{mixingdef}
\lim_{t\rightarrow \infty} \mu(\varphi_t(A)\cap B)=\mu(A)\mu(B).
\ee

\paragraph{Main Result.}
The main result of this paper is the following.
\begin{thm}\label{mixing}
The suspension flow $\{\varphi_t \}_{t\in \mathbb{R}}$ built over a typical IET $T$ with the help of a roof function $f$ having a single asymmetric logarithmic singularity at the origin is mixing. \end{thm}
The notion of \emph{typical} IET is undestood from the measure theoretical point of view. More precisely, for every irreducible $\pi$, Theorem \ref{mixing} holds for a.e. lengths vector $\underline{\lambda}\in \Delta_{d-1}$ with respect to the Lebesgue measure on the simplex $\Delta_{d-1}$. 

\subsection{A criterium for mixing.}\label{mixingcriteriumsec}
\paragraph{Partial partitions and rectangles.} By a \emph{partial partition} $\eta$ of $I^{(0)}$ into intervals we mean a collection of disjoint intervals $I=[a,b[$. We do not require that the union of these intervals is the whole $I^{(0)}$. All the partitions in this paper will be partial partitions into a finite number of intervals.
Denote by $Leb$ the Lebesgue measure on the Borel subsets of $I^{(0)}$. By using the notation $Leb(\eta)$ we mean the total measure of a partition  $\eta$, i.e.   
$Leb(\eta) \doteqdot \sum_{I\in \eta} Leb(I)$.
The \emph{mesh} of the partition $\eta$ is given by 
$mesh(\eta) \doteqdot \sup_{I\in \eta }Leb(I) $.
We will consider 1-parameter families of partial partitions $\eta(t)$, $t\in \mathbb{R}$.

Call \emph{rectangle} of base $b(R)\subset I^{(0)}$ and height $h=h(R) < m_f$ the set $R$ of points $(x,y)$ such that $0\leq y\leq h$ and $x\in b(R)$. 
 Rectangles and their shifts $\varphi_t(R)$ generate the Borel $\sigma$-algebra of  $(X_f,\mu)$.

\paragraph{Mixing criterium.}
In order to show mixing it enough to verify the following criterium, similar to the one also used in \cite{Ko:mix, Fa:ana}.
\begin{lemma}[Mixing criterium] \label{mixingcriterium}
If, given any rectangle $R$, any $\epsilon>0$ and any $\delta >0$, one can find $t_0>0$ such that for each $t\geq t_0$ one can define a partial partition $\eta(t)$ of $I^{(0)}$ into intervals such that
\be \label{partitiontendstopoint}
Leb(\eta(t))> 1-\delta, \qquad mesh (\eta(t) ) \leq \delta 
\ee
and for each $I\in \eta(t)$
\be \label{mainestimate}
Leb (I \cap \varphi_{-t} (R)) \geq (1-\epsilon) Leb (I) \mu(R),
\ee
then the flow $\{ \varphi_t \}_{t\in \mathbb{R}}$ is mixing.
\end{lemma}
\begin{proof}
Mixing means that for any two measurable sets $A$ and $B$ and any $\epsilon>0$, for big enough positive $t$,
\be \label{mixingineq}
\mu(A\cap \varphi_{-t}(B))>(1-\epsilon)\mu(A)\mu(B),
\ee
since applying (\ref{mixingineq}) to $A^{C}$ and $B$ one gets 
\bes 
\mu(A \cap \varphi_{-t}(B))<(1+\epsilon)\mu(B)\mu(A)+ \epsilon\mu(B)
\ees
and therefore (\ref{mixingdef}). For $t<0$, it is enough to exchange the roles of $A$ and $B$ and use $\mu$-invarance of $\varphi_t$.
Moreover, it is enough to verify (\ref{mixingineq}) for $A$ and $B$ rectangles, since any measurable set can be approximated by a finite union of rectangles and their shifts under the flow.

Let $b(A)$ be the base of a rectangle $A$. For each $\delta>0$ and $t\geq t_0$,  
there exists a finite number of intervals $I^{(t)}_{k} \in \eta(t)$, $k=0,\dots, K(t)$, such that
$Leb \left( b(A)\Delta \cup_{k=0}^{K(t)} I^{(t)}_k \right) \leq 3 \delta$. Here $\Delta$ denotes the symmetric difference of sets. To see this, consider all intervals of $\eta(t)$ which intersect $b(A)$ and use (\ref{partitiontendstopoint}).
Let 
\bes
A'\doteqdot  \bigcup_{0\leq y \leq h(A)} \left( \bigcup_{k_y=0}^{K(t+y)} I^{(t+y)}_{k_y}  \times \{ y \} \right).
\ees
Choosing $\delta \leq \frac{ \epsilon }{3} Leb (b(A)) \mu(B)$, 
by Fubini theorem, 
$\mu \left( A \, \Delta \, A' \right) \leq  \epsilon \mu(A)\mu(B)$ .
Remarking the inequality $y\leq h(A) <  m_f$,  
we have, for each slice of $A'$,
\bes
\left( \cup_{k_y=0}^{K(t+y)} I^{(t+y)}_{k_y}   \times \{ y \} \right) \cap \varphi_{-t} (B) = \varphi_{y} \left( \left( \cup_{k_y=0}^{K(t+y)} I^{(t+y)}_{k_y}  \times \{0 \} \right) \cap \varphi_{-t-y} (B) \right) .
\ees
Moreover $\varphi_y$ preserves $Leb$ on each slice and therefore one can assume that the hypothesis (\ref{mainestimate}) in which we set $B=R$ holds for all slices. Thus, combining these estimates and applying again Fubini theorem, we get, for $t\geq t_0$,
\bes
\begin{split}
&\mu ( A \cap \varphi^{-t} (B))  \geq  
\mu \left( 
A' \cap \varphi^{-t} (B) \right) - 3\delta h(A)
\geq  \\ & \int_{0}^{h(A)} (1-\epsilon)  Leb 
  \bigcup_{k_y=0}^{K(t+y)} I^{(t+y)}_{k_y} 
\mu(B) dy - \epsilon \mu(A) \mu(B) 
 \geq \left ((1-\epsilon)^2 - \epsilon 
\right) \mu(A) \mu(B),
\end{split}
\ees
hence proving the lemma.
\end{proof}

\paragraph{Intuitive explanation of the mixing mechanism.} \label{geommech}
In this section we explain the geometric mechanism that produces mixing for this type of suspension flows. 
\begin{figure}
\centering
\subfigure[ $\varphi_t(I)$ ]{\label{verticalstripg}
\includegraphics[width=0.55\textwidth]{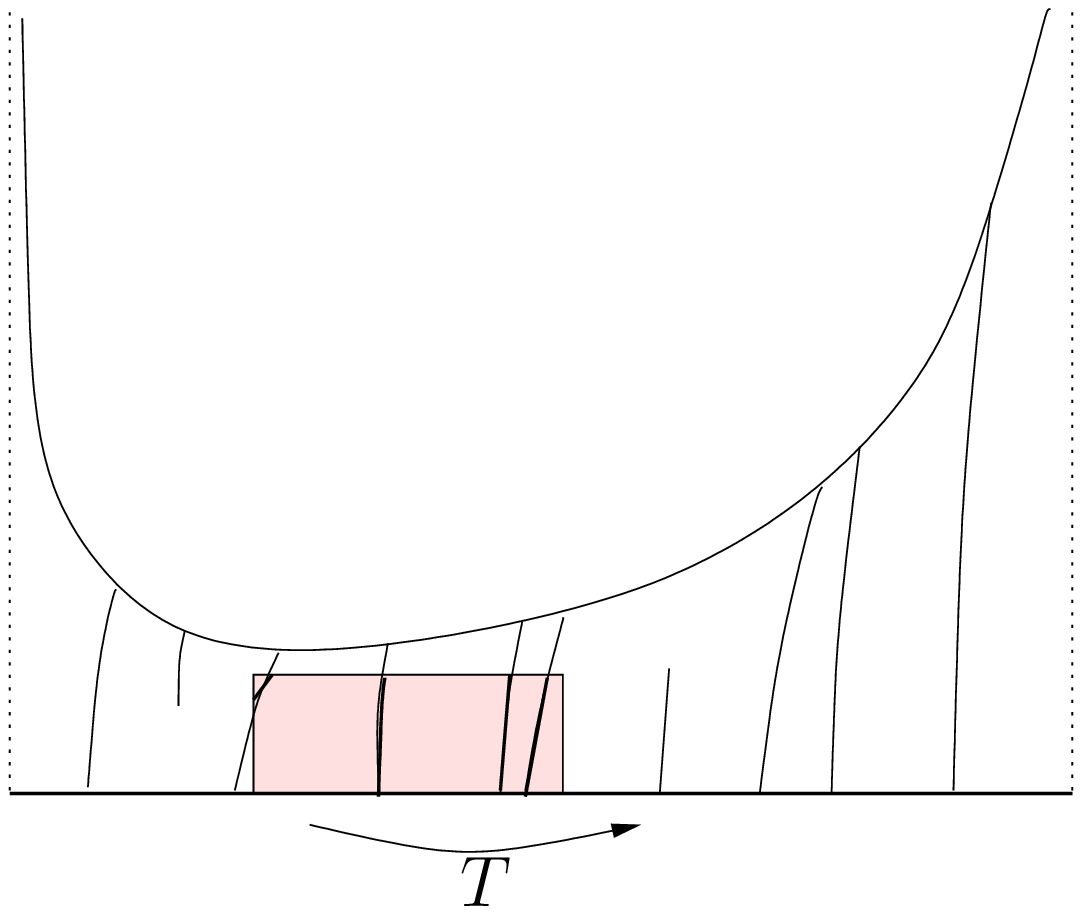}}
\hspace{5mm}
\subfigure[ Toy model ]{\label{toyflow}
\includegraphics[width=0.3\textwidth]{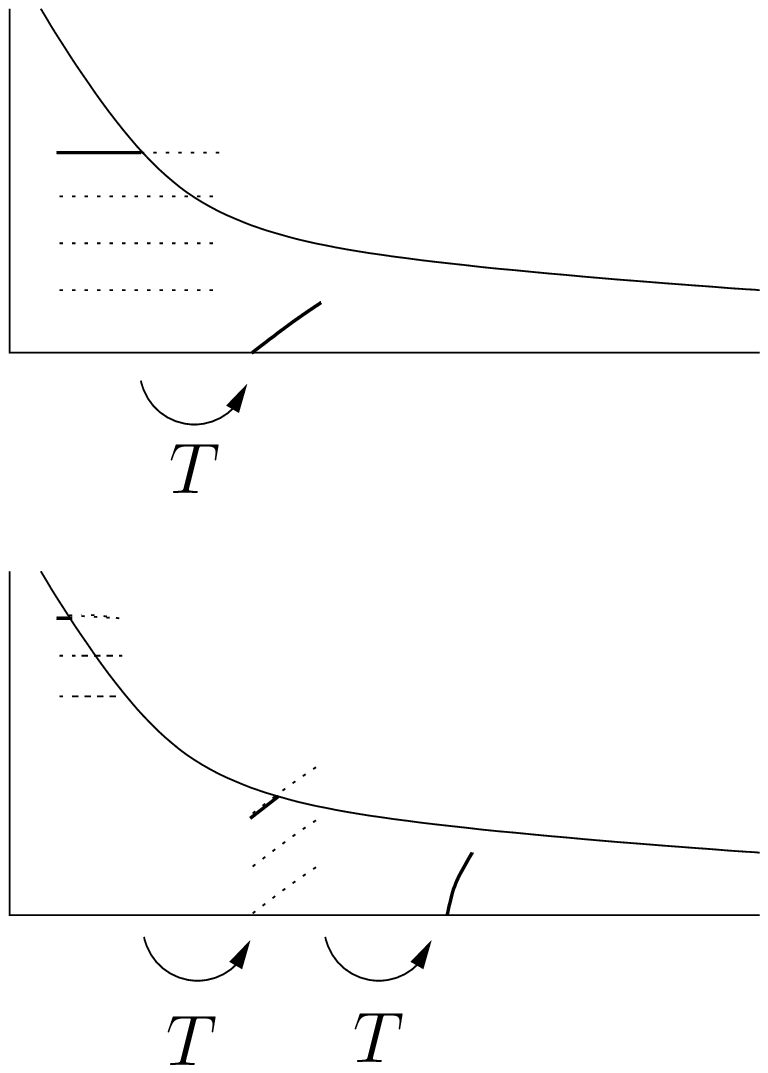}}
\end{figure}

Consider a sufficiently small segment $I=[a,b]\subset I^{(0)}$ and let us understand how its image under the flow, $\varphi_t(I)$, looks like for very large $t$. We claim that $\varphi_t(I)$ will consist of many almost vertical curves, as shown in Figure \ref{verticalstripg}. 

Assume as a toy example that $f$ has only a one-sided logarithmic singularity at the origin and is monotonically decreasing, as in Figure \ref{toyflow}. Notice first that until $t < m_f$, $\varphi_t(I)$ is still a horizontal segment, while as $t= f(x_0) $ for some $x_0\in I$, $\varphi_t(I)$ splits into two curves: one is a still an horizontal segment, while the other piece will project over $T ([x_0,b])$ and be a translate of the graph of $-f|_{[x_0,b]}$, as can be seen by (\ref{flowdef}). See Figure \ref{toyflow}. More generally, from (\ref{flowdef}), each of the curves in which $\varphi_t(I)$ will split is a graph of a translate of the Birkhoff sum $\BS{f}{r}$ restricted over a small interval of the form $T^{r}([x_i,x_{i+1}))$, where $[x_i,x_{i+1}) \subset I$. Noticing that $f'<0$ and the integral of $f'$ is divergent, one can prove in this toy model that the slopes of these curves, which are given by $-\BS{f'}{r}$, are growing to infinity, i.e. they are becoming almost vertical. Hence the increasingly big delay between different points causes $\varphi_t(I)$ to split into many curves, which are distributed over the orbit of $T^n(I)$. Using unique ergodicity of $T$ on the base and the fact that each strip can be approximated by a straight line, one can show that the fraction of $x\in I$ such that $\varphi_t(x)\cap R \neq \emptyset$ is proportional to $\mu (R)$.

When the singularity is asymmetric the same phenomenon happens and one can show that delays accumulated from visits on one side are stronger than opposite delays accumulated from the other, causing $\BS{f'}{r}(x)$ to diverge as in the presence of a one-sided singularity 
 for most of the points. 

\paragraph{Outline of the proof of Theorem \ref{mixing}.}
In order to prove mixing for the suspension flow, we use the criterium in Lemma \ref{mixingcriterium}. Given a rectangle $R$ and $\epsilon , \delta >0$,
our goal is to construct, for any sufficiently large $t$, a partial partition $\eta(t)$ of $I^{(0)}$ into intervals which satisfy (\ref{partitiontendstopoint}, \ref{mainestimate}). Each of the intervals of these partitions behaves under the flow as explained in the previous paragraph. The construction of the partition is carried out in several steps, explained in Section \ref{mixingpartitionsec}.
In order to get the final estimate (\ref{mainestimate}), in Section \ref{areaestimatesec}, the key step is to get a good estimate of the rate of growth of the first two derivatives of $\BS{f}{r}$. Such estimates, presented in Section \ref{growthBSsec}, are based on some property of the renormalization cocycle for IETs introduced by Rauzy, Veech and Zorich. The definition and some properties of this cocycle are recalled in Section \ref{algorithmssec}.


\section{Renormalization algorithms for IETs.} \label{algorithmssec}
Rauzy, Veech and Zorich (see \cite{Ra:ech, Ve:pro,{Zo:fin}}) developed a renormalization algorithm for IET which is a multi-dimensional generalization of the continued fraction algorithm. 

In what follows let $T=(\underline{\lambda}, \pi)$ be an IET. We assume that $\pi$ is \emph{irreducible}, i.e. if the subset $\{ 1, 2, \dots i \}$ is $\pi$-invariant, then $i=d$, since this is a necessary condition for minimality.
We also assume that $T$ satisfies the Infinite Distinct Orbit Condition (IDOC) introduced by M. Keane in \cite{Ke:int}. 
We say that $T$ has the (IDOC) if, denoting by $\beta_0 =0$ and $\beta_j \doteqdot \sum_{i=1}^{j} \lambda_i$ for $j=1, \dots , d$
 the discontinuities of $T$, the orbits $\mathscr{O}(\beta_j) \doteqdot \{ T^{n}(\beta_j) |\, n \in \mathbb{N}\}$, $1\leq j\leq d-1$, are infinite and disjoint, i.e. $\mathscr{O}(\beta_j)\cap \mathscr{O}(\beta_i)=\emptyset$ for any $i\neq j$. 
As it was shown by Keane in \cite{Ke:int}, the (IDOC) implies minimality.

\subsection{Rauzy-Veech and Zorich algorithms and cocycles.}\label{induction}
Starting with $T=T^{(0)}$, the Rauzy-Veech algorithm produces a sequence of IET $T^{(r)}$ which are induced maps of $T$ onto a sequence of nested subintervals $I^{(r)}\subset I$. It is easy to see that in general the induced first return map of $T$ on a subinterval $I'\subset I$ is again an IET, of at most $d + 2$ intervals. One Rauzy-step is defined so that $T^{(1)}$ is an exchange of exactly the same number $d$ of subintervals.

\subsubsection{The induction algorithms.}
\paragraph{One step of Rauzy-Veech algorithm.}
At the first step, compare the lengths of $I_{d}$ and $I_{\pi^{-1}(d)}$, i.e. of the last subintervals before and after the transformation. It follows from the (IDOC) that $\lambda_{d} \neq \lambda_{\pi^{-1}d}$. Hence there can be two cases:
\begin{itemize}
\item[(a)]
$\lambda_{d} < \lambda_{\pi^{-1}d}$. In this case we consider the new interval $I^{(1)}\doteqdot [0, 1-\lambda_{d}[$. Define $T^{(1)}$ to be the induced map, i.e. the first return map of $T^{(0)}$ onto $I^{(1)}$. It is important that it is again an IET of the same number $d$ of exchanged intervals.

\item[(b)]
$\lambda_{d} > \lambda_{\pi^{-1}d}$. In this case we consider the new interval $I^{(1)}\doteqdot [0, 1-\lambda_{\pi^{-1}d}[$ and, as before, define $T^{(1)}$ to be the induced map on $I^{(1)}$. Also in this case $T^{(1)}$ is again an IET of $d$ intervals.
\end{itemize}

\begin{figure}
\centering
\subfigure[Case $\lambda_{d} < \lambda_{\pi ^{-1}(d)}$]{\label{figTa}\includegraphics[width=0.44\textwidth]{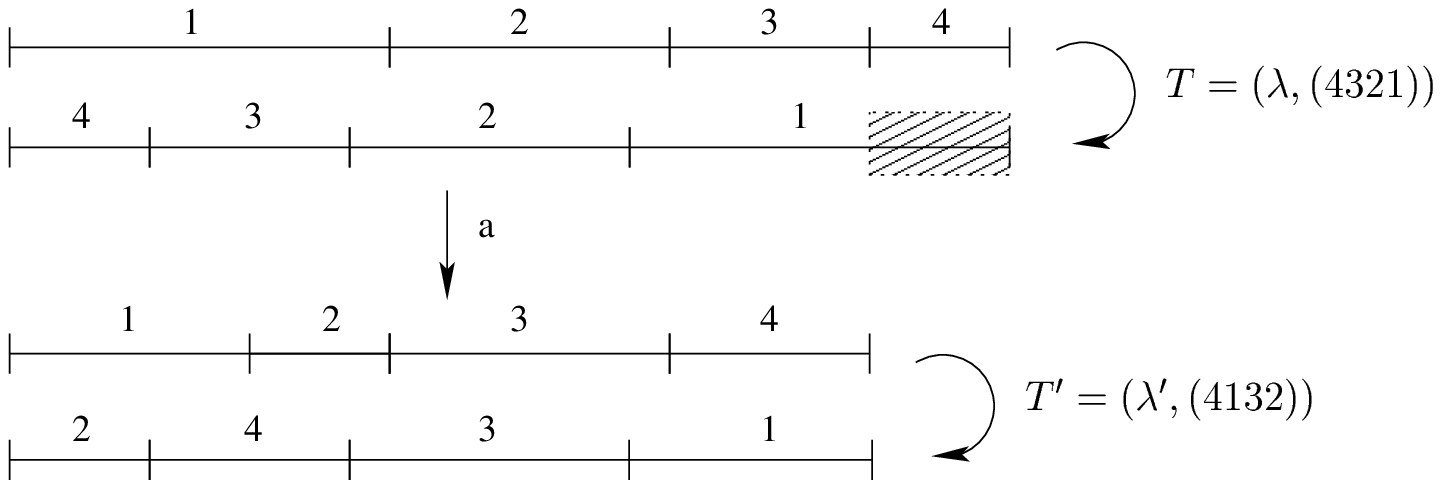}}
\hspace{5mm}
\subfigure[Case $\lambda_{d} > \lambda_{\pi ^{-1}({d})}$]{\label{figTb}\includegraphics[width=0.44\textwidth]{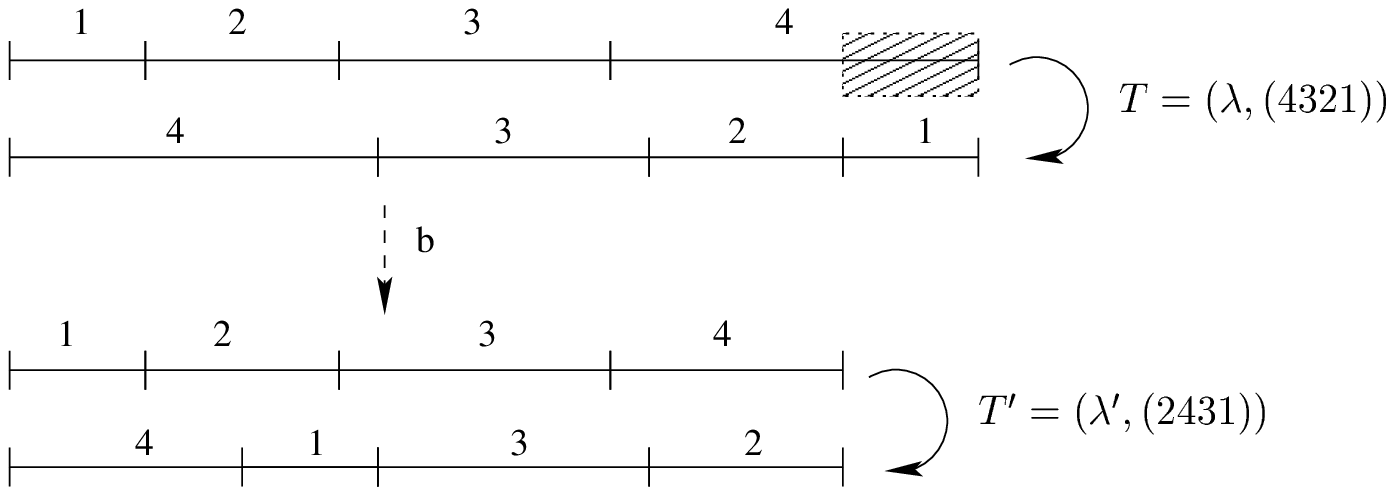}}
\caption{One step of Rauzy-Veech algorithm on $T=(\lambda, (4321))$.}
\label{passoVeechfig2}
\end{figure}

Since $T^{(1)}$ is again an exchange of $d$ intervals, we can write $T^{(1)}=(\lambda^{(1)}, \pi^{(1)})$ (see Figure \ref{figTa},\ref{figTb}), defining in this way a new lengths vector and a new permutation.
One can explicitly write the expressions for two combinatorial operators $a$ and $b$ on $\Sym{d}$, where $\Sym{d}$ is the space of permutations of $d$ elements, such that $\pi^{(1)} = a \pi$ or $ b \pi$ respectively. Explicitly:
\bes
a\pi (j)=  \left\{  \begin{array}{ll}
			\pi (j) & j\leq \pi^{-1}(d) ; \\
			\pi (d) & j = \pi^{-1}(d) + 1 ;\\
			\pi (j-1) & \mathrm{otherwise} ;  \\
		     \end{array} \right.
\ees
\bes
b\pi (j)=  \left\{  \begin{array}{ll}
			\pi (j) & j\leq \pi(d) ; \\
			\pi (j) + 1 & \pi(d) < \pi (j) < d ;\\
			\pi (d) + 1 & \pi (j) = d . \\
		     \end{array} \right.
\ees
Introduce the following matrices to describe the new lengths. Denote by {\it Id}
 the identity $d \times d$ matrix and by $E_{i,j}$ the matrix whose only non zero entry is $(E_{i,j})_{ij}=1$. Introduce the auxiliary permutation $\tau_s \in S^{d}$, $\tau_s = (1\, \, 2\, \, \dots \, s \, \,  s\!+\!2\, \, \,  \dots \, \,  d\, \, s \!+\!1 )$ if $ 1\leq s <d-1$ and $\tau_{d-1} =   \mathrm{id}$, which rotates cyclically all elements after the $s^{th}$-one. Denote by $P(\tau_s )$ the matrix associated to the permutation, i.e. $P(\tau_s )_{ij}=\delta _{i \tau_s(i)}$.
The two \emph{Rauzy-Veech elementary matrices} associated to $(\underline{\lambda},\pi)$ are defined by
\be\label{matriciaeb}
\left\{  \begin{array}{rcl}
			A(\pi ,a) &= & ({\it Id} + E_{\pi ^{-1}(d),d} ) \cdot P(\tau _{\pi ^{-1}(d)} ) ; \\
			A(\pi ,b) &= & {\it Id} + E_{d,\pi^{-1}(d)} . \\
		     \end{array} \right.
\ee

\noindent The induced IET $T^{(1)}$ is then given by
\begin{equation} \label{T1}
(\underline{\lambda}^{(1)} , \pi^{(1)} ) \doteqdot \left\{  \begin{array}{ll}
			\left( A^{-1}(\pi , \, a)\cdot \underline{\lambda} ,\, a(\pi ) \right), & \lambda_{d} < \lambda_{\pi^{-1}({d})} ; \\
			\left( A^{-1}(\pi , \, b)\cdot \underline{\lambda} ,\, b(\pi ) \right), & \lambda_{d} > \lambda_{\pi^{-1}({d})} . \\	
		     \end{array} \right.
\end{equation}
Remark that both $A(\pi ,a)$ and $A(\pi ,b)$ belong to $SL(d, \mathbb{Z})$ and have non-negative entries.

Define inductively $T^{(r)}=(\underline{\lambda}^{(r)}, \pi^{(r)})$ to be the induced map of $T^{(r-1)}$ on $I^{(r)}$. It can be seen that the (IDOC) assures that the algorithm is well defined at each step, i.e. that $\lambda^{(r-1)}_{d} \neq \lambda^{(r-1)}_{ (\pi^{(r-1)} ) ^{-1}(d)}$ for any $r \in \mathbb{N}$.

\paragraph{Renormalized Rauzy-Veech map.}
The \emph{Rauzy class} of $\pi$, denoted by $\mathscr{R}(\pi)\subset S_{d}$, is the set of all permutations obtained iterating the operators $a$ and $b$ starting from $\pi$. 
Using the norm $|\underline{\lambda}|=\sum_{i=1}^{d}\lambda_i$, assume that the initial lengths belong to the simplex $\Delta_{d-1}$ of vectors $\underline{\lambda} \in \mathbb{R}^{d}_{+}$ such that $|\underline{\lambda}|=1$. 
Let us denote by $\Delta(\mathscr{R}) = \Delta_{d-1} \times \mathscr{R}(\pi)$ the space of IETs on the unit interval corresponding to a given Rauzy class $\mathscr{R}$. 

Consider the map on  $\Delta(\mathscr{R})$
which associates to $T$ the induced IET after one step of the algorithm including the following renormalization:
\bes \R \left ( (\underline{\lambda}^{(0)},\pi^{(0)}) \right) \doteqdot \left( \frac{\underline{\lambda}^{(1)}}{|\underline{\lambda}^{(1)}|},\pi^{(1)} \right).
\ees 
Let us call it the \emph{Rauzy-Veech map}. 
Veech proved that $\mathscr{R}$ admits an invariant measure $\mu_{\V}$, absolutely continuous with respect to the Lebesgue measure, which is infinite. 
The main result proved by Veech in \cite{Ve:gau} is that the map $\R$ is conservative. As a consequence, he proves that given $\pi \in \Sym{d}$, for a.e. $\underline{\lambda} \in \Delta_{d-1}$, the IET $T=(\underline{\lambda}, \pi)$ is uniquely ergodic.


\paragraph{Zorich acceleration.}
Take an IET $T$ and consider its Rauzy-Veech orbit $\{\R^n T \}_{n\in \mathbb{N}}$. 
In a typical situation one can find an integer $z_o=z_o(T)>0$ so that $T, \R T,  \dots , \R^{z_0-1}(T) $ all correspond to the same case $(a)$ or $(b)$ while $\R^{z_o}(T) $ corresponds to the other one. Grouping together these $z_0$ steps of Rauzy induction, we get a new transformation $\Z$ on the space of IET, where the letter $\Z$ is chosen in honor of A. Zorich who introduced this map in \cite{Zo:fin}. 
Zorich showed in \cite{Zo:fin} that $\Z$  has an absolutely continuous \emph{finite} invariant measure. 
We will denote the \emph{Zorich invariant measure} by $\mu_{\Z}$ . 

\subsubsection{Rauzy-Veech lengths cocycle.}
As was explained above, to each $T$ one can associate an elementary matrix $A(T)$ in $SL(d, \mathbb{Z})$ defining $A(T)\doteqdot A(\pi ,a)$ or $A(T)\doteqdot A(\pi ,b)$ respectively. 
Let $\RL{r} = \RL{r} (T) \doteqdot  A(\R^r T)$. Then for each $r$ we can associate to $T$ the product:
\bes
\RLp{r} \doteqdot \RL{1} \dots \RL{r}. 
\ees
We can easily see that the map  $A^{-1}$: $\Delta(\mathscr{R}) \rightarrow SL(d,\mathbb{Z})$ is a cocycle over $\R$, which we call the \emph{Rauzy-Veech lengths cocycle}.
Iterating the lengths relation in (\ref{T1}) we get the formula for the lengths vector of $T^{(r)}$:
\begin{equation} \label{lengthsrelationR}
\underline{\lambda} ^{(r)} = \left( \RLp{r} \right) ^{-1 }  \underline{\lambda} .
\end{equation}

Let us also introduce the following notation useful to consider more generally products of Rauzy-Veech cocycle matrices from $m$ to $n$, $m<n$:
 \bes \RLp{m,n} \doteqdot \RL{m} \cdot  \RL{m+1} \cdot \, \dots \, \cdot \RL{n-2} \cdot \RL{n-1}.
\ees

\subsubsection{Hilbert metric and projective contractions.} \label{Hilbertsec}
Consider on the simplex $\Delta_{d-1} \subset \mathbb{R}_+^d$ the \emph{Hilbert distance} $d_H$, defined as follows.
\be
d_H(\lambda, \lambda') \doteqdot \log \frac{\max_{i=1,\dots , d}\frac{\lambda_i}{\lambda'_i}}{\min_{i=1 , \dots ,d}\frac{\lambda_i}{\lambda'_i}}.
\ee
We denote the diameter with respect to $d_H$ of a projective subset $\Lambda \subset \Delta_{d-1}$ by
\be \label{diameter}
diam_H (\Lambda ) \doteqdot \sup_{\lambda, \lambda' \in \Lambda} d_H( \lambda, \lambda' ) .
\ee
Remark that if its closure $\overline{\Lambda} \subset \Delta_{d-1}$, 
  then $diam_H (\Lambda )$ is finite. 

Let us write $A\geq 0$ if $A$ has non negative entries and $A>0$ is $A$ has strictly positive entries.
Recall that to each $A\in SL(d,\mathbb{Z})$, $A\geq 0$, one can associate a projective transformation $\widetilde{A}: \Delta_{d-1} \rightarrow \Delta_{d-1}$ given by
\bes
\widetilde{A} \lambda = \frac{A \lambda }{|A \lambda |}.
\ees
When $A \geq 0$, $d_H(\widetilde{A} \lambda, \widetilde{A} \lambda' ) \leq d_H(\lambda, \lambda')$. Furthermore, if $A>0$, then we get a contraction. More precisely, $A>0$ is equivalent to the closure $ \overline{ \widetilde{A} \left( \Delta_{d-1} \right)}$ being contained in $\Delta_{d-1}$, hence defining 
\be \label{D_Adef}
D(A) \doteqdot diam_H \left( \, \overline{ \widetilde{A}  \left( {\Delta}_{d-1} \right) }\, \right) , 
\ee
we have $D(A) < \infty$. Then
\be \label{contraction}
 d_H(\widetilde{A} \lambda, \widetilde{A} \lambda' ) \leq (1-e^{-D(A)}) d_H(\lambda, \lambda').
\ee



\paragraph{Paths on Rauzy classes.}\label{rauzyclasses}
Rauzy classes can be visualized in terms of directed labeled graphs, the \emph{Rauzy graphs}. Vertices are in one-to-one correspondence with permutations of  $\mathscr{R}(\pi)$; arrows connect permutations obtained one from the other by applying $a$ or $b$ and are labeled according to the type, $a$ or $b$ respectively. 
Each vertex is the starting point and the ending point of exactly two arrows, one of each type. We will denote by $\gamma_i=\gamma_i(\pi', a)$ ($\gamma_i(\pi', b)$) the arrow of type $a$ (type $b$) coming out from the vertex $\pi'$. 
 
A \emph{path} $\gamma = (\gamma_1, \dots ,\gamma_r)$ is a sequence of compatible arrows on the Rauzy graph, i.e. such that the starting vertex of $\gamma_{i+1}$ is the ending vertex of $\gamma_{i}$, $i=1,\dots, r-1$. 
Given a path $\gamma$, we can associate to it a matrix 
\bes
A(\gamma) \doteqdot  A(\gamma_1)\cdot \, \dots \, \cdot A(\gamma_r),
\ees
where  $A(\gamma_i)= {A(\pi' , a)}$ if $\gamma_i=\gamma_i(\pi_i, a)$ and ${A(\pi_i, b)}$ if $\gamma_i=\gamma_i(\pi_i, b)$. Associate to $\gamma$ also the subsimplex:
\be
\Delta(\gamma) \doteqdot \left\{ \widetilde{A(\gamma)} \, \underline{\lambda}\, | \quad  \underline{\lambda} \in \Delta_{d-1}  \right\} \subset \Delta_{d-1} .
\ee
Using the induction, one can easily verify the following.
\begin{rem}\label{sequence}
If $T=(\underline{\lambda}, \pi)$ and $\underline{\lambda} \in \Delta(\gamma)$ where $\gamma= (\gamma_1,\dots, \gamma_r)$ is a path starting at $\pi$, $\gamma_i =\gamma_i (\pi_i,c_i)$, $c_i\in\{ a,b \}$,  
the sequence of types and permutations obtained in the first $r$ steps of Rauzy-Veech induction is determined by $\gamma$, 
 i.e. $\RL{i}(T)=A(\gamma_i)$ and $\pi^{(i)} = \pi_i$. 
\end{rem}

\subsubsection{The natural extension.}
The natural extension $\hat{\R}$ of the map $\R$ was introduced by Veech \cite{Ve:gau} and admits a geometric interpretation in terms of the space of zippered rectangles. 
We use the simpler choice of coordinates for zippered rectangles, adopted by \cite{Bu:dec, MMY:coh, AGY:exp}. 

Consider the following polyhedral cones $\Theta_{\pi}\subset \mathbb{R}^{d}$, where $\pi \in \mathscr{R}$.
\bes
\Theta_{\pi} \doteqdot \{ \underline{\tau}=(\tau_1,\dots ,\tau_d)  \in   \mathbb{R}^{d}| \quad \sum_{i=1}^{k} \tau_i >0 , \, \sum_{i=1}^{k} \tau_{\pi^{-1}i} <0, \, k=1,\dots, d-1 \},
\ees
which is non-empty since if $\tau_i\doteqdot \pi(i) -i$, $\underline{\tau}\in \Theta_{\pi}$.

The real valued function $Area(\cdot)$ associates to $(\underline{\lambda}, \pi, \underline{\tau}) \in   \Delta_{d-1} \times \{\pi\} \times \Theta_{\pi}$, 
\bes
Area(\underline{\lambda}, \pi, \underline{\tau} ) \doteqdot \sum_{k=1}^{d} \lambda_k \, \left( \sum_{i=1}^{k-1} \tau_i  - \sum_{i=1}^{\pi(k)-1} \tau_{\pi^{-1}i}  \right). 
\ees
$Area(\cdot)$ has a geometric interpretation as the area of the zippered rectangle associated to the datas $(\underline{\lambda}, \pi, \underline{\tau})$ (see e.g. \cite{MMY:coh}). 

Consider the following space as domain of the natural extension.
\bes
\hat{\Upsilon}^{(1)}_{\mathscr{R}}\doteqdot \{ (\underline{\lambda}, \pi, \underline{\tau} ) | \quad (\underline{\lambda}, \pi) \in \Delta(\mathscr{R}), \, \underline{\tau} \in \Theta_{\pi} ,\, Area((\underline{\lambda}, \pi, \underline{\tau} ) ) = 1   \}.
\ees
The map $\hat{\R}:\hat{\Upsilon}^{(1)}_{\mathscr{R}} \rightarrow  \hat{\Upsilon}^{(1)}_{\mathscr{R}} $ is defined as follows.
\footnote{More precisely $\hat{\R}$ is defined on triples $(\underline{\lambda}, \pi, \underline{\tau} )$ such that $(\underline{\lambda}, \pi)$ belong to the domain of $\R$.}
\bes
\hat{\R} \left( (\underline{\lambda}^{(0)}, \pi^{(0)}, \underline{\tau}^{(0)} ) \right) = \left( \R (\underline{\lambda}, \pi),  |\underline{\lambda}^{(1)}| \underline{\tau}^{(1)}  \right) = \left( \frac{\underline{\lambda}^{(1)}}{|\underline{\lambda}^{(1)}|},\pi^{(1)} ,|\underline{\lambda}^{(1)}|\underline{\tau}^{(1)}  \right) ,
\ees
where $(\underline{\lambda}^{(1)},\pi^{(1)})$ is defined in (\ref{T1}) and, analogously, 
\bes
\underline{\tau}^{(1)}  \doteqdot \left\{  \begin{array}{ll}
			 A^{-1}(\pi , \, a)\cdot \underline{\tau} , & \lambda_{d} < \lambda_{\pi^{-1}({d})} , \\
			A^{-1}(\pi , \, b)\cdot \underline{\tau} , & \lambda_{d} > \lambda_{\pi^{-1}({d})} . \\	
		     \end{array} \right.
\ees
The map $\hat{\R}$ preserves an invariant measure $\hat{m}$ which is the restriction to $\hat{\Upsilon}^{(1)}_{\mathscr{R}}$ of the Lebesgue measure. Denote by $p$ the projection 
\bes
p: \hat{\Upsilon}^{(1)}_{\mathscr{R}}  \rightarrow \Delta(\mathscr{R}), \qquad p  (\underline{\lambda}, \pi, \underline{\tau} )  = (\underline{\lambda}, \pi ).
\ees
The measure $p\hat{m}$ is absolutely continuous w.r.t. Lebesgue on $\Delta(\mathscr{R})$ and it is exactly the $\R$-invariant measure $\mu_{\V}$ constructed by Veech.

If $\gamma$ is a path on $\mathscr{R}$, starting at $\pi$, denote by 
\bes
\Theta(\gamma) \doteqdot \left\{ {A(\gamma)}^{-1}\, \underline{\tau}\, | \quad  \underline{\tau} \in \Theta_{\pi}  \right\} \subset \mathbb{R}^{d}.
\ees
If $\gamma$ is an arrow starting at $\pi$ and ending at $\pi'$, then $\hat{\R}$ maps
\bes \left( \Delta(\gamma) \times \{\pi\} \times \Theta_{\pi}\right) \cap  \hat{\Upsilon}^{(1)}_{\mathscr{R}}  \xrightarrow{\hat{\R}} \left(\Delta_{d-1} \times \{\pi' \} \times \Theta(\gamma) \right) \cap  \hat{\Upsilon}^{(1)}_{\mathscr{R}} .
\ees
As $\underline{\lambda}$ determines the \emph{future} induction steps 
(see Remark \ref{sequence}), similarly $\underline{\tau}$ determines the \emph{past} ones. 
More precisely, let $(\underline{\lambda}^{(-i)}, \pi^{(-i)}, \underline{\tau}^{(-i)} ) \doteqdot \hat{\R}^{-i}(\underline{\lambda}, \pi, \underline{\tau} )$, for $i\in \mathbb{N}$.
\begin{rem}\label{rempast}
If $\pi'$ is the ending vertex of $\gamma=(\gamma_1,\dots, \gamma_r)$ where $\gamma_i=\gamma_i(\pi_i,c_i)$, $c_i \in \{ a,b \}$, and $\underline{\tau} \in \Theta(\gamma)$, the sequence of types and permutations obtained in the past $r$ steps of $\hat{\R}$ is determined by $\gamma$, i.e.  
$\RLs (\lambda^{(-i)}, \pi^{(-i)} )=A(\gamma_{r-i+1})$ and $\pi^{(-i)} = \pi_{r-i+1}$ for $i=1,\dots,r$.
\end{rem}


\subsection{Towers construction and heights vectors.}\label{towers}
The initial interval exchange $T$ can be seen as a suspension over each of the induced $T^{(r)}$ obtained at the $r^{th}$ step of Rauzy-Veech algorithm. In this subsection we define the towers which allow to retrieve $T$ from $T^{(r)}$ and $\RLp{r}$.

Note that the entries of $\RLp{r}$ have a \emph{dynamical meaning} in terms of return times. Namely, denote $I^{(r)}_j$, $1\leq j\leq d$, the subintervals of $T^{(r)}$. 
\begin{rem} \label{dynmean}
The entry $\RLp{r}_{ij}$ is equal to the number of visits of the orbit of any point $x\in I^{(r)}_j$ to the interval $I^{(0)}_i$ of the original partition before its first return in $I^{(r)}$. 
\end{rem}
Therefore, the norm $h^{(r)}_j$ of the $j^{th}$ column of $\RLp{r}$, i.e. $h^{(r)}_j \doteqdot \sum_{i=1}^{d} \RLp{r}_{ij}$ gives the return time of any $x\in I^{(r)}_j$ to $I^{(r)}$.


\paragraph{The towers.}\label{towersec}
Define 
\begin{equation}
Z^{(r)}_j \doteqdot \bigcup _{l=0}^{h^{(r)}_j-1} T^l I^{(r)}_j.
\end{equation}
When $T$ is ergodic, $\bigcup_{j=1}^{d} Z^{(r)}_j$ is a non-trivial $T$-invariant set, therefore the sets $Z^{(r)}_j$, $1\leq j\leq d$ give a partition of the whole $I$. Each $Z^{(r)}_j$ can be visualized as a tower over $I^{(r)}_j\subset I^{(r)}$, of height  $h^{(r)}_j$ (see Figure \ref{stacking}). A floor of the tower, denoted by  $Z^{(r)}_{j,l}$, is defined by
$Z^{(r)}_{j,l} \doteqdot  T^l I^{(r)}_j$ , 
$l=0, \dots, h^{(r)}_j-1$.
The original $T$ is an integral map over $I^{(r)}$; under the action of $T$ every floor $Z^{(r)}_{j,l}$, but the top one ($l\neq h^{(r)}_j$), moves one step up, while $T \left( Z^{(r)}_{j,h^{(r)}_j} \right) = T^{(r)} \left( I^{(r)}_j \right)$.

\paragraph{Heights cocycle.} \label{Rauzyheights}
Let $\underline{h}^{(0)}$ be the column vector $\underline{e} \doteqdot (1, \dots , 1)^T \in \mathbb{Z}^d$ and $\underline{h}^{(n)}$ the column vector whose components are the heights $(h ^{(1)}_{1}, \dots , h ^{(1)}_{d})^T$ of the towers after the first step of the Rauzy-Veech algorithm. If we write 
$\underline{h}^{(1)} = R \,\, \underline{h}^{(0)}$,
where $R=R(T)$ is a matrix in $SL(d, \mathbb{Z})$, it is easy to see that $R(T) = A(\pi,a)^T$ or $R(T) = A(\pi,b)^T$, depending on whether the Rauzy step is of type $(a)$ or $(b)$.

Hence, comparing with (\ref{lengthsrelationR}), the cocycle that determines how the vectors of the heights transform is given by the inverse transpose of the Rauzy-Veech cocycle. More precisely, if $\underline{h} ^{(r)}$ is the vector of the heights after $r$ iterations of the algorithm, then
\be \label{hZtrasf}
\underline{h}^{(r)} = (\RLp{r})^T \, \,  
\underline{e}; \qquad\underline{h}^{(s+r)} = {\RLp{r}(\R^s(T) )}^T \underline{h}^{(s)} \,  .
\ee



\subsubsection{Recurrent structure of the towers.}\label{partitionsec}
\begin{figure}
\centering
\subfigure[$\lambda_{d} < \lambda_{\pi ^{-1}(d)}$]{\label{stacka}{\includegraphics[width=0.45\textwidth]{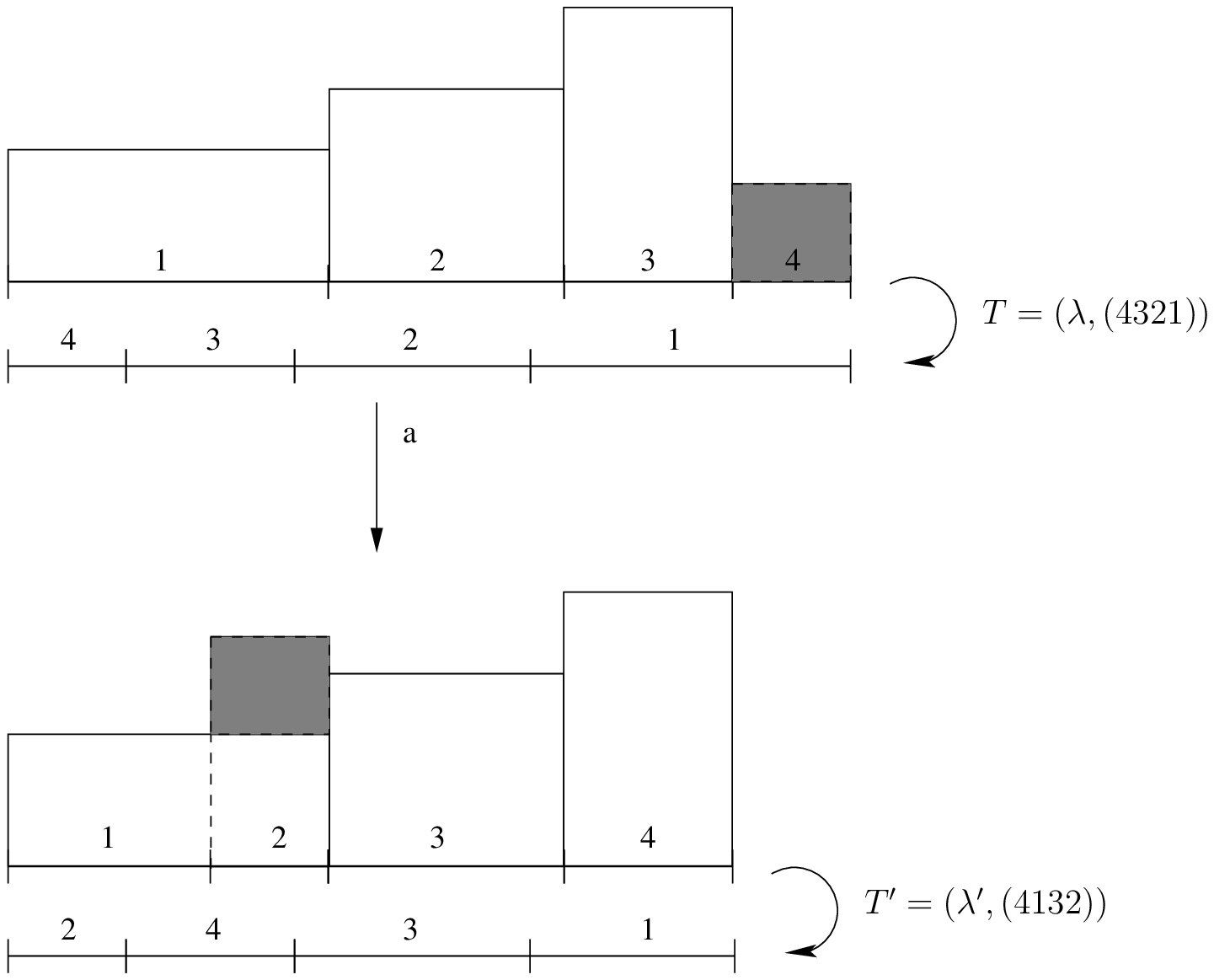}}}
\hspace{5mm}
\subfigure[$\lambda_{d} > \lambda_{\pi ^{-1}(d)}$]{\label{stackb}{\includegraphics[width=0.45\textwidth]{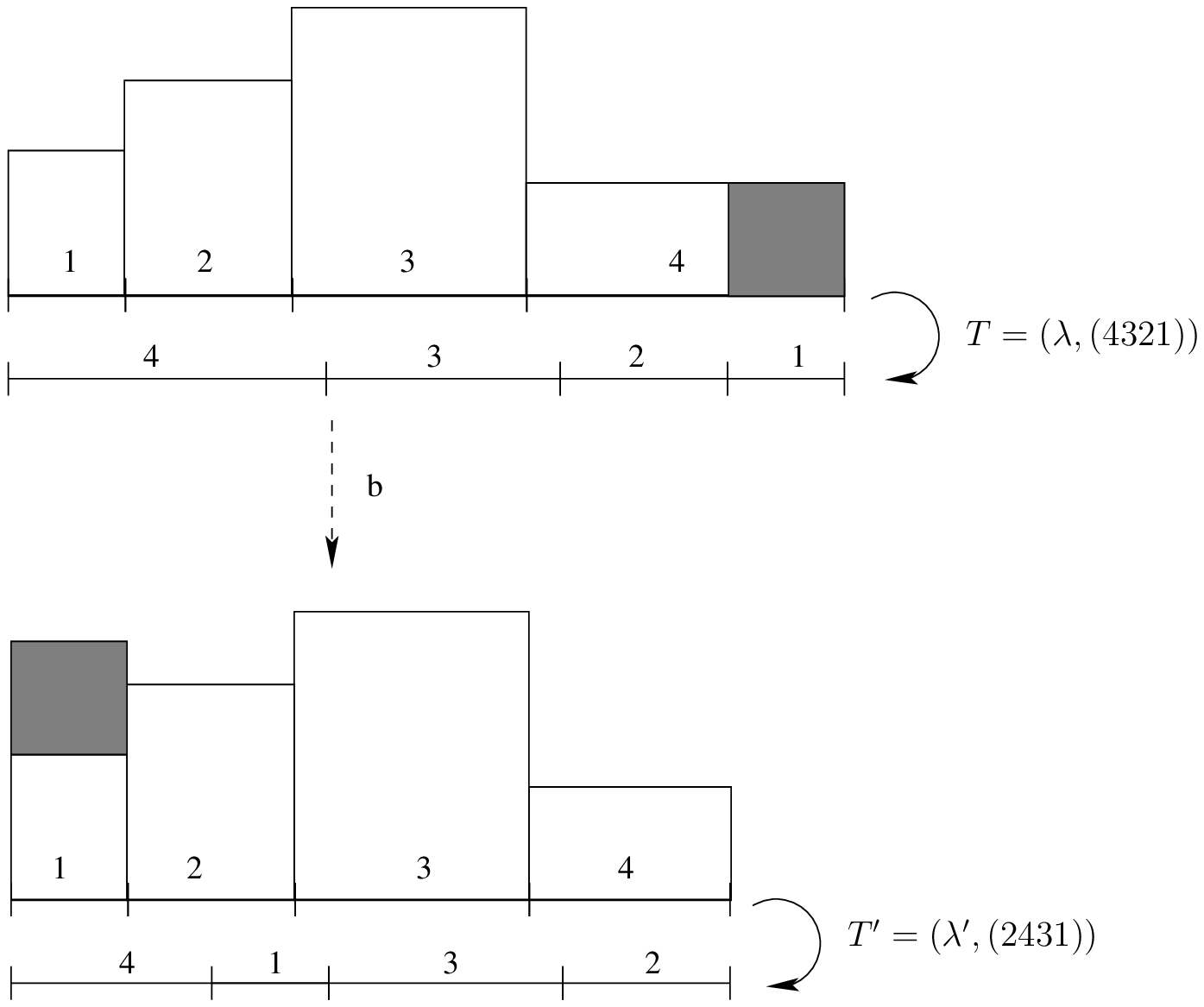}}}
\caption{ Stacking for $T=(\lambda, (4321))$.}
\label{stacking}
\end{figure}

\paragraph{Algorithm action on towers.}
The Rauzy-Veech algorithm can be visualized as acting on the towers, in terms of \emph{stacking} towers. One step corresponds to cutting the last tower before the permutation, i.e. $Z^{(r)}_{d}$, and stacking it over $Z^{(r)}_{\pi^{-1}d}$. In the case $a$, when ${\lambda}^{(r)}_{d} < \lambda^{(r)}_{\pi^{-1}d}$, $Z^{(r)}_{d}$ is completely cut and stacked above $Z^{(r)}_{\pi^{-1}d}$, at its right end (see Figure \ref{stacka}). In the case $b$, $\lambda^{(r)}_{d} > \lambda^{(r)}_{\pi^{-1}d}$, only the right portion of $Z^{(r)}_{d}$ of width $\lambda^{(r)}_{\pi^{-1}d}$ is cut and stacked completely above $Z^{(r)}_{\pi^{-1}d}$ (see Figure \ref{stackb}).
 
It is clear from the stacking description of the algorithm that each tower $Z^{(r)}_{j_{0}}$ consists of pieces of towers $Z^{(s)}_{j}$. 

\paragraph{Towers partitions.}
Define the following system of measurable partitions $\xi_s =\xi_s ( Z^{(r)}_{j_{0}} )$ of the tower $Z^{(r)}_{j_{0}}$ in terms of the subtowers $Z^{(s)}_{i}$, $0\leq s \leq r$.
The elements of the partition are complete blocks of floors of $Z^{(r)}_{j_{0}}$ which are all contained inside the same tower $Z^{(s)}_{j}$: namely, for each floor $Z^{(r)}_{j_{0},l}$ which is contained in $I^{(s)}$, construct an element $Z \in \xi_s$ in the following way. If $Z^{(r)}_{j_{0},l} \subset I^{(s)}_{j}$,  
\begin{equation}
Z \doteqdot  \bigcup_{i=0}^{h^{(s)}_{j}-1} T^i Z^{(r)}_{j_{0},l}.
\end{equation}
The set of all such $Z$ gives a partition $\xi_s$ of $Z^{(r)}_{j_{0}}$. Clearly for each  $Z \in \xi_s$ there is a unique $j$ such that $Z \subset Z^{(s)}_{j}$.
Partitions $\xi_{s'}$, $s'<s$, are refinements of $\xi_s$. 



The entries of $A^{(m,n)}$ have the following meaning for the partition $\xi_m ( Z^{(n)}_{j})$.
For $m<n$, $\ZLp{m,n}_{ij}$ gives the number of visits of $x \in I^{(n)}_j$ to $I^{(m)}_i$ under the action of $T^{(m)}$ before the first return to $I^{(m)}$. Hence 
\be \label{entriesBmn}
 A^{(m,n)}_{ij}= \# \{ Z \in \xi_m ( Z^{(n)}_{j}) |  \quad Z \subset Z^{(m)}_{i} \}.
\ee

\subsection{Notation Summary.}
For the convenience of the reader, we list in this section the notation relative to the induction algorithm for IETs introduced in the last paragraphs and used in the following ones.
\\
\begin{tabular}{ll}
$\R$ & Rauzy induction, $\R$: IETs $\rightarrow$ IETs; \\
$A^{-1}$ &  Rauzy lengths cocycle, $A^{-1}$: $\Delta(\mathscr{R}) \rightarrow SL(d,\mathbb{Z})$; \\
$\RL{n}$ & $ = \RLs(\R^n T)$ $n^{th}$ elementary Rauzy matrix; \\
$\RLp{n,m}$ & $=\RL{n} \dots \RL{m-1} $ product of Rauzy matrices; \\
$\RLp{n}$ & $=\RLp{0,n}$; 
\end{tabular}\\
\begin{tabular}{ll}
$I^{(n)}$ & $n^{th}$ inducing subinterval; \\
$T^{(n)}$ & IET obtained as first return on $I^{(n)}$;  \\
$I^{(n)}_j$ & $j=1, \dots , d$ subintervals exchanged by $T^{(n)}$;\\
$\lambda_j^{(n)}$ & = $|I^{(n)}_j|$ length of $I^{(n)}_j$, w.r.t. the Lebesgue measure on $I$;\\
$\underline{\lambda}^{(n)}$ &= $(\lambda_1^{(n)}, \dots, \lambda_j^{(n)}, \dots, \lambda_d^{(n)})$ lengths vector; \\
${\lambda}^{(n)}$ & = $| \underline{\lambda}^{(n)} | = \sum_{j=1}^d \lambda_j^{(n)}$ lenght $|I^{(n)}|$ of the $n^th$ inducing subinterval;
\end{tabular}\\
\begin{tabular}{ll}
$\RLp{ n ,m }_{ij }$ & number of visits of $x_0 \in I^{(m)}_{j} $ to $I^{(n)}_i$ before the first return to $ I^{(m)}$;\\
$h^{(n)}_j$ & $= \sum_i \RLp{ n ,m }_{ij }$, $j=1, \dots , d$ first return time of any $x\in I^{(n)}_j $ to $I^{(n)}$;\\
$Z^{(n)}_j$ & $= \cup_{k=1 }^{h^{(n)}_j-1} T^k  I^{(n)}_j$, for $j=1, \dots , d$ towers representing $T$ over $T^{(n)}$;\\
${h}^{(n)}$ & $= \max_{j=1,\dots, d} h^{(n)}_j$ maximum of the towers heights;\\
$\xi_{n}$ & partition into floors of step $n$, i.e $T^k(I^{(n)}_j)$, $k\leq h^{(n)}_j$, $j=1,\dots,d$;\\
$\xi_{n} (Z) $ & restriction of the partition $\xi_{n}$ to the set $Z$;
\end{tabular}\\
\begin{tabular}{ll}
$\|A\|$ & where $A\in SL(d,\mathbb{Z})$, $A\geq 0$, is given by $\|A\|=\sum_{ij} A_{ij} ; $\\
$A> 0$ & (or $A \geq0) $ if $A\in SL(d,\mathbb{Z})$ has positive (non-negative) entries; \\
$d_H$ & Hilbert distance on $\Delta_{d-1}$; \\
$diam_H$ & diameter w.r.t. the Hilbert distance; \\
$\widetilde{A}$ & $= A\underline{\lambda}/ |A\underline{\lambda}|$ projective transformation associated to $A\geq 0$;\\
$D(A)$&$= diam_h (\widetilde{A} (\Delta_{d-1}))$.
\end{tabular}\\


\section{Growth of Birkhoff sums of derivatives
.}\label{growthBSsec}
Let us introduce two auxiliary functions $u$, $v$ defined on $I^{(0)}$:
\bes
u(x) \doteqdot \frac{1}{ x};\qquad v(x) \doteqdot \frac{1}{1- x}.
\ees

\begin{prop}\label{reducetouv}
Assume $T$ is uniquely ergodic. There exists a sequence $\alpha_r$ such that $\alpha_r \rightarrow 0$ as $r \rightarrow \infty$ and for all $x$ distinct from singularities of $\BS{f}{r}$,
\bes
\BS{f'}{r}(x) = (  - C^+ + \alpha_r^+ )  \BS{u}{r}(x) + (  C^- + \alpha_r^- ) \BS{v}{r}(x), 
\ees
where $|\alpha_r^{\pm}| \leq \alpha_r$.
\end{prop}
\begin{proof}
See Theorem 3.1 in \cite{Ko:nonI}. The same proof applies also for uniquely ergodic IETs.
\end{proof}

In Section \ref{growthBSusec} we prove estimates on the growth of the Birkhoff sums for $u$ and $v$ for a typical IET and then we use them in Section \ref{BSf'andf''sec} to derive some information about the growth of $\BS{f'}{r}$ and $\BS{f''}{r}$. It is enough to get estimates from $u$, since estimates from $v$ can be easily derived from the following observation.
Let $\I(x)=1-x$, 
 be the reflection on the interval $I^{(0)}$. Since $v(x)=u(\I x)$, 
$v\cdot T^n   = u\cdot (\I \cdot T\cdot \I ^{-1})^n \cdot \I  $. 
Let us denote by $T^{\I } \doteqdot \I \cdot T \cdot \I ^{-1}$. Hence the Birkhoff sums for $v$ with respect to $T$ and the ones for $u$ with respect to $T^{\I}$ are related by
\be\label{BSuvrel}
\BS{v,T}{r}(x)= \BS{u,T^{\I}}{r}(1-x).
\ee 
Remark that if $T=((\lambda_1,\lambda_2, \dots, \lambda_n),\pi)$, then $T^{\I}= ((\lambda_n,\lambda_{n-1}, \dots, \lambda_1),\pi^{\I})$ where $\pi^{\I}\doteqdot (n \, n\!-\!1\dots 2\, 1)\cdot \pi \cdot (n\, n\!-\!1\dots 2\, 1)$. Hence the map $T\mapsto T^{\I}$ from $\Delta_{d-1} \times \mathscr{R}(\pi) \rightarrow \Delta_{d-1} \times \mathscr{R}(\pi^{\I})$ preserves the Lebesgue measure.




\subsection{A diophantine-type condition for IETs.} \label{existencebalancedtimessec}
In this section we define the set of full measure of IETs for which we prove Theorem \ref{mixing}.
Proposition \ref{existencebalancedtimes} shows that for typical $T$ one can find a subsequence $\{ n_l\}_{l\in \mathbb{N}}$ of induction times such that the corresponding IETs $\{ \R^{n_l} T \}_{l\in \mathbb{N}}$ in the Rauzy orbit $\{ \R ^n T \}_{n\in \mathbb{N}}$ enjoy some good properties (listed in the Proposition \ref{existencebalancedtimes} below), which we call \emph{balance}; moreover it gives a control on their frequencies.  
In the next Sections we will use these balanced induction times 
in order to estimate the growth of $\BS{u}{r}$.   

Balanced times are related to occurrences of some positive matrices in the renormalization cocycle. Conditions on the frequencies of occurrence of such matrices play for interval exchanges a role analogous to diophantine conditions for rotations. Different type of estimates in this spirit appear in the works of \cite{Ke:sim, MMY:coh, Bu:dec, AGY:exp}. The condition that we use is derived from \cite{AGY:exp}.



\subsubsection{Existence of balanced return times.}
\begin{prop}\label{existencebalancedtimes}
Let $1<\tau <2$. For each $\pi$ and for Lebesgue a.e. $\lambda\in  \Delta_{d-1}$, there exist a subsequence $\{ n_l\}_{l\in \mathbb{N}}$ of induction times, $\nu>1$, $\kappa >1$, $0<D < \infty$ and $\overline{l} \in \mathbb{N}$, such that the following properties hold for all $l\in \mathbb{N}$:
\begin{enumerate}
\item \emph{( $\nu$ -balance of  lengths )}
\be \label{lengthsbalance} 
\frac{1}{\nu} \leq \frac{\lambda_i^{(n_l)}}{\lambda_j^{(n_l)}} \leq \nu, \quad  \forall \, 1 \leq i, \, j \leq d;  
\ee
\item \emph{($\kappa$ -balance of  heights)}
\be \label{heightsbalance} 
\frac{1}{\kappa} \leq \frac{h_i^{(n_l)}}{h_j^{(n_l)}} \leq \kappa, \quad  \forall \, 1 \leq i, \, j \leq d;  
\ee
\item \emph{(positivity)} 
\be \label{positive} 
\ZLp{n_{l\phantom{\overline{l}}},n_{l+\overline{l}} }> 0, \quad \mathrm{and}\footnote{Recall that $D(A)$ was defined in (\ref{D_Adef}).}\quad D({\ZLp{n_{l\phantom{\overline{l}}},n_{l+\overline{l}}}}) \leq D .  
\ee
\item \emph{(integrability)}
\be \label{integrability} 
\lim _{l \rightarrow + \infty} \frac {\|\ZLp{n_l,n_{l+1}} \|}{l^{\tau}} =0 .
\ee
\end{enumerate}
\end{prop}

A return time which satisfies Properties $1$ and $2$  will be called \emph{balanced return time}. 
A balanced return time is such that 
the lengths and the heights of the induction towers
are approximately of the same size. 
Property $3$ gives some uniform distribution of subintervals of time $n_{l+\overline{l}}$ inside the subintervals of the previous balanced time $n_l$.
Property $4$ 
 is the diophantine condition which guarantees some control of the frequencies of occurrence of balanced times. It will be deduced from the power integrability of a certain induced cocycle, proved in \cite{AGY:exp}. 

We remark that (\ref{integrability}) is analogous to the diophantine condition used for rotations in \cite{SK:mix}, i.e. $k_l = o( l^{\tau})$, where $\{ k_l \}_{l \in \mathbb{N}}$ are the entries of the continued fraction and the exponent $\tau$ satisfies the same assumption $1 < \tau <2$.  

\begin{defn}\label{Mdef}
Let $\M^+ = \M^+ (\Delta_{d-1}\times \mathscr{R}(\pi))$ be the set of IETs in $\Delta_{d-1}\times \mathscr{R}(\pi)$ such that Proposition \ref{existencebalancedtimes} hold and $\M^- = \M^- ( \Delta_{d-1} \times \mathscr{R}(\pi) )$ be the set of $T\in \Delta_{d-1} \times \mathscr{R}(\pi)$ such that $T^{\I} \in \M^+ (\Delta_{d-1} \times \mathscr{R}(\pi^{\I}))$.
Denote $\M = \M^+ \cap \M^-.$
\end{defn}
The IETs in $\M$ are the ones for which we prove mixing of the suspension flows having one asymmetric logarithmic singularity.
\begin{rem} The set $\M$ has full measure. 
Indeed, $\M^+$ has full measure by Proposition \ref{existencebalancedtimes} and also $\M^-$ has full measure since, as already remarked, $T \mapsto T^{\I}$ preserves the Lebesgue measure.
\end{rem}
\begin{rem}
The IETs in $\M$ are uniquely ergodic, as it follows from Property (\ref{positive}) in Proposition \ref{existencebalancedtimes} with the help of techniques used by Veech in \cite{Ve:gau, Ve:pro}.
\end{rem}



In the remaining part of this section we derive Proposition \ref{existencebalancedtimes} from a result in \cite{AGY:exp} and then prove some simple corollaries, which will be used later.

\subsubsection{Proof of Proposition \ref{existencebalancedtimes}.}
If $Y\subset \hat{\Upsilon}^{(1)}_{\mathscr{\R}}$, let $\RLs _{Y}$ denote the induced cocycle of the Rauzy-Veech lengths cocycle associated to first returns to $Y$ under $\hat{\R}$, i.e. for $(\underline{\lambda}, \pi, \underline{\tau} )\in {Y}$,
\bes
\RLs_{Y}(\underline{\lambda}, \pi, \underline{\tau}  ) \doteqdot \RLp{r_{Y}} \left( (\underline{\lambda}, \pi ) \right), \quad \mathrm{where} \quad r_{Y}\doteqdot \min \{ r\in \mathbb{N}^+| \quad \hat{\R}^r (\underline{\lambda}, \pi, \underline{\tau}  ) \in Y \}.
\ees

The following result is proved in \cite{AGY:exp}.
\begin{thm}[Avila, Gou\"{e}zel, Yoccoz]\label{AGY}
For every $\delta>0$ there exists a finite union
\bes
\hat{Z}^{(1)}
 \doteqdot \left( \bigcup_{i=1}^{n} \Delta(\gamma_{s_i}) \times \{ \pi_i \} \times \Theta(\gamma_{e_i})\right) \cap \hat{\Upsilon}^{(1)}_{\mathscr{R}}, 
\ees
where $\pi_i$ is both the initial permutation of the path $\gamma_{s_i}$ and the final one of the path $\gamma_{e_i}$ and where $A(\gamma_{s_i})>0$ and $A(\gamma_{e_i})>0$ for all $i=1, \dots, n$, such that
\be\label{integrabilityAGY}
\int_{\hat{Z}^{(1)}
} \| \RLs_{\hat{Z}^{(1)}
}  \|^{1-\delta} \ud \hat{m} < \infty.
\ee
\end{thm}
\noindent Theorem \ref{AGY} is a reformulation of Theorem $4.10$ in \cite{AGY:exp}. The original statement claims the integrability of $e^{(1-\delta)r_{\hat{Z}^{(1)}}}$, where $r_{\hat{Z}^{(1)}}$ is the first return time of $(\underline{\lambda} , \pi, \underline{\tau}) \in \hat{Z}^{(1)}$ under the Veech flow, which is given by
\bes
r_{\hat{Z}^{(1)}}( (\underline{\lambda} , \pi, \underline{\tau})) \doteqdot  - \log | \RLs_{\hat{Z}^{(1)}
}^{-1} \, \underline{\lambda} \, | =  
  \log | \RLs_{\hat{Z}^{(1)}
}\, \underline{\lambda }'  \, |, 
\ees
where $(\underline{\lambda}',\pi')= R^{r_{\hat{Z}^{(1)}}} (\underline{\lambda} , \pi)$. The second equality follows by taking norms of $\RLs_{\hat{Z}^{(1)}}\, \underline{\lambda }' = \underline{\lambda } /   | \RLs_{\hat{Z}^{(1)}}^{-1} \, \underline{\lambda} \, |$. Since $\underline{\lambda }'$ belong to the compact set $\cup_{i=1}^{n} \Delta(\gamma_{s_i}) $, 
\bes
 \log | \RLs_{\hat{Z}^{(1)}
} \underline{\lambda }'  | \geq   \log ( \min_i\lambda'_i \,  \| \RLs_{\hat{Z}^{(1)}
 } \| ) \geq const + \log \| \RLs_{\hat{Z}^{(1)}
}\| .
\ees 
Hence (\ref{integrabilityAGY}) follows from the integrability of $e^{(1-\delta)r_{\hat{Z}^{(1)}}}$.
Positivity of  $A(\gamma_{s_i})$ and $A(\gamma_{e_i})$ is clear from the proof of Thm. 4.10, in which $\gamma_{s_i}$ and  $\gamma_{e_i}$ are chosen minimal and $2d-3$ complete and hence positive by Lemma $3.3$.  

\begin{proofof}{Proposition}{existencebalancedtimes}
Given $1<\tau<2$, let $\delta\doteqdot 1-\tau^{-1}>0$. Let $\hat{Z}^{(1)}$ be the corresponding set given by Theorem \ref{AGY}. Let $\overline{l}$ be the maximum length of the paths $\gamma_{s_i}$ and $\gamma_{e_i}$ for $i=1,\dots, n$.

Given $(\underline{\lambda},\pi)$, choose any $\underline{\tau} \in \Theta_{\pi}$. Let $\{n_l \}_l\in \mathbb{N}$ be the subsequence of visits of the $\hat{\R}$ orbit of  $(\underline{\lambda},\pi , \underline{\tau})$ to $\hat{Z}^{(1)}$ given by
\begin{eqnarray}\label{sequencedef}
n_0 &\doteqdot &\min \{ n\in \mathbb{N}^+ |\, n\geq \overline{l}, \,  \hat{\R}^n (\underline{\lambda},\pi , \underline{\tau}) \in \hat{Z}^{(1)} \}; \\
n_{l+1} &\doteqdot &\min \{ n \in \mathbb{N}^+ |\, n> n_l , \, \hat{\R}^n (\underline{\lambda},\pi , \underline{\tau}) \in \hat{Z}^{(1)} \}; \qquad l \in \mathbb{N^+}. 
\end{eqnarray}
Notice that (\ref{sequencedef}) is independent on $\underline{\tau}$, since as soon as $n$ is bigger than the maximum length $\overline{l}$ of the paths $\gamma_{e_i}$, visits to $\hat{Z}^{(1)}$ are determined by $\underline{\lambda}$ only (see Remark \ref{sequence} and \ref{rempast}).
Let us show that Properties $1$, $2$ (balance) and $3$ (positivity) of Proposition \ref{existencebalancedtimes} automatically hold for $(\underline{\lambda}, \pi)$ and the sequence  $\{n_l \}_{l\in \mathbb{N}}$. 

Since  by definition $\hat{\R}^{n_l} (\underline{\lambda},\pi, \underline{\tau}) \in \Delta(\gamma_{s_j})\times \{ \pi_j \} \times \Theta ( {\gamma_{s_j}})  $ for some $j$, in particular $\lambda^{(n_l)} / | \lambda^{(n_l)}| \in  \Delta(\gamma_{s_j}) $. By positivity of the $A( \gamma_{s_i} ) >0$, the union $\cup_i \Delta(\gamma_{s_i}) $ is compact and hence (see subsection \ref{Hilbertsec}) contained in a ball for the Hilbert metric $d_H$, centered at $(1/d, \dots, 1/d)$, of some radius $r_{s}>0$. Hence, 
\bes
d_H \left( \frac{\underline{\lambda}^{(n_l)} }{ | \underline{\lambda}^{(n_l)}|}  ,\left( \frac{1}{d}, \dots, \frac{1}{d}  \right) \right) \leq r_s \quad \mathrm{or}\, \mathrm{equivalently} \quad \frac{  \max_i \lambda^{(n_l)}_i }{ \min_i \lambda^{(n_l)}_i } \leq e^{r_s},
\ees
which, setting $\nu\doteqdot e^{r_s}> 1$,  is $\nu$-balance of lengths.

Similarly, from $\underline{\tau}^{(n_l)} \in \Theta ( {\gamma_{e_j}})$ we get by Remark \ref{rempast} that $\underline{\lambda}^{(n_l)} = A(\gamma_{e_j})^{-1} \underline{\lambda}^{(n_l-L)}$, where $L$ is length of $\gamma_{e_j}$. Since the heights transform according to (\ref{hZtrasf}), $\underline{h}^{(n_l)} = A(\gamma_{e_j})^T \underline{h}^{(n_l-L)}$. Arguing as above, by compactness, the union $\cup_i \widetilde{A(\gamma_{e_i})^T} \Delta_{d-1}$ is contained in a ball centered at $(1/d, \dots, 1/d)$ of some radius $r_{e}>0$ and this gives $\kappa\doteqdot e^{r_e}$ balance of the heights.

For Property $3$, 
since $\overline{l}$ is the maximum lengths of the paths $\gamma_{s_i}$ and $n_{l+\overline{l}}\geq n_{l}+\overline{l}$, by Remark \ref{sequence} we have $\RLp{n_l, n_{l+\overline{l}}} = A(\gamma_{s_j}) A $ for some $A\geq 0$. Hence $\RLp{n_l, n_{l+\overline{l}}}>0$ and $D(\RLp{n_l, n_{l+\overline{l}}}) \leq D\left( A(\gamma_{s_j}) \right) \leq 2 r_s$.

Let us show first that Property $4$ holds for typical $(\underline{\lambda},\pi,\underline{\tau})\in \hat{Z}^{(1)}$. Remark that $ \RLp{n_l, n_{l+1}} ( \underline{\lambda}, \pi, \underline{\tau } ) = \RLs_{\hat{Z}^{(1)}} (\hat{R}^{n_l} (\underline{\lambda},\pi, \underline{\tau}) ) $.
For each $\epsilon_k\doteqdot 1/k$, by $\hat{\R}$ invariance of $\hat{m}$, 
\bes
\begin{split}
\hat{m} \{ ( \underline{\lambda}, \pi, \underline{\tau } ) \in \hat{R}^{n_l} | &\,\quad \| \RLp{n_l, n_{l+1}}( \underline{\lambda}, \pi, \underline{\tau } )  \|  \geq \epsilon_k l^{\tau} \}
= \\ &=  \hat{m} \{  ( \underline{\lambda}, \pi, \underline{\tau } ) \in \hat{R}^{n_l} | \quad \|  \RLs_{\hat{Z}^{(1)}} ( \underline{\lambda}, \pi, \underline{\tau } )   \| ^{\tau^{-1}} \epsilon_k ^{-\tau^{-1}} \geq l  \} .
\end{split}
\ees
Since we chose $\tau^{-1}=1-\delta$, the integrability condition (\ref{integrabilityAGY}) implies that $\sum_l  \hat{m} \{ \|  \RLs_{\hat{Z}^{(1)}}  \| ^{\tau^{-1}}  \epsilon_k ^{-\tau^{-1}} \geq l  \} < \infty$ for each $\epsilon_k$. Hence, it follows by Borel-Cantelli that there exists a subset $\hat{Z}'\subset \hat{Z}^{(1)}$ with $\hat{m}(\hat{Z}')= \hat{m} (\hat{Z}^{(1)})$ such that for $(\underline{\lambda}, \pi, \underline{\tau })\in \hat{Z}'$ the sequence $n_l$ satisfies (\ref{integrability}). 

Consider the projection $p\hat{Z}'\subset \Delta(\mathscr{R})$.
By indipendence of the definition of $\{n_l\}$ on $\tau$ , for each $(\underline{\lambda},\pi ) \in p\hat{Z}'$ we have (\ref{integrability}). Moreover, $p\hat{Z}'$ has $p\hat{m}$ full measure in $p \hat{Z}^{(1)}= \cup_i \Delta(\gamma_{s_i})\times \{ \pi_i \}$ and in particular positive $p\hat{m}=\mu_{\V}$ measure and hence also $\mu_{Z}$ positive measure. 

To conclude, let $\M^+$ be the set of $T\in \Delta(\mathscr{R})$ such that there exists $\overline{n}$, for which $R^{\overline{n}}T\in p\hat{Z}'$. Clearly if $T\in \M^+$, all the Properties are satisfied by the sequence $n_l \doteqdot \widetilde{n}_{l}+\overline{n}$, where $\widetilde{n}_l$ is the sequence associated to $\R^{\overline{n}}T$. To see that $\M^+$ has full measure, it is enough to use ergodicity of $\Z$ and the fact that $\mu_{\Z} (p\hat{Z}' )>0$, remarking  that $\Z$ orbits are subsets of $\R$ orbits. The formulation in Proposition \ref{existencebalancedtimes} follows by absolute continuity of $\mu_{\Z}$ w.r.t. Lebesgue. 

\end{proofof}

\subsubsection{Some consequences of balance.}
We will also frequently use the following simple lemmas. Recall that $A>0$ means that $A$ has strictly positive entries.
\begin{lemma}\label{productpositive}
Let $A_i >0$, $A_i \in SL(d, \mathbb{Z})$ for $i=0, \dots, n$. If 
$\underline{\lambda} = A_0 \cdot \underline{\lambda}'$, then $\sum_j\lambda_j > d \sum_j \lambda_j'$. 

If $\underline{h} = A_0 \cdot \underline{h}'$, then $\min_j h_j > d \min_j h_j'$. In particular if $\underline{h} = A_1 \cdots A_n \, \underline{e}$, then $\min_j {h_j} \geq d^n $.
\end{lemma}
\begin{proof}
All the properties follow easily from $A_{ij} \geq 1$.
\end{proof}

Recall that we defined $ \lambda^{(n_l)} \doteqdot \sum_j  \lambda^{(n_l)}_j $ and $ h^{(n_l)} \doteqdot \max_j  h^{(n_l)}_j$.

\begin{cor} \label{superexp} For each $L\in \mathbb{N}$, $\log(h^{(n_{L \overline{l}})}) \geq {L} \log d$.
In particular,
\be \label{superexplimit}
\lim_{l\rightarrow +\infty} \frac{1}{\log h^{(n_{l})} }=0.
\ee
\end{cor}

\begin{proof}
For the first property, apply Lemma \ref{productpositive} to $h^{(n_{L \overline{l}})}= {\ZLp{n_{L \overline{l}}}}^T \underline{e}$ and remark that ${\ZLp{n_{L \overline{l}}}}^T $ is the product of at least $L$ positive matrices by Property 3 (see (\ref{positive})) of the sequence $\{n_l\}_{l\in \mathbb{N}}$. It follows that $\log h^{(n_{l})} \geq \log d  [{l}/{\overline{l}}]  
$
 and hence we get (\ref{superexplimit}).
\end{proof}

Recall that we defined $ \lambda^{(n_l)} \doteqdot \sum_j  \lambda^{(n_l)}_j $ and $ h^{(n_l)} \doteqdot \max_j  h^{(n_l)}_j$.
\begin{lemma}\label{balancehlambda}
If $n_l$ has $\nu$-balanced lengths (\ref{lengthsbalance}) and $\kappa$-balanced heights (\ref{heightsbalance}), then for each $j=0,\dots, d$,
\begin{eqnarray}
\frac{1}{\kappa \lambda^{(n_l)} } &\leq &h^{(n_l)}_j \leq \frac{\kappa}{ \lambda^{(n_l)} } ;\label{hcomparedtolambda} \\ 
\frac{1}{d \kappa \nu h^{(n_l)} } &\leq & \lambda^{(n_l)}_j \leq \frac{\kappa}{h^{(n_l)} } .\label{lambdacomparedtoh} 
\end{eqnarray}
\end{lemma}
Hence, if $n_l$ is a balanced time, $\lambda^{(n_l)} \approx 1/ h^{(n_l)}$ up to constants.
\begin{proof}
Since by Kac lemma $\sum_i h^{(n_l)}_i \lambda^{(n_l)}_i =1$, remark that $\min_i h^{(n_l)}_i \lambda^{(n_l)} \leq 1 $ and $\max_i h^{(n_l)}_i \lambda^{(n_l)} \geq 1 $. Together with heights balance (\ref{heightsbalance}), this gives
\bes
\frac{1}{\kappa \lambda^{(n_l)} } \leq \frac{1}{\kappa}  \max_i h^{(n_l)}_i  \leq h^{(n_l)}_j \leq {\kappa} \min_i h^{(n_l)}_i \leq \frac{\kappa}{ \lambda^{(n_l)} }.
\ees 
To show (\ref{lambdacomparedtoh}), let $i$ be such that $ h^{(n_l)}_i \lambda^{(n_l)}_i = \max_j h^{(n_l)}_j \lambda^{(n_l)}_j \geq  1/d$. Also, for each $j$, $h^{(n_l)}_j \lambda^{(n_l)}_j <1$. Using also lengths balance (\ref{lengthsbalance}),
\bes
\frac{1}{\kappa \nu d h^{(n_l)}_{\phantom{j}} } \leq  \frac{1}{\nu d h^{(n_l)}_i }  \leq  \frac{\lambda^{(n_l)}_i}{\nu} \leq \lambda^{(n_l)}_j \leq  \frac{1 }{h^{(n_l)}_j} \leq \frac{\kappa}{ h^{(n_l)}_{\phantom{j}} }.
\ees 
\end{proof}
\begin{lemma} \label{loglimit}
For each fixed $L\in \mathbb{N}$, 
\be
\lim_{l\rightarrow + \infty } \frac{ \log \|\ZLp{n_l, n_{l+L}} \|}{\log h^{(n_l)}} = \lim_{l\rightarrow + \infty } \frac{ \log \|\ZLp{n_{l-L}, n_l} \|}{\log h^{(n_l)}} = 0.
\ee
\end{lemma}
\begin{proof}
By Corollary \ref{superexp}, $\log h^{(n_{l })} \geq {[l/\overline{l}]} \log d$, where $[\cdot]$ denotes the integer part.
Hence
\be
\lim_{l\rightarrow + \infty } \frac{ \log \|\ZLp{n_l, n_{l+L}} \|}{\log h^{(n_l)} } \leq  \lim_{l\rightarrow + \infty }
\frac{\sum_{i=l}^{l + L-1} \log \left\| \ZLp{n_{i},n_{i+1} } \right\|} {(l/\overline{l} -1 )\log d  } .
\ee
Using the property (\ref{integrability}), each of the $L$ terms in the sum in the RHS can be bounded for $l\gg1 $ by $\tau \log i \leq \tau \log (l+L-1)$ and hence the first limit is zero. The second one is analogous.
\end{proof}

\subsection{Growth of Birkhoff sums of $u$.}\label{growthBSusec}
Assume in this section that the roof function is $u(x)=1/x$. 
\subsubsection{Growth of Birkhoff sums along a balanced tower.}\label{BSalongtowersec}
In order to understand the asymptotic growth of $\BS{u}{r}(x)$, we first consider $\BS{u}{r}(x)$ when $x \in I^{(n)}_j$ is a point in the base of the tower $Z^{(n)}_j$, $r=h^{(n)}_j$ is exactly the height of the same tower and $n$ is one of the balanced return times constructed in Section \ref{existencebalancedtimes}. This preliminary estimate is used in Section \ref{generalx} as a building block to get an estimate for any $r$ and most of the other points $x$.    

\begin{prop}\label{alongtower}
Assume $T\in \M^+$. Let $n_{l_0}$ be a balanced return time given by Proposition \ref{existencebalancedtimes}. Let $x_0 \in I^{(n_{l_0})}_{j_0}$ be a point belonging to the base of the tower $Z^{(n_{l_0})}_{j_0}$ and let $r_0= h^{(n_{l_0})}_{j_0}$ be the corresponding tower height. Given $\varepsilon>0$, there exists $l(\varepsilon)$ such that for $l_0 \geq l(\varepsilon)$,
\be \label{alongtowerestimate}
(1-\varepsilon)  h^{(n_{l_0})}_{j_0} \ln 
\left( 
h^{(n_{l_0})} 
\right) 
\leq \BS{u}{{r_0}}(x_0) - \frac{1}{x_0}
\leq  (1+ \varepsilon) h^{(n_{l_0})}_{j_0} \ln \left( h^{(n_{l_0})} \right).
\ee
\end{prop}

A Birkhoff sum of the form $\BS{u}{{r_0}}(x_0)$ where $x_0$ and $r_0$ satisfy the hypotheses of Proposition \ref{alongtower} will be referred as \emph{Birkhoff sum along a tower}. The proposition shows that each sum along a tower gives a contribution of order $r_0\log(r_0)$, plus the contribution of the closest point to the singularity, $x_0$, which could be arbitrary large and will be estimated separately when using this sums as building blocks in Section \ref{generalx}.


\begin{proofof}{Proposition}{alongtower}
Consider the inducing intervals $I^{(n_{l})} = [0,\lambda^{(n_{l})}]$ 
where $\{n_l\}_{l\in \mathbb{N}}$ is the sequence of balanced induction times constructed in Section \ref{existencebalancedtimessec}.
Given any $\epsilon > 0$, let $D$ and $\overline{l}$ be given by Proposition \ref{existencebalancedtimes}. Choose $L_1 \in \mathbb{N}$ such that 
\be \label{L_1def}
(1-e^{-D})^{L_1-1} D  
 <\epsilon .
\ee
We remark that $diam(\widetilde {\ZLp{n_{\overline{l}}}} \Delta_{d-1}) < \infty$ by (\ref{positive}). Choose also $ L_2 \in \mathbb{N}$ such that $1/d^{L_2} <\epsilon$. Assume $l_0 \geq \overline{l}(1 + L_1 + L_2)$. For convenience, introduce the notation:
\be \label{nl1nl2def}
{l_{-1}} \doteqdot {l_0 -L_1 \overline{l}};\qquad {l_{-2}} \doteqdot {l_{-1} -L_{2}\overline{l}} = {l_0  - ( L_1 + L_2) \overline{l} } ; 
\ee
The past induction times $n_{l_{-1}}$ and $n_{l_{-2}}$ will play the following role in the proof: $n_{l_{-1}}$ is such that the elements of the orbit $\{ T^r(x_0) \} _{0 \leq r < r_0}$ are uniformly distributed inside the elements of the partition $\xi_{n_{l_{-1}}}$; $n_{l_{-2}}$ is such that the main contribution to $\BS{u}{{r_0}}(x_0)$ comes from visits to $(\lambda^{(n_{l_{-2}})},1]$.

Denote the points along the $T$-orbit of $x_0$ by $x_i= T^i (x_0) $, $0\leq i < r_0$. 
Since the original interval can be partitioned as
\bes
I^{(0)} =  I^{(n_{l_0})} \cup  I^{(n_{l_{-2}})} \backslash  I^{(n_{l_0})} \cup (\lambda^{(n_{l_{-2}})},1],
\ees
and $x_i \notin I^{(n_{l_0})} $, for $1 \leq i < r_0$ because $r_0$ is by definition the first return time of $x_0$ to $I^{(n_{l_0})}$, the Birkhoff sums can be decomposed as follows:
\be \label{BSdecomp}
\BS{u}{{r_0}}(x_0) = \sum_{i=0}^{r_0 -1 } \frac{1}{x_i} =  \frac{1}{x_0} + \sum_{x_i \in I^{(n_{l_{-2}})} \backslash  I^{(n_{l_0})}} \frac{1}{x_i} +\sum_{x_i \in (\lambda^{(n_{l_{-2}})},1]} \frac{1}{x_i} .
\ee
We will refer to the first term in the RHS of (\ref{BSdecomp}) as \emph{singular error}, to the sum which appears as a second term as \emph{gap error},  while the sum which appears as a last term determines the \emph{main contribution}.

 
\paragraph{Uniform ergodic convergence.}
The following Lemma gives an uniform bound for the equidistribution of the orbits of points $x\in I^{(n_{l_0})}$ inside the elements of the past partition $\xi_{n_{l_{-1}}}$. 
Recall that $\ZLp{ n_{l_{-1}},n_{l_0} }_{i j}$ gives the number of visits of $x \in I^{(n_{l_0})}_{j} $ to $I^{(n_{l_{-1}})}_{i}$ before the time $h^{(n_{l_0})}_j$ of first return to $I^{(n_{l_0})}$.
\begin{lemma}[Uniform distribution] \label{uniformdistributionlemma}
For each $1\leq i,\, j\leq d $,
\be \label{uniformdistribution}
e^{-2\epsilon} \lambda^{(n_{l_{-1}})}_i  \leq \frac{\ZLp{ n_{l_{-1}},n_{l_0} }_{i j}}{h^{(n_{l_0})}_j} \leq e^{2 \epsilon} \lambda^{(n_{l_{-1}})}_i .
\ee
\end{lemma}
\begin{proof}
Consider the sets $\widetilde {\ZLp{ n_{l_{-1}},n } } \Delta_{d-1} \subset \Delta_{d-1}$, for $n > n_{l_{-1}}$, which form a nested sequence of compact sets. 
By the transformation formula (\ref{lengthsrelationR}) for lengths vectors $\underline{\lambda}^{(n_{l_{-1}})} = \ZLp{ n_{l_{-1}},n }  \underline{\lambda}^{(n)}$, the renormalized vector
\bes
\frac{ \underline{\lambda}^{(n_{l_{-1}})}} { \lambda^{(n_{l_{-1}})} } \in \bigcap_{n > n_{l_{-1}}} \widetilde {\ZLp{ n_{l_{-1}},n } } \Delta_{d-1} .
\ees
When $n= n_{l_0}$, since $l_0 = l_{-1} + L_1 \overline{l}$, applying $L_1$ times Property 3 (positivity) in Proposition \ref{existencebalancedtimes} through the contraction property (\ref{contraction}), we get
\be \label{contracteddiamater}
diam_H ( \widetilde {\ZLp{ n_{l_{-1}},n_{l_0} } } \Delta_{d-1} ) \leq (1-e^{-D})^{L_1-1} D 
\leq \epsilon , 
\ee
where the last inequality 
follows by the choice (\ref{L_1def}) of $L_1$.

Denote by $\underline{e}_j$ the unit vector $(\underline{e}_j)_i=\delta_{ij}$ ($\delta$ is here the Kronecker symbol). Since both the vectors $\widetilde { \ZLp{ n_{l_{-1}},n_{l_0} }} \underline{e}_j $ and $ \frac{ \underline{\lambda}^{(n_{l_{-1}})}} { \lambda^{(n_{l_{-1}})} } $ belong to the closure of ${\ZLp{ n_{l_{-1}},n_{l_0} } } { \Delta_{d-1} }$,
it follows by (\ref{contracteddiamater}), using compactness, that
\bes
d_H \left(  \frac{ \underline{\lambda}^{(n_{l_{-1}})}} {\lambda^{(n_{l_{-1}})} }  , 
\widetilde { \ZLp{ n_{l_{-1}},n_{l_0} } } \underline{e_j}  \right) = \log \frac{\max_{i=1,\dots , d}\frac{\ZLp{ n_{l_{-1}},n_{l_0} }_{ij}}{{\lambda_i^{(n_{l_{-1}})}} }} {\min_{i=1 , \dots ,d}\frac{\ZLp{ n_{l_{-1}},n_{l_0} }_{ij}}{{\lambda_i^{(n_{l_{-1}})}} }}    \leq \epsilon , 
\ees
where we also used the invariance of the distance expression by multiplication of the arguments by a scalar.
Equivalently, for each $1\leq i ,\, k \leq d$, 
\be
\label{projdistancecons}
e^{-\epsilon} \left( {\ZLp{ n_{l_{-1}},n_{l_0} }_{kj}}{{\lambda_i^{(n_{l_{-1}})}} }\right)
 \leq {\ZLp{ n_{l_{-1}},n_{l_0} }_{ij}}{{\lambda^{(n_{l_{-1}})}}_k } \leq e^{\epsilon}\left( {\ZLp{ n_{l_{-1}},n_{l_0} }_{kj}}{{\lambda_i^{(n_{l_{-1}})}} } \right)
\ee
and summing over $k$ we get
\be\label{1stratio}
e^{-\epsilon}\leq  \frac{{\ZLp{ n_{l_{-1}},n_{l_0} }_{ij}}{{\lambda^{(n_{l_{-1}})}} } }{{\sum_k \ZLp{ n_{l_{-1}},n_{l_0} }_{kj}}{{\lambda_i^{(n_{l_{-1}})}} } }  \leq e^{\epsilon} .
\ee
If we multiply (\ref{projdistancecons})
 by $h^{(n_{l_{-1}})}_i$ and then sum also over $i$, using that $\sum_i h^{(n_{l_{-1}})}_i \lambda ^{(n_{l_{-1}})}_i=1$ by Kac's lemma and that $\underline{h}^{(n_{l_0})}= {\ZLp{ n_{l_{-1}},n_{l_0} }}^T \underline{h}^{(n_{l_{-1}})}$,
\be\label{2ndratio}
e^{-\epsilon} \leq  \frac{ {h}^{(n_{l_0})}_j {\lambda^{(n_{l_{-1}})}}}{ { \sum_k \ZLp{ n_{l_{-1}},n_{l_0} }_{kj}}}  \leq e^{\epsilon} .
\ee
The combination of (\ref{1stratio}, \ref{2ndratio}) 
 gives (\ref{uniformdistribution}).
\end{proof}

\paragraph{Estimate of the Main Contribution.}\label{estimatemainsec}
The following lemma shows that the main contribution in (\ref{BSdecomp}) determines the order of the Birkhoff sum in Proposition \ref{alongtower}.
\begin{lemma}[Main contribution.]\label{estimatemainlemma}
For each $\epsilon>0$, if $l_0 > l_{m}(\epsilon)$,
\be \label{estimatemain}
e^{-2\epsilon} {(1-\epsilon)^2 }  {h_{j_0}^{(n_{l_0})} }  \log h^{(n_{l_0})}  \leq 
\sum_{x_i \in (\lambda^{(n_{l_{-2}})},1]} \frac{1}{x_i} 
\leq  
e^{2\epsilon} (1+ \epsilon)^2    {h_{j_0}^{(n_{l_0})} }  \log h^{(n_{l_0})}
\ee 
\end{lemma}
\begin{proof}
Consider the partition $\xi_{n_{l_{-1}}}\left( (\lambda^{(n_{l_{-2}})},1]\right)$ of $(\lambda^{(n_{l_{-2}})},1]$. Recall that the elements $F_{\alpha}\in \xi_{n_{l_{-1}}}$ are floors of the towers of step $n_{l_{-1}}$, i.e. for some $1\leq j_{\alpha} \leq d$, $F_{\alpha}= T^k( I^{(n_{l_{-1}})}_{j_{\alpha}})$, where $0\leq k< h^{(n_{l_{-1}})}_{j_{\alpha}}$. In particular $Leb(F_{\alpha})= \lambda^{(n_{l_{-1}})}_{j_{\alpha}}$. For each $F_{\alpha}$ choose, by the mean value theorem, a point $\bar{x}_{\alpha}$ such that
\be\label{xalphadef}
\frac{1}{\bar{x}_{\alpha}} \doteqdot \frac{1}{ \lambda^{(n_{l_{-1}})}_{j_{\alpha}}} \int_{F_{\alpha}} \frac{1}{s} \ud{s}.
\ee

\begin{lemma}\label{ratios}
If $x_i \in F_{\alpha}$ and $F_{\alpha} \subset (\lambda^{(n_{l_{-2}})},1]$,
\bes
1-\epsilon \leq \frac{1/x_i}{1/ \bar{x}_{\alpha}}\leq 1+\epsilon.
\ees
\end{lemma}
\begin{proof}
Let $F_{\alpha}=[a,b)$; then $a\leq x_i, \, \bar{x}_{\alpha} <b$. Since by assumption $b-a \leq \lambda^{(n_{l_{-1}})}$ and $a\geq \lambda^{(n_{l_{-2}})} $, 
\begin{eqnarray}
\frac{\bar{x}_{\alpha}}{ x_i} &\leq& \frac{b}{a}= \frac{a+(b-a)}{a} \leq 1 + \frac{\lambda^{(n_{l_{-1}})}}{\lambda^{(n_{l_{-2}})}};
\nonumber \\
\frac{\bar{x}_{\alpha}}{ x_i} &\geq& \frac{a}{b} = \frac{b-(b-a)}{b} \geq 1 - \frac{\lambda^{(n_{l_{-1}})}}{\lambda^{(n_{l_{-2}})}} .\nonumber
\end{eqnarray}
Let us show that $\lambda^{(n_{l_{-1}})}/ \lambda^{(n_{l_{-2}})}< \epsilon$. Since $l_{-1}=l_{-2}+{L_2}\overline{l}$,
\bes
\underline{\lambda}^{(n_{l_{-2}})} = \prod_{i=0}^{L_2-1} \ZLp{ n_{l_{-2} +i \overline{l}}, n_{l_{-2} + (i+1) \overline{l} }} \underline{\lambda}^{(n_{l_{-1}})}
\ees
and each of the matrices in the product has 
positive entries by Property 3 in Proposition \ref{existencebalancedtimes}. Hence by iterated application of Lemma \ref{productpositive} and by the choice of $L_2$, we get $\lambda^{(n_{l_{-1}})}/ \lambda^{(n_{l_{-2}})}< 1/d^{L_2} < \epsilon$.  
\end{proof}
Rearranging the main contribution in (\ref{BSdecomp}) by floors, i.e.
\bes
\sum_{x_i \in (\lambda^{(n_{l_{-2}})},1]}  \frac{1}{x_i}  = \sum_{F_{\alpha} \in  \xi_{n_{l_{-1}}} \left( (\lambda^{(n_{l_{-2}})},1]\right) } \sum_{x_i \in F_{\alpha}}  \frac{1}{x_i},
\ees
and applying Lemma \ref{ratios}, we get
\be \label{mainbetweenFalphas}
 \sum_{F_{\alpha} \subset  (\lambda^{(n_{l_{-2}})},1] } \sum_{x_i \in F_{\alpha}}  \frac{1-\epsilon}{\bar{x}_{\alpha}}  \leq \sum_{F_{\alpha} \subset  (\lambda^{(n_{l_{-2}})},1] } \sum_{x_i \in F_{\alpha}}  \frac{1}{x_i}\leq  
 \sum_{F_{\alpha} \subset  (\lambda^{(n_{l_{-2}})},1] } \sum_{x_i \in F_{\alpha}}  \frac{1+\epsilon}{\bar{x}_{\alpha}} .
\ee
Consider now 
\be \label{termmain}
\sum_{x_i \in F_{\alpha}}  \frac{1}{\bar{x}_{\alpha}} = \#\{ x_i \in F_{\alpha} \}  \frac{1}{\bar{x}_{\alpha}}.
\ee
Recall that $x_0\in I^{(n_{l_0})}_{j_0} $; if $F_{\alpha}= T^k( I^{(n_{l_{-1}})}_{j_{\alpha}})$ is a floor of the $j^{th}$ tower $Z^{(n_{l_{-1}})}_{j_{\alpha}}$, 
\be \label{cardinality}
\#\{ x_i \in F_{\alpha} \} =  \#\{ x_i \in I^{(n_{l_{-1}})}_{j_{\alpha}}  \} = \ZLp{ n_{l_{-1}},n_{l_0} }_{j_{\alpha} j_0}.
\ee
In (\ref{cardinality}) we used the dynamical meaning of $\ZLp{ n_{l_{-1}},n_{l_0} }_{j_{\alpha} j_0}$ 
together with the fact that $Z^{(n_{l_0})}_{j_0}$ is decomposed in a whole number of elements of $\xi_{n_{l_{-1}}}$ corresponding to towers of the previous step $n_{l_{-1}}$,  hence visits to a floor $T^k 
I^{(n_{l_{-1}})}_{j_{\alpha}} 
$ of the tower are in one to one correspondence with visits to its base $I^{(n_{l_{-1}})}_{j_{\alpha}}$. 



From Lemma \ref{uniformdistributionlemma} and (\ref{cardinality}),
\be
e^{-2\epsilon} {\lambda^{(n_{l_{-1}})}_{j_{\alpha}} } h^{(n_{l_0})}_{j_0} \leq \#\{ x_i \in F_{\alpha} \} \leq {e^{2\epsilon}} {\lambda^{(n_{l_{-1}})}_{j_{\alpha}} } {h^{(n_{l_0})}_{j_0} }.
\ee
Using this bound and recalling the definition (\ref{xalphadef}) of $\overline{x}_{\alpha}$, we get
\bes
e^{-2\epsilon} { h^{(n_{l_0})}_{j_0} }  \int_{F_{\alpha}} \frac{1}{s} \ud  s  \leq \#\{ x_i \in F_{\alpha} \}  \frac{1}{\bar{x}_{\alpha}} \leq {e^{2\epsilon}} {h^{(n_{l_0})}_{j_0} }  \int_{F_{\alpha}} \frac{1}{s} \ud  s.
\ees
Summing it over $F_{\alpha} \subset  (\lambda^{(n_{l_{-2}})},1]$ and using
\bes
 \sum_{F_{\alpha} \subset  (\lambda^{(n_{l_{-2}})},1] }   \int_{F_{\alpha}} \frac{1}{s} \ud  s =  \int_{(\lambda^{(n_{l_{-2}})},1]} \frac{1}{s}\ud s = \log \frac{1}{\lambda^{(n_{l_{-2}})}} , 
\ees
we get by (\ref{mainbetweenFalphas})
\be \label{mainwlambda}
{e^{-2\epsilon}}(1-\epsilon)   h^{(n_{l_0})}_{j_0}   \log \frac{1}{\lambda^{(n_{l_{-2}})}}  \leq \sum_{F_{\alpha} \subset  (\lambda^{(n_{l_{-2}})},1] } \sum_{x_i \in F_{\alpha}}  \frac{1}{x_i}\leq  
 {e^{2\epsilon}}  (1+\epsilon)  h^{(n_{l_0})}_{j_0}   \log \frac{1}{\lambda^{(n_{l_{-2}})}} .
\ee

In order to get the estimate (\ref{estimatemain}) of Lemma \ref{estimatemainlemma} from (\ref{mainwlambda}) it is enough to compare $1/\lambda^{(n_{l_{-2}})}$ with $h^{(n_{l_0})}$. Since $\lambda^{(n_{l_{-2}})}> \lambda^{(n_{l_0})} \geq 1/(\kappa h^{(n_{l_0})} )$ by $\kappa$-balance of heights (see Lemma \ref{balancehlambda}), we get 
\bes
\log (1/ \lambda^{(n_{l_{-2}})}) \leq   \log(h^{(n_{l_0})}) \left( 1+ \log{\kappa}/{ \log(h^{(n_{l_0})}) } \right)
\ees
and, if $l_0\geq l_m$ for some $l_m(\epsilon)>0$, the upper estimate in (\ref{estimatemain}) by Corollary \ref{superexp}.

For the lower bound, adding and subtracting $h ^{(n_{l_0})}_{j_0} \log  h^{(n_{l_0})}$,
\be \label{lambdavsh}
 h^{(n_{l_0})}_{j_0}   \log  \frac{1}{\lambda^{(n_{l_{-2}})}}  =   h^{(n_{l_0})}_{j_0} \log h^{(n_{l_0})}  \left( 1  -   \frac{ \log \left( h^{(n_{l_0})}_{j_0}  \lambda^{(n_{l_{-2}})}\right)}{  \log h^{(n_{l_0})}} \right).
\ee
In order to estimate the very last term in (\ref{lambdavsh}), notice that again by Lemma \ref{balancehlambda} 
and balance, $\lambda^{(n_{l_{-2}})} \leq \kappa /  h^{(n_{l_{-2}})} $. Hence, using that $h^{(n_{l_0})}/  h^{(n_{l_{-2}})} \leq \| \ZLp{n_{l_{-2}},n_{l_0} } \|$,
\be \label{logsumisepsilon}
\frac{ \log \left( h^{(n_{l_0})}  \lambda^{(n_{l_{-2}})}\right)}{  \log h^{(n_{l_0})}}  \leq 
\frac{\log \kappa + \log  \| \ZLp{n_{l_{-2}}
,n_{l_0} } \| }{  \log h^{(n_{l_0})}}.
\ee
Enlarging  $ l_{m}$ if necessary, the RHS is less than $\epsilon$ for $l_0\geq l_m$  by Lemma \ref{loglimit} 
(recall that the difference $l_{0}-l_{-2}=(L_1 + L_2) \overline{l}$ is fixed) and Corollary \ref{superexp}.

Combining (\ref{logsumisepsilon}) and (\ref{lambdavsh}) to estimate the LHS of (\ref{mainwlambda}) from below,
we get the lower estimate that complete the proof of the Lemma.
\end{proof}

\paragraph{Estimate of the Gap Error.}\label{estimategapsec}
\begin{lemma}[Gap error.]\label{estimategaplemma}
For each $\epsilon$, if $l_0 > l_{g}(\epsilon)$,
\be \label{estimategap}
0 \leq  \sum_{x_i \in I^{(n_{l_{-2}})} \backslash  I^{(n_{l_0})}}   \frac{1}{x_i}\,  \leq
 \epsilon    
\left( h^{(n_{l_0})}_{j_0}  \log h^{(n_{l_0})} \right).
\ee 
\end{lemma}
\begin{proof}
The bound below is trivial since $1/x_i >0$. Since we are considering $0\leq i < r_0=h^{(n_{l_0})}_{j_0}$, it follows from the tower construction (see Section \ref{towersec}) that the points $x_i= T^ix_0$ of the orbit of $x_0\in I^{(n_{l_0})}_{j_0}$ belong to different floors of the tower $Z^{(n_{l_0})}_{j_0}$ and that their minimum distance is bounded from below by
\bes
\min_{0\leq i,j < r_0} |x_i - x_j | \geq  
\lambda^{(n_{l_0})}_{j_0} 
\geq  
\frac{1}{d \kappa \nu h^{(n_{l_0})}_{j_0}} ,
\ees
where in the last inequality we used that $n_{l_0}$ is balanced and Lemma \ref{balancehlambda}.

Remarking also that $x_0$ is the closest point to the singularity, it follows that if we rearrange the $x_i$ in increasing order and relabel them $\tilde{x}_i$ ($\tilde{x}_{i}< \tilde{x}_{i+1}$), we have 
\bes
\tilde{x}_i \geq x_0 + 
\frac{i}{d \kappa \nu h^{(n_{l_0})}_{j_0}} 
, \qquad i=0, \dots , r_0-1.
\ees

Since the roof function $1/x$ is monotonically decreasing, the gap error can be bounded from above by
\be \label{arithmeticbound}
 \sum_{x_i \in I^{(n_{l_{-2}})} \backslash  I^{(n_{l_0})}}   \frac{1}{x_i}\,  \leq     \sum_{k=1}^{K}  \frac{1}{x_0 + k 
/ d \kappa \nu h_{j_0}^{(n_{l_0})}}
,
\ee 
where $K=\# \{ x_i \in  I^{(n_{l_{-2}})} \backslash  I^{(n_{l_0})} \} $ and $k\geq 1$ since $x_0 \in  I^{(n_{l_0})}$.

The following lemma is proved by Kochergin in \cite{Ko:nonI} as Lemma 5.1.
\begin{lemma}\label{Kocharithmetic} Let $h>0$ and $x>0$.
\bes
\sum_{k=1}^{K} \frac{1}{x+k h} = \frac{1}{h} \log \left( \frac{t_0 + K}{t_0+1} \right) + \frac{1}{h} R_K(t_0) ,
\ees
where $t_0=x/h$ and $0 <  R_K(t_0) < 1/(t_0+1)$.
\end{lemma}
Applying Lemma \ref{Kocharithmetic} to (\ref{arithmeticbound}) and using that $\log \left( \frac{t_0 + K}{t_0+1} \right)$ is decreasing in $t_0$, so it reaches its maximum $\log K$ at $t_0=0$, 
\be \label{arithmetic}
 \sum_{x_i \in I^{(n_{l_{-2}})} \backslash  I^{(n_{l_0})}}   \frac{1}{x_i}\, \leq 
 { d\kappa \nu}{h^{(n_{l_0})}_{j_0}} (\log {K}  + 1). 
\ee 
The cardinality $K$ of points $x_i \in  I^{(n_{l_{-2}})} \backslash  I^{(n_{l_0})} $ can be bounded by Remark \ref{dynmean} in terms of the cocycle matrices, by
\be \label{K}
K \leq \sum_{j=1}^{d} \ZLp{n_{l_{-2}},n_{l_0}}_{j j_0} \leq \| \ZLp{n_{l_{-2}},n_{l_0}} \| .
\ee
Hence, applying (\ref{arithmetic}) and (\ref{K}) we get
\be \label{gapratio}
\frac{ \sum_{x_i \in I^{(n_{l_{-2}})} \backslash  I^{(n_{l_0})}}   \frac{1}{x_i}  }{  { h^{(n_{l_0})}_{j_0} }  \log {h^{(n_{l_0})}}  } \leq  d \kappa \nu  \frac{ \log {K}  + 1} {\log  {h^{(n_{l_0})}} } \leq   
 d \kappa \nu  \frac{ \log {\| \ZLp{n_{l_0 -{(L_1 + L_2)\overline{l}}},n_{l_0}} \|   }  + 1} {\log  {h^{(n_{l_0})}} }.
\ee
The RHS can be made smaller than $\epsilon$ using again Lemma \ref{loglimit} and Corollary \ref{superexp} as long as $l_0\geq l_g$ for some $l_g(\epsilon)\in \mathbb{N}$.
\end{proof}

Recalling the decomposition (\ref{BSdecomp}) of the Birkhoff sums, the estimates of the main contribution and of the gap error in Lemma \ref{estimatemainlemma} and Lemma \ref{estimategaplemma} combined together, for $l_0\geq$ $l(\epsilon)\doteqdot  $$ \max\{ l_m, l_g\}$, yield the estimate in Proposition \ref{alongtower}. 
\end{proofof}

\subsubsection{Growth of Birkhoff sums for other points.}\label{generalx}
In this section we get an estimate for $\BS{u}{r}(x)$ using the estimate 
found in Section \ref{BSalongtowersec} as a fundamental block, i.e. decomposing $\BS{u}{r}(x)$ into pieces which correspond to Birkhoff sums along a tower. It turns out that singular errors from points in the bottom floors of the towers (see the terminology introduced just after (\ref{BSdecomp})) could prevent from getting an estimate of order $r\log r$. In order to get this type of asymptotic, it is necessary to throw away a set of initial points $x$ which has an arbitrarily small measure. 
The integrability condition (\ref{integrability}) of the sequence of balanced times is used in its full strength only in this part.

\paragraph{Preliminary notation.}
Let $\{ n_l \}_{l \in \mathbb{N}} $ be the sequence of balanced times in Proposition \ref{existencebalancedtimes}.
Assume $h^{(n_l)} \leq r < h^{(n_{l+1})}$. 
Define the following sequence $\{ \sigma_l \}_{l\in \mathbb{N}}$, used in the proof of Proposition \ref{BSgeneralxprop} below as a threshold to determine whether $r$ is closer to $h^{(n_l)}$ or to $h^{(n_{l+1})}$.
Let $\tau'$ be such that $\tau/2 < \tau' < 1$, where $\tau$ is the diophantine exponent in (\ref{integrability}) given by Proposition \ref{existencebalancedtimes} and $\tau '$ is well defined since $\tau <2$. Let
\be \label{sigmadef}
\sigma_l = \sigma_l (T) \doteqdot \left( \frac{ \log \| \ZLp{n_{l} , n_{l+1} } \|} {\log h^{(n_l)}}  \right)^{\tau'} ,  \qquad \frac{\tau}{2} < \tau' < 1,
\ee
Clearly $\sigma_l$ depends on the IET $T$ we start with, since the sequence $\{n_l\}_{l\in \mathbb{N}}$ does.

\begin{lemma}\label{sigmarlemma}
The sequence $\{ \sigma_l \}_{l \in \mathbb{N}}$ satisfies the following properties:
\begin{itemize}
\item[(i)] 
$\lim_{l\rightarrow +\infty} \sigma_l =0 $;
\item[(ii)] $ \lim_{l\rightarrow +\infty} \left( \frac{\log \| \ZLp{n_{l} , n_{l+1} } \|}{\log h^{(n_l)} }  \right) \frac{1}{ \sigma_l} =0 $;
\item[(iii)] $\lim_{l \rightarrow +\infty} \sigma_l  \log h^{(n_l)}  = + \infty $;
\item[(iv)] $ \lim_{l\rightarrow +\infty} \sigma_l ^2  h^{(n_{l+1})} \lambda^{(n_l)}  =0$.
\end{itemize}
\end{lemma}

\begin{proof}
Both $(i)$ and $(ii)$ follow from Lemma \ref{loglimit}. 
To show $(iii)$, remark that $ \log \| \ZLp{n_{l} , n_{l+1} } \| \geq \log d \geq 1$, so $\sigma_l \log h^{(n_l)} \geq (\log h^{(n_l)}) ^{1-\tau'} $ and apply Corollary \ref{superexp}.
For $(iv)$, using in order balance (see Lemma \ref{balancehlambda}), the transformation relation for heights, the definition of $\sigma_l$, the diophantine Property (\ref{integrability}) in Proposition \ref{existencebalancedtimes} and Corollary \ref{superexp}, we get 
\bes
\begin{split}
&\sigma_l ^2  h^{(n_{l+1})} \lambda^{(n_l)}  \leq \sigma_l ^2 \kappa \frac{h^{(n_{l+1})}}{h^{(n_l)}} \leq  \sigma_l ^2 \kappa   \| \ZLp{n_{l} , n_{l+1} } \| \leq \\
& \leq   \lim_{l\rightarrow +\infty}   \frac{\kappa ( \log \| \ZLp{n_{l} , n_{l+1} } \| )^{2\tau'}  \| \ZLp{n_{l} , n_{l+1} } \| } { (\log h^{(n_l)})^ {2\tau'}} \leq
  \lim_{l\rightarrow +\infty} 
const \frac{ (\log l) ^{2\tau'} o( l^{\tau} )}{l^{2 \tau ' }} =0,
   \end{split}
\ees
where the last limit is zero since $2\tau'>\tau$.
\end{proof}




\begin{defn} Let $\Sigma^+_{l}=\Sigma^+_l (T)$ be the following set, where $[\cdot]$ denotes the fractional part:
\be \label{badsetforBSu}
\Sigma^+_l (T) 
\doteqdot \bigcup_{i=0}^{[\sigma_l h^{(n_{l+1})}]} T^{-i} [0, \sigma_l \lambda^{(n_l)}].
\ee
\end{defn}
\noindent Remark that, by Property $(iv)$ of Lemma \ref{sigmarlemma},
\be\label{LebSigma}
Leb(\Sigma^+_l) \leq (\sigma_l h^{(n_{l+1})})( \sigma_l \lambda^{(n_l)}) \xrightarrow{l\rightarrow +\infty} 0.
\ee

\begin{prop}[Growth of Birkhoff sums for general points]\label{BSgeneralxprop} Let $T\in \M^+$.
For any $\varepsilon > 0 $ there exists $l_o >0$ such that for $l\geq l_o$, then for each $r \in \mathbb{N}$ and $x\in I^{(0)}$ such that 
\begin{eqnarray} \label{assumptiongeneralxprop}
h^{(n_l)}\leq r < h^{(n_{l+1})} \, &\mathrm{\it and}& \, x \notin  \Sigma^+_{l}(T) , \\
 \label{BSgeneralxestimate}
(1-\varepsilon )r\log r \leq & \BS{u}{r}(x) & \leq 
 (1+\varepsilon ) r\log r + \frac{\kappa +1}{x_m},
\end{eqnarray}
where $x_m \doteqdot \min_{0\leq i < r} T^i x$ and $\kappa$ is given by Proposition \ref{existencebalancedtimes}.
\end{prop}
By adding a small measure set to the excluded set $\Sigma^+_l$ of initial points,  as in Section \ref{mixingpartitionsec}, one can take into account also the term $\kappa/x_m$ and get the asymptotic $r\log r$ for $\BS{u}{r}(x)$.

\begin{figure}\label{OinZ}
\centering
\subfigure[$\Orb{x}{r} \prec Z_1\wedge \dots \wedge Z_4$, ]{\label{OinZ1Z4}{\includegraphics[width=0.45\textwidth]{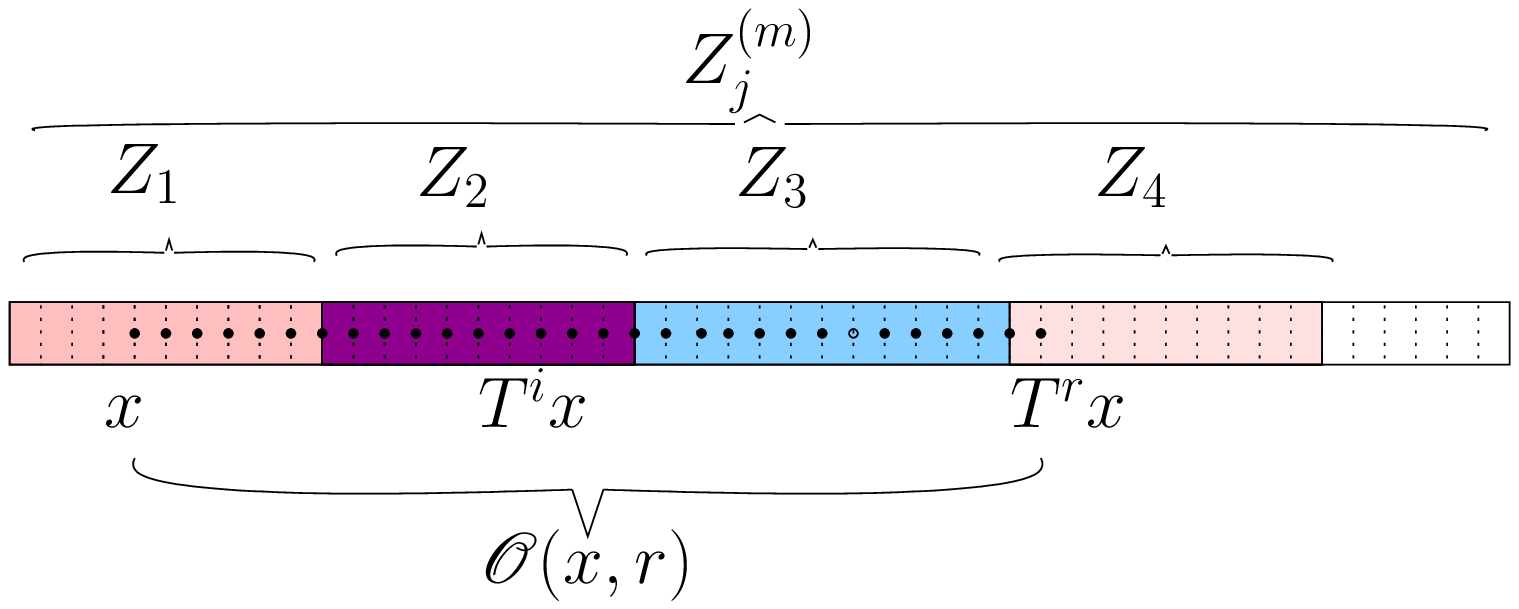}}}
\hspace{5mm}
\subfigure[$ Z_1\wedge Z_2 \prec \Orb{x}{r} $]{\label{Z1Z2inO}{\includegraphics[width=0.45\textwidth]{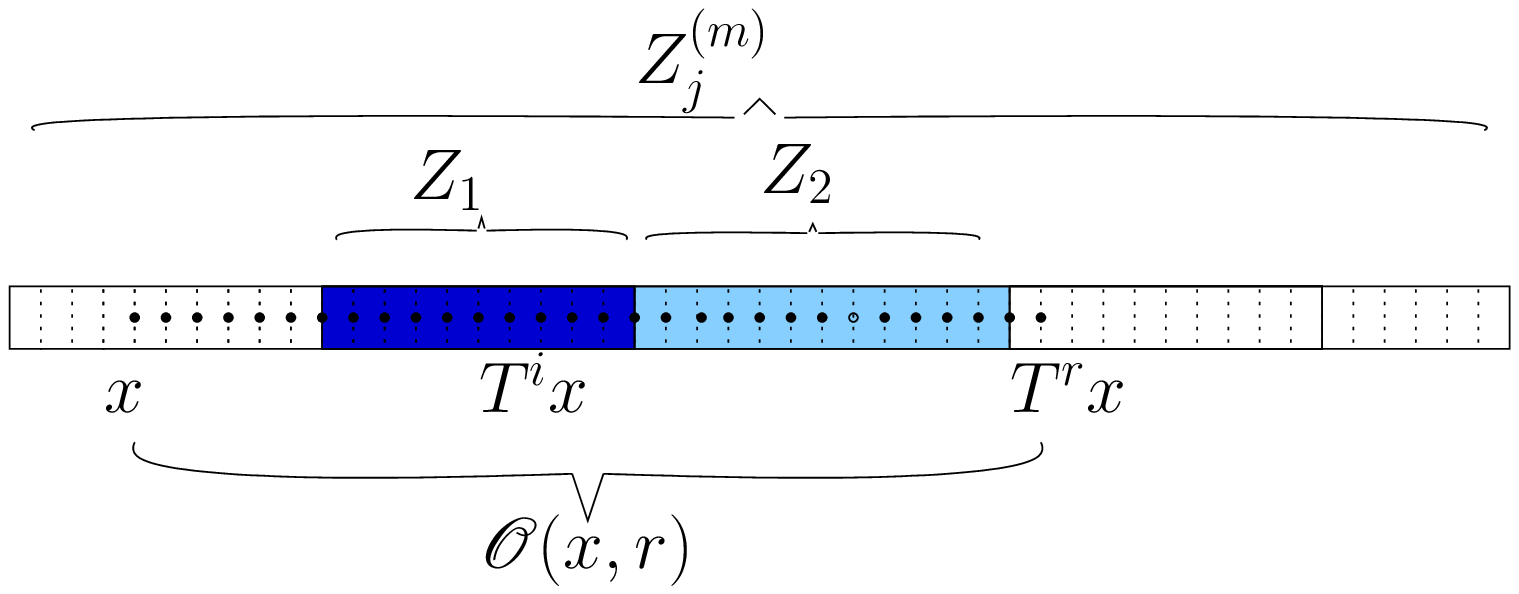}}}
\caption{$\Orb{x}{r} \prec Z^{(m)}_j$  }
\end{figure}

The following notation is used in the proof of Proposition \ref{BSgeneralxprop}.
\paragraph{Notation for approximation by towers.}
Denote by $\Orb{x}{r}$ the orbit segment $\{ T^i x$, $0\leq i < r \}$. Consider a tower $Z^{(m)}_j$. 
In what follows, we write 
\begin{displaymath} 
\Orb{x}{r} \prec Z^{(m)}_j \quad \mathrm{iff} \quad \exists \,k  \,| \quad 0\leq k \leq h_j^{(m)}-r , \quad 
 T^i x \in  T^{k+i}  I^{(m)}_j 
, \quad  0\leq i < r.
\end{displaymath}
When $\Orb{x}{r} \prec Z^{(m)}_j$, 
each point of $\Orb{x}{r}$ is contained in a different floor of $Z^{(m)}_j$ and $T$ acts on the orbit points $T^ix$ ($0\leq i < r-1$) by shifting them to the next floor. (see Figure \ref{OinZ}, where the tower $Z^{(m)}_j$ is drawn horizontally.)

Assume $\Orb{x}{r} \prec Z^{(m)}_j$. 
We write $\Orb{x}{r} \prec Z_1 \wedge Z_2 \wedge \dots \wedge Z_N$, where $Z_i \in \xi_n (Z^{(m)}_j)$, $n\leq m$ (see Section \ref{partitionsec}) for $i=1, \dots, N$, if $\Orb{x}{r} \subset  \cup_{i=1}^{N} Z_i$, $\Orb{x}{r} \cap Z_i \neq \emptyset$ for $1\leq i \leq N$ and moreover $Z_i$ are \emph{consecutive} partition elements of $\xi_n (Z^{(m)}_j)$, i.e., if $h_i$ denotes the height of $Z_i$, remarking that $Z_i\cap I^{(n)}$ is the base of $Z_i$, we have $T^{h_i} Z_i \cap I^{(n)} = Z_{i+1} \cap  I^{(n)}$ for $i=1,\dots, N-1$ (see e.g. Figure \ref{OinZ1Z4}). 
On the other hand, 
we write
 $Z_1 \wedge Z_2 \wedge \dots \wedge Z_N \prec \Orb{x}{r}$, where  $Z_i \in \xi_n (Z^{(m)}_j)$, $n<m$, if $\# \{ Z_i \cap \Orb{x}{r} \} = h_i$ for all $i=1\, \dots, N$, i.e. there is exactly one point of $\Orb{x}{r}$ in each floor of each $Z_i$ and moreover $Z_i$ are, as above, consecutive partition elements. (See Figure \ref{Z1Z2inO}.)




\begin{proofof}{Proposition}{BSgeneralxprop} 
It is always possible to assume that $ \Orb{x}{r} \prec \overline{Z} \doteqdot Z^{(n_{\overline{l}})}_{\overline{j}}$ for some $\overline{l}\geq l+1$, since by choosing $\overline{l}$ such that $\lambda^{(n_{\overline{l}})} < x_m=\min_{0\leq i< r}T^ix$, we assure that $\Orb{x}{r} \cap I^{(n_{\overline{l}})} = \emptyset$.
\paragraph{Orbit decomposition into sums along towers.}
Let us approximate $\Orb{x}{r}$ with elements of $\xi_{n_{l+1}}(\overline{Z})$ and $\xi_{n_{l}}(\overline{Z})$.
Using the assumption $r<h^{(n_{l+1})}$, the cardinality $\# \Orb{x}{r} \cap  I^{(n_{{l+1}})}$ is bounded by $[\kappa]+1$: since the return time to $ I^{(n_{{l+1}})}$ is at least $\min_j{h^{(n_{l+1})}}\geq h^{(n_{l+1})}/\kappa$ by $\kappa$-balance of heights, 
there cannot be more than $r/ \min_j h_j^{(n_{l+1})} + 1 \leq [\kappa]+1 $ returns.

Hence there exists $Z^{l+1}_0, Z^{l+1}_1, \dots, Z^{l+1}_{a^{l+1}+1} \in \xi_{n_{l+1}} ( \overline{Z})$, with $a^{l+1} \leq [\kappa]$ such that (see Figure \ref{orbitdecompfigure})
\be \label{approxl+1}
\Orb{x}{r} \prec Z_0^{l+1} \wedge Z_1^{l+1} \wedge \dots \wedge Z_{a^{l+1}+1}^{l+1}.
\ee
\begin{figure}
\centering
\includegraphics[width=0.85\textwidth]{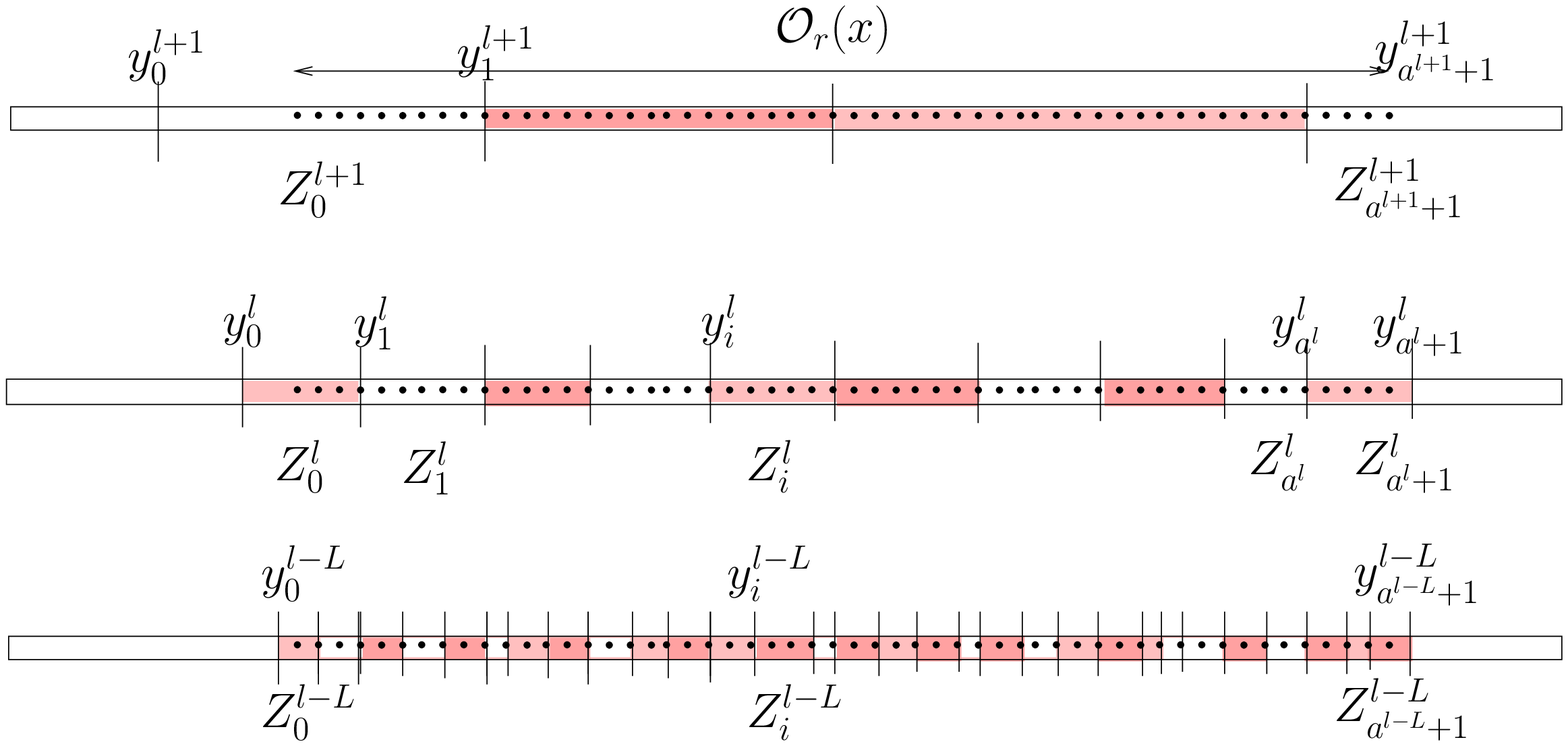}
\hspace{5mm}
\caption{Approximation of $ \Orb{x}{r}$ by elements of $\xi_{n_{l+1}}$, $\xi_{n_{l}}$ and $\xi_{n_{l-L}}$. }
\label{orbitdecompfigure}
\end{figure}
\noindent Approximating $ \Orb{x}{r} $ also with elements $Z^{l}_i \in \xi_{n_{l}} ( \overline{Z})$ (see again Figure \ref{orbitdecompfigure}):
\be \label{apprl}
 Z_1^{l} \wedge Z_2^{l} \wedge \dots \wedge Z_{a^{l}}^{l}  \prec \Orb{x}{r} \prec Z_0^{l} \wedge Z_1^{l}  \wedge \dots \wedge Z_{a^{l}}^{l} \wedge  Z_{a^{l}+1}^{l}.
\ee 
%
We still need another level of approximation. Let $L \doteqdot l_1 \overline{l} \in \mathbb{N}$ where $l_1$ is such that $2\kappa /d^{l_1} < \epsilon$. We can find elements  $Z^{l-L}_i \in \xi_{n_{l-L}} ( \overline{Z})$ such that (Figure \ref{orbitdecompfigure})
\be\label{apprl-L} 
 Z_1^{l-L} \wedge Z_2^{l-L} \wedge \dots \wedge Z_{a^{l-L}}^{l-L}  \prec \Orb{x}{r} \prec Z_0^{l-L} \wedge Z_1^{l-L}  \wedge \dots \wedge Z_{a^{l-L}}^{l-L} \wedge  Z_{a^{l-L}+1}^{l-L}.
\ee
Denote by $h_i^{l}$ and $h_i^{l-L}$ the heights of $Z_i^{l}$ and $Z_i^{l-L}$ respectively. Remark that 
\be \label{restimate}
 \sum_{i=1}^{a^{l}} h_{i}^{l} \leq    \sum_{i=1}^{a^{l-L}} h_{i}^{l-L} \leq r \leq \sum_{i=0}^{a^{l-L}+1} h_{i}^{l-L} \leq \sum_{i=0}^{a^{l}+1} h_{i}^{l} .
\ee



Let us truncate $\Orb{x}{r}$ into segments contained in different elements $Z^{l-L}_i \in \xi_{l-L}(\overline{Z})$. Remarking that $I^{(n_{l-L})}$ contains all the bases of the towers $Z_i^{l-L}$, let us denote 
\be \label{defyil-L}
y^{l-L}_i \doteqdot  \Orb{x}{r} \cap  Z_i^{l-L} \cap I^{(n_{l-L})}; \quad i=1,\dots, a^{l-L}.
\ee
Since $h^{l-L}_i$ is exactly the first return time of $y^{l-L}_{i}$ to $I^{(n_{l-L})}$, $y^{l-L}_{i+1} = T^{h^{l-L}_i}  y^{l-L}_{i} =  T^{(n_{l-L})}  y^{l-L}_{i}$. 
Add also the two auxiliary points:
\be \label{defy0l-L}
y^{l-L}_0 \doteqdot (T^{(n_{l-L})})^{-1}  y^{l-L}_{1}; \qquad y^{l-L}_{ a^{l-L}+1} \doteqdot T^{(n_{l-L})}  y^{l-L}_{ a^{l-L}}.
\ee
From (\ref{apprl-L}), 
$ \bigcup_{i=1}^{a^{l-L}} \Orb{y_{i}^{l-L}}{h^{l-L}_{i}} \subset \Orb{x}{r} \subset \bigcup_{i=0}^{a^{l-L}+1} \Orb{y_{i}^{l-L}}{h^{l-L}_{i}}.$
As a consequence, since $u>0$, we get the following estimate for $\BS{u}{r}(x)$:
\be\label{BSgeneralxestimate2}
\sum_{i=1}^{a^{l-L}} \BS{u}{h_i^{l-L}}(y_i^{l-L}) \leq \BS{u}{r}(x) \leq \sum_{i=0}^{a^{l-L}+1} \BS{u}{h_i^{l-L}}(y_i^{l-L}) .
\ee
Each term in the summations in (\ref{BSgeneralxestimate2}) is a Birkhoff sum along a tower of step $n_{l-L}$. Hence we can apply Proposition \ref{alongtower} to each term and find $l_0\doteqdot l(\epsilon)+L$ such that for each $l\geq l_0$ we get
\begin{eqnarray}
\BS{u}{r}(x) &\geq & (1-\epsilon)
\sum_{i=1}^{a^{l-L}} h_{i}^{l-L}   \log {{ h^{(n_{{l}-L})}}} + \sum_{i=1}^{a^{l-L}} \frac{1}{y_i^{l-L}} \label{BSgenxbb} \\
\BS{u}{r}(x)  &\leq& (1+\epsilon)
\sum_{i=0}^{a^{l-L}+1} h_{i}^{l-L}   \log   h ^{(n_{{l}-L})}  
+  \sum_{i=0}^{a^{l-L}+1} \frac{1}{y_i^{l-L}}. \label{BSgenxba}
\end{eqnarray}
Let us refer to the first term in the LHS of (\ref{BSgenxbb}) or (\ref{BSgenxba}) as \emph{ergodic term} and to the last term, i.e. the contributions of points in the bottom floors, as \emph{resonant term}\footnote{There is again an analogy with the terminology used by \cite{Ko:nonI, Ko:nonII}.}.
\paragraph{Ergodic term.}
Taking the ratio of the ergodic term over $r\log r$ and applying the bounds (\ref{restimate}) for $r$:
\begin{eqnarray}
\frac{\left( \sum_{i=1}^{a^{l-L}} h_{i}^{l-L} \right)  \log  h^{(n_{{l}-L})}}{r\log r } \geq  \left( 1- \frac{2  h^{(n_{l-L})}}{r} \right) \frac{  \log  h^{(n_{{l}-L})}}{\log r }  \label{ergodicb} \\
\frac{\left( \sum_{i=0}^{a^{l-L}+1} h_{i}^{l-L}\right) \log  h^{(n_{{l}-L})}}{r\log r } \leq  \left( 1+ \frac{2  h^{(n_{l-L})}}{r} \right) \frac{  \log  h^{(n_{{l}-L})}}{\log r }  \label{ergodica}
\end{eqnarray}
By assumption (\ref{assumptiongeneralxprop}) on $r$, Property $3$, Lemma \ref{productpositive}, balance and choice of $L=l_1\overline{l}$, 
\bes \frac{2  h^{(n_{l-L})}}{r} \leq \frac{2  h^{(n_{l-l_1\overline{l}})}}{h^{(n_{l})}} \leq \frac{2\kappa}{d^{l_1}} < \epsilon.
\ees
Hence the first factors in the RHS of (\ref{ergodicb},\ref{ergodica}) are bounded respectively below by $(1-\epsilon)$ and above by $(1+ \epsilon)$.  
The second factor in the RHS of (\ref{ergodica}) is trivially less than one. 
From $ r < h^{(n_{{l}+1})} $,
\be \label{ergodicerror}
 \frac{  \log  h^{(n_{{l}-L})}}{\log r } \geq \frac{ \log  h^{(n_{{l+1}})} - \log \frac{ h^{(n_{{l+1}})}}{ h^{(n_{{l}-L})}}}{\log  h^{(n_{{l}+1})} } .
\ee
Since by the heights transformation formula and Lemma \ref{loglimit} 
\bes
\frac{ \log \frac{ h^{(n_{{l+1}})}}{ h^{(n_{{l}-L})}}}{\log  h^{(n_{{l}+1})} } \leq  \frac{\log \| \ZLp{n_{{l-L}} ,  n_{{l+1}}} \|}{\log  h^{(n_{{l}+1})} }  \xrightarrow{l\rightarrow \infty} 0,
\ees
also the second factor in the RHS of (\ref{ergodica}) is bounded from below by $(1-\epsilon)$ if $l\geq l_o$ for some $l_o\geq l_0$.


So far, combining (\ref{BSgenxbb}, \ref{BSgenxba}) with (\ref{ergodicb}, \ref{ergodica}) and (\ref{ergodicerror}) and using the fact that  
 the resonant term is positive, we proved that
\be \label{ergodicestimate}
(1-\epsilon)^2 r \log r
 \leq { \BS{u}{r}(x)}   \leq  (1+ \epsilon) r\log r  + { \sum_{j=0}^{a^{l-L}+1}  \frac{1}{y_j^{l-L}} }.
\ee

\paragraph{Resonant term.}
We want to prove the following estimate for the resonant term:
\be \label{estimateresonant} 
0\leq  \sum_{j=0}^{a^{l-L}+1}  \frac{1}{y_j^{l-L}} \leq \epsilon r\log r + \frac{\kappa+1}{x_m}.
\ee
Let us first group $\{ y_j^{l-L} \}_{j=0,\dots , a^{l_l}+1}$ according to visits to different elements of the partition $\xi_{n_l}(\overline{Z})$. 
Since 
\bes
\bigcup_{j=0}^{a^{l-L}+1} y_j^{l-L} \subset \bigcup_{j=0}^{a^{l-L}+1} Z_j^{l-L} \subset \bigcup_{i=0}^{a^{l}+1}  Z_i^{l}. 
\ees
we have the estimate
\bes
\sum_{j=0}^{a^{l-L}+1}  \frac{1}{y_j^{l-L}} \leq \sum_{i=0}^{a^{l}+1} \sum_{{\begin{subarray}{c} y_j^{l-L} \in Z_i^{l} \\ j=0,\dots,a^{l-L}+1  \end{subarray}} }   \frac{1}{y_j^{l-L}}.
\ees
Each of the points $y_j^{l-L} \in Z_i^{l}$ belongs to a different floor of a tower of step $n_l$. Hence, arguing as in the Lemma \ref{estimategaplemma} about the Gap Error, each of the terms $\sum_{y_j^{l-L} \in Z_i^{l} }   \frac{1}{y_j^{l-L}}$ can be bounded from above by applying Lemma \ref{Kocharithmetic} to an auxiliary arithmetic progression with step $1/d \kappa \nu h^{l}_i$. The cardinality of points in each group, by (\ref{entriesBmn}), is bounded by
$\# \{ y_j^{l-L}  | y_j^{l-L} \in Z_i^{l} \} \leq \# \{  Z_j^{l-L} \in \xi_{n_{l -L}}( Z_i^{l} ) \} \leq \| \ZLp{n_{l-L}, n_l} \|.
$ 
The initial point $\min \{ y_j^{l-L} | \,  y_j^{l-L} \in Z_i^{l} \}$ is given by the only visit to the base $ Z_i^{l} \cap I^{(n_{l})}$. Denote, as above (see (\ref{defyil-L}, \ref{defy0l-L})),
\bes
y^{l}_i \doteqdot  \Orb{x}{r} \cap  Z_i^{l} \cap I^{(n_{l})}, \quad i=1,\dots , a^{l}; \quad y^{l}_0 \doteqdot (T^{(n_{l})})^{-1}  y^{l}_{1}; \qquad y^{l}_{ a^{l}+1} \doteqdot T^{(n_{l})}  y^{l}_{ a^{l}}.
\ees
Hence we get
\be \label{resonantgrouping}
\sum_{j=0}^{a^{l-L}+1}  \frac{1}{y_j^{l-L}} \leq  \sum_{i=0}^{a^{l}+1}   
d \kappa \nu h^{l}_i ( \log \| \ZLp{n_{l-L}, n_l} \|+1) +
 \sum_{i=0}^{a^{l}+1} \frac{1}{y_i^{l}}.
\ee
Comparing the first term in the RHS of (\ref{resonantgrouping}) to $r\log r$ and recalling (\ref{restimate}) we get
\bes\label{secondaryresonaces}
\frac{ d \kappa \nu \left( \sum_{i=0}^{a^{l}+1}  h^{l}_i \right) ( \log \| \ZLp{n_{l-L}, n_l} \|+1) }{r\log r} \leq d \kappa \nu \left(1+\frac{2h^{(n_l)}}{r}\right) \frac{ \log \| \ZLp{n_{l-L}, n_l} \|+1  }{\log h^{(n_l)}} 
\ees
where $(1+2h^{(n_l)}/r)\leq 3$, so the last term, enlarging $l_o$ if necessary, is less than $\epsilon$ when $l\geq l_o$ by Lemma \ref{loglimit} and Corollary \ref{superexp}. 

The second term in the RHS of (\ref{resonantgrouping}) is bounded in two different ways, according to the ratio between $r$ and $h^{(n_{l+1})}$, using the quantity $\sigma_l$ defined in (\ref{sigmadef}) as a threshold.

{\it Case 1}. Assume $\sigma_l h^{(n_{{l}+1})} \leq r < h^{(n_{{l}+1})}$. 
Recalling (\ref{approxl+1}), 
\bes
\{ y_i^{l} |\,\, i=0 , \dots , a^{l}+1  \} \subset \bigcup_{j=0}^{a^{l+1}+1} \{  y_i^{l} |\,\, y_i^{l}  \in  \cap  Z_j^{l+1} , \,  i=0 , \dots , a^{l}+1 \}.
\ees
To estimate the contribution from each of the sets in the RHS of (\ref{resonantgrouping}), arguing as above, consider an auxiliary arithmetic progression of step $ d \kappa \nu h^{l+1}_j$. The closest point of each set is given by the visit to $I^{(n_{l+1})}$. Remarking that the number of points in each is bounded by
$\# \{ y_i^{l}|\, y_i^{l}  \in   Z_j^{l+1} ,\, i=0,\dots, a^{l}+1   \} \leq \# \xi_{n_l}( Z_j^{l+1} )  \leq \| \ZLp{n_{l}, n_{l+1}} \|,
$ 
we get
\be \label{mainres1} 
 \sum_{i=0}^{a^{l}+1} \frac{1}{y_i^{l}}
\leq  \sum_{j=0}^{a^{l+1}+1}    
d \kappa \nu h^{l+1}_j ( \log \| \ZLp{n_{l}, n_{l+1}} \|+1)  + 
\sum_{\begin{subarray}{c} y_i^{l} \in I^{(n_{l+1})} \\ i=0,\dots, a^{l}+1 \end{subarray}}  \frac{1}{y_i^{l}}.
\ee
For the first term in the RHS of (\ref{mainres1}), by the assumptions on $r$, $ h^{l+1}_j/r \leq  h^{(n_{l+1})}/r \leq 1/\sigma_l$ and $a^{l+1} \leq \kappa $, 
\bes
\frac{\sum_{j=0}^{a^{l+1}}    d \kappa \nu h^{l+1}_j ( \log \| \ZLp{n_{l}, n_{l+1}} \|+1)}
{r \log r}\leq \frac{(\kappa+2 ) d \kappa \nu ( \log \| \ZLp{n_{l}, n_{l+1}} \|+1)}{\sigma_l \log h^{(n_{l})}} .
\ees
The last expression, enlarging again $l_o$ if necessary, is smaller than $\epsilon$ if $l \geq l_o$ by Property $(ii)$ and $(iii)$ in Lemma \ref{sigmarlemma}.

The second term in the LHS of (\ref{mainres1}) can be just estimated with $(\kappa +1)/  x_m$ since, as remarked at the beginning of this proof, $\# \Orb{x}{r} \cap  I^{(n_{l+1})} \leq [\kappa] $. This completes the proof of (\ref{estimateresonant}) in this case.

{\it Case 2}. Assume $ h^{(n_{l})} \leq  r < \sigma_l h^{(n_{l+1})} $. In this case, use the trivial estimate 
\bes \label{trivialresbnd}
 \sum_{i=0}^{a^{l}+1} \frac{1}{y_i^{l}}\leq  (a^{l}+2) \frac{1}{x_m}.
\ees
Since $x_m= T^i x$ for some $0\leq i < r$ and in this case $r< \sigma_l h^{(n_{l+1})}$, by the assumption (\ref{assumptiongeneralxprop}) and the definition (\ref{badsetforBSu}) of $\Sigma^+_l$,
\bes
x_m \geq \sigma_l \lambda^{(n_l)} \geq  \sigma_l  \frac{1}{\kappa h^{(n_l)}}, 
\ees
where the last inequality uses the balance of $n_l$ (see Lemma \ref{balancehlambda}).
Moreover, from (\ref{restimate}) and $\kappa$-balance of heights,  
\bes
a^{l}\leq \frac{r}{\min_j h^{(n_{l})}_j} \leq  \frac{\kappa r}{ h^{(n_{l})}}.
\ees
Hence,
\bes
\frac{  \sum_{i=0}^{a^{l}+1} \frac{1}{y_i^{l}}  }{r\log r} \leq \frac{ (\frac{\kappa r}{ h^{(n_{l})}} +2)  \frac{\kappa h^{(n_l)}}{\sigma_l} }{r\log r } \leq \frac{\kappa^2+2\kappa}{\sigma_l \log h^{(n_l)}},
\ees
which, again enlarging $l_o$, is smaller than $\epsilon$ for $l\geq l_o$ by Property $(iii)$ in Lemma \ref{sigmarlemma}.

In both cases we proved the estimate (\ref{estimateresonant}) for the resonant term. Together with (\ref{ergodicestimate}), for an appropriate choice of $\epsilon$, this completes the proof of Proposition \ref{BSgeneralxprop}.
\end{proofof}

\begin{cor}\label{BSr^2}
Let $T\in \M^+$. For each $\varepsilon >0$ there exists $r_0$ such that for all $r\geq r_0$, $x\in I^{(0)}$, 
\be
\BS{u}{r}(x) \leq \varepsilon r^2 + \frac{\kappa +1}{x_m}.
\ee
\end{cor}
The estimate in the Corollary is worst than (\ref{BSgeneralxestimate}) in Proposition \ref{BSgeneralxprop}, but holds for all points and is used in Section \ref{preliminarypartitionsec}.
\begin{proof} Let $h^{(n_l)}\leq r < h^{(n_{l+1})}$
Remark that equations (\ref{ergodicestimate}), (\ref{resonantgrouping}) and (\ref{mainres1}) in the proof of Proposition \ref{BSgeneralxprop} were obtained without using the assumption (\ref{assumptiongeneralxprop}) and hence still hold if $r\geq r_0\doteqdot h^{(n_{l_o})}$. Terms estimated by $r\log r$ are clearly less than $\varepsilon  r^2$ choosing $r_0$ large enough. The second term in the RHS of (\ref{mainres1}) is estimated by $(\kappa+1)/x_m$. Let us estimate the first term in the RHS of (\ref{mainres1}) by
\bes
\frac{h^{(n_{l+1})} \log \| \ZLp{n_{l}, n_{l+1}}  \|}{r^2} \leq \frac{h^{(n_{l+1})} \log \| \ZLp{n_{l}, n_{l+1}} \ \|}{{h^{(n_{l})}}^2} \leq \frac{  \| \ZLp{n_{l}, n_{l+1}}  \| \log \| \ZLp{n_{l}, n_{l+1}} \ \|}{const \, d^l}
\ees 
where we used Corollary \ref{superexp} in the last bound. The limit of this ratio as $l\rightarrow \infty$ is zero by (\ref{integrability}) in Property 4. 
\end{proof}

\subsection{Growth of Birkhoff sums of the derivatives.}\label{BSf'andf''sec}
Let $h_\I^{(n)}$ and $\lambda_\I^{(n)}$ be the sequences of heights and lenghts of towers for $T^{\I}$. For $T^{\I} \in \M^+$, let $\{n'_{l'}\}_{l'\in \mathbb{N}}$ be the sequence of balanced times for $T^{\I}$ given by Proposition \ref{existencebalancedtimes}. Let $\sigma_{l'}=\sigma_{l'}(T^{\I})$. Define
\be\label{DefSigma-}
\Sigma^-_{l'} (T) 
\doteqdot \bigcup_{i=0}^{[\sigma_{l'} h_\I^{(n'_{l'+1})}]} T^{-i} [1- \sigma_{l'} \lambda_\I^{(n'_{l'})},1].
\ee
\begin{cor}[Growth of  $\BS{v}{r}$]\label{BSv} Let $T\in \M^-$. 
For any $\epsilon > 0 $ there exists $l_o' >0$ such that, if $l'\geq l_o'$, for any $r \in \mathbb{N}$ and $x\in I^{(0)}$ such that 
\bes
h_\I^{(n'_{l'})} < r < h_\I^{(n'_{{l'}+1})}  \qquad \mathrm{and} \qquad x \notin  \Sigma^-_{l'}(T) ,
\ees 
denoting by $x_M \doteqdot \max_{0\leq i < r}  T^i x$, we get \footnote{Here $\kappa'$ is the same given by Proposition \ref{existencebalancedtimes}.}
\bes 
(1-\epsilon)r\log r \leq  \BS{v}{r}(x)  \leq 
 (1+\epsilon) r\log r + \frac{\kappa'+1}{1-x_M}.
\ees
\end{cor}
\begin{proof} Corollary \ref{BSv} is simply obtained by restating Proposition \ref{BSgeneralxprop} for $u$ and $T^{\I}$ and using the relation with $v$ given by (\ref{BSuvrel}). Remark that $ \min_i \left( T^{\I} \right)^i(1-x) = \min_i \I T^i (x) = \min_i (1-T^i x) = 1-x_M$
and that $(1-x)\in \Sigma^+_{l'}(T^{\I})$ if and only if $x\in \Sigma^-_{l'}(T)$.
\end{proof}
We can combine the estimates for $\BS{u}{r}$ and $\BS{v}{r}$ as follows.
\begin{cor}[Growth $\BS{f'}{r}$]\label{BSf'cor}
For each $T\in \M$ and $C^+\neq C^-$, there exists $C_1,C'_1 , C_2,C_3 >0$ and $r_o$ such that for $r\geq r_o$, if 
\be\label{assumptionsBSf'}
\begin{array}{ccc}
\begin{array}{rcl}
 h^{(n_l)}&\leq r < &h^{(n_{l+1})}, \\
 h_\I^{(n'_{l'})}&\leq r < &h_\I^{(n'_{{l'}+1})} 
\end{array} & \mathrm{and} &
x \notin  \Sigma^+_{l}(T) \cup \Sigma^-_{l'}(T)  ,
\end{array}
\ee 
and $x$ is not a singularity of $\BS{f}{r}$,\footnote{Here $x_m$, $x_M$, $\kappa$ and $\kappa'$ are as in Corollary \ref{BSv}.}
\begin{eqnarray}\label{BSf'}
\BS{f'}{r}(x)  &\leq& 
 -C_1 r\log r + \frac{C^+(\kappa'+1)}{1-x_M}, \qquad \mathrm{if} \,\, C^+>C^-;
\\
 \BS{f'}{r}(x)  &\geq& 
 C_1 r\log r - \frac{C^-(\kappa+1)}{x_m}, \qquad \mathrm{if} \, \, C^+ < C^- ; 
\\
| \BS{f'}{r}(x) | &\leq& \label{BSf'upperbound}
 C_2 r\log r +\frac{C_2 (\kappa+1)}{x_m} + \frac{C_2 (\kappa'+1)}{1-x_M} .
\end{eqnarray}
Moreover, for all $x\in I^{(0)}$ different from singularities of $\BS{f}{r}$,
\be \label{roughboundf'}
| \BS{f'}{r}(x)|  \leq 
 C_3 r^2 +  \frac{C_3( \kappa+1)}{x_m} + \frac{C_3 (\kappa'+1)}{1-x_M}.
\ee
\end{cor}
\begin{proof}
Assume $C^+>C^-$. Consider the sequence  $\alpha_r \rightarrow 0$ in Proposition \ref{reducetouv}. 
One can choose $r_1$ so that, for each $r \geq r_1$, we have $(C^+ - \alpha_r ) > (C^- + \alpha_r)$.
Hence it is also possible to choose $\epsilon >0$ so that, for $r \geq r_1$, 
$C_1 \doteqdot (C^+ - \alpha_r )(1-\epsilon)  - ( C^- + \alpha_r ) (1+ \epsilon)  > 0$.  
By Proposition \ref{reducetouv}, Proposition \ref{BSgeneralxprop} and Corollary \ref{BSv}, which can be applied by the assumptions (\ref{assumptionsBSf'}) when $r\geq r_o$, $r_o\doteqdot \max \{ r_1, h^{(n_{l_o})},  h_\I^{(n'_{l'_o})} \}$, 
\bes
\begin{split}
&\BS{f'}{r}(x)  \leq -( C^+\!-\alpha_r )  \BS{u}{r}(x) + (  C^-\!+\alpha_r ) \BS{v}{r}(x) \leq  -( C^+\!-\alpha_r ) {(1-\epsilon)} r \log r + \\& + (  C^-\!+\alpha_r ) (1+ \epsilon)   r \log r + \frac{(C^-\!+\alpha_r) (\kappa'+1)}{1-x_M}  \leq - C_1  r \log r + \frac{C^+  (\kappa'+1)}{1-x_M} .
\end{split}
\ees
The case $C^+<C^-$ can be treated analogously.

Also (\ref{BSf'upperbound}) follows similarly: enlarging $r_1$ so that if $r\geq r_1$, $\alpha_r \leq \min (C^+, C^- )$, by Proposition \ref{reducetouv},
\be\label{bounduv}
| \BS{f'}{r}(x) | \leq ( C^+ + C^- ) (  \BS{u}{r}(x) + \BS{v}{r}(x) )
\ee
and $\BS{u}{r}$, $\BS{v}{r}$ can be again estimated by Proposition \ref{BSgeneralxprop} and Corollary \ref{BSv}.

For (\ref{roughboundf'}), apply to (\ref{bounduv}) the rough estimate on $\BS{u}{x}$ for all points given by Corollary \ref{BSr^2}  and the analogous one for $\BS{v}{x}$ which follows from $T\in \M^-$.

\end{proof}
\begin{cor}[Growth $\BS{f''}{r}$]\label{growthBSf''}
For each $T\in \M$ and $C^+\neq C^-$, there exists $C_4 >0$ and $r_o$ such that for $r\geq r_o$, if (\ref{assumptionsBSf'}) holds,
\bes
| \BS{f''}{r}(x) | \leq C_4
 \max \left\{\frac{1}{x_m}, \frac{1}{1-x_M} \right\} \left( r\log r + \frac{\kappa+1}{x_m} + \frac{\kappa'+1}{1-x_M} \right).
\ees
\end{cor}
\begin{proof}
By definition of logarithmic singularity, there exists $\delta>0$ such that $0\leq f''(x)\leq 2C^+/x^2$ if $x<\delta$ and $0\leq f''(x)\leq 2C^-/(1-x)^2$ if $x >1-\delta$. Let $M_{f''}$ be the maximum of $|f''|$ on $[\delta, 1-\delta]$. Hence, for each $x$ not a singularity of $\BS{f}{r}$,
\be \label{trivialbndf''}
\begin{split}
|\BS{f''}{r}(x)| & \leq 2C^+ \BS{{1}/{x^2}}{r}(x) +  2C^- \BS{{1}/{(1-x)^2}}{r}(x) + r M_{f''} \leq \\
& \leq  2C^+ \frac{1}{x_m} \BS{u}{r}(x) +  2C^-  \frac{1}{1-x_M} \BS{v}{r}(x)  + r M_{f''}   .
\end{split}
\ee 
Applying Proposition \ref{BSgeneralxprop} and Corollary \ref{BSv} one get the desired estimate.  
\end{proof}


\section{Construction of the mixing partitions.}\label{mixingpartitionsec}
In this Section we construct the partitions $\eta_{m}(t)$ which verify the mixing criterium (see Lemma \ref{mixingcriterium}). 
The construction is carried out in three main steps, formulated in Section \ref{propertiessec} as Propositions \ref{preliminarypartitionprop}, \ref{growthBSf'partition} and \ref{mainpropertiespartitionprop}. Their proofs are in Sections \ref{preliminarypartitionsec}, \ref{stretchingpartitionsec} and \ref{mainpropertiespartitionsec} respectively. We anticipate in Section \ref{areaestimatesec} the final area estimates, which conclude the proof of Theorem \ref{mixing}.

\subsection{Partitions properties.}\label{propertiessec}
As a preliminary step, we construct in Section \ref{preliminarypartitionsec} partial partitions which satisfy the following proposition.
Denoting $[\cdot]$ the integer part, let 
\be \label{RMdef}
R_M(t) \doteqdot \left[ t/m_f
 \right] + 2.
\ee
\begin{prop}[Preliminary partitions.]\label{preliminarypartitionprop}
For each $0<\delta<1$ and $M>1$, there exist $t_0>0$ and partial partitions $\eta_p(t)$ for $t \geq t_0$, such that $Leb(\eta_p(t)) > 1-\delta$ 
and the following properties hold for each $I=[a,b[ \, \in \eta_p(t)$.
\begin{enumerate}
\item \emph{Continuity intervals:}

$T^j$ is continuous on $[a,b]$ for each $0\leq j \leq  R_M(t)$; 
\item \emph{Control of interval sizes:}  
\bes \frac{1}{t \log\log t} \leq Leb(I) \leq \frac{2}{t \log\log t };
\ees
\item \emph{Control of the distance from singularities:}
\bes dist (T^j I, 0 ) \geq \frac{M}{t \log \log t}, \quad dist (T^j I, 1 ) \geq \frac{M}{t \log \log t}, \quad  0\leq j \leq  R_M(t);
\ees
\item \emph{Control of the number of discrete iterations:} 
\bes
\frac{t}{3} \leq r(x,t) \leq R_M(t) \leq \frac{2}{{m_f}} t, \qquad \forall x \in I  .
\ees
\end{enumerate}
\end{prop}


Assume now that $T\in \M$. 
For definiteness, assume also that the asymmetry constants of the roof function satisfy $C^+>C^-$.  
Using the estimates on the growth of Birkhoff sums obtained in Section \ref{growthBSsec}, we can refine the partitions $\eta_p(t)$ to get the following.
\begin{prop}[Stretching partitions.]\label{growthBSf'partition}
For $T\in \M$ and $C^+> C^-$, there exists $C_1', C_2', C'' >0$ such that for each $0<\delta<1$, $M>1$, if $\eta_p(t)$ are the corresponding partitions in Proposition \ref{preliminarypartitionprop}, 
there exists a sequence of refined partitions $\eta_{s}(t) \subset \eta_p(t)$ with $Leb(\eta_{s}) > Leb(\eta_{p}) -\delta$ and there exists $t_1>t_0$ such that, when $t\geq t_1$, for any $x \in \eta_{s}(t)$ and 
 integer $r$
 with  $t/3  \leq r \leq 2t / m_f $,
\begin{eqnarray} \label{f'decreases}
\BS{f'}{r}(x)  &\leq &- C_1'\,
r \log r ;\\
| \BS{f'}{r}(x) |  &\leq & C_2'\,
r \log r   \label{SBf'boundabove}; \\
\label{BSf''}
\BS{f''}{r}(x) & \leq &\frac{ C''}{M} r^2 (\log r) ( \log \log r ).
\end{eqnarray}
\end{prop}
Let us show that Proposition \ref{growthBSf'partition} implies in particular that $r(\cdot,t)$ is an increasing function on each interval $I\in \eta_{s}(t)$ for $t\geq t_1$. Assume $x< y$ are points of $I$. 
Since in particular $\BS{f'}{r(x,t)}<0$, the function $\BS{f}{r(x,t)}$ is strictly decreasing. Hence
$\BS{f}{r(x,t)}(y) <  \BS{f}{r(x,t)}(x) \leq t $. 
Using again the definition of $r(\cdot,t)$, we get $r(y,t)\geq r(x,t)$.


\paragraph{Geometric description of the dynamics of partition elements.} \label{geometricdescriptionsec}
Let $I=[a,b[$ be an element of the partition $\eta_{s}(t)$.
Consider $\varphi_t I$ and let us give first a geometric description of $\varphi_t I$ for $t \gg 1$. 
For $C^+ > C^-$, 
as just proved, 
 $r(\cdot,t)$ is an increasing function on each $[a,b[\in \eta_p(t)$. Hence, let\footnote{The dependence on $t$ will be omitted when $t$ is clear from the context.}
\begin{eqnarray} \label{rarbmaxmin}
&&r(a)= r(a,t)= \inf_{x\in [a,b[} r(x,t); \qquad r(b) = r(b,t)=\sup_{x\in [a,b[} r(x,t);
\\ \nonumber
&&J=J([a,b[,t) \doteqdot r(b,t) -r(a,t) +1 .
\end{eqnarray}
The image $\varphi_t I$ splits into several curves and $J$ gives exactly their cardinality. More precisely, consider the equation $\BS{f}{r}=t$ on $I$, which has a solution exactly for $r=r(a)+1, \dots, r(b)$, unique by monotonicity. Denote by $y_j$ the solution of
\be \label{splittingpointsdef}
\BS{f}{r(a)+j}(y_j)=t , \qquad j = 1 , \dots, J-1, 
\ee
so that by (\ref{flowdef}), $r(y_j,t)=r(a)+j$ and 
$
\varphi_t (y_j, 0 ) =( T^{r(y_j,t)}(y_j) , 0).
$ 
The points $a\doteqdot y_0 \leq y_1 \leq  \dots \leq y_j \leq y_{j+1} \leq y_{J-1} \leq y_{J} \doteqdot b $
 are splitting points, meaning that the image $\varphi_t I$ consists of $J$ curves, which are the graphs of $t - \BS{f}{r(a)+j}$ restricted to $I_j \doteqdot [y_{j}, y_{j+1}[$, $j=0,\dots, J-1$. Each curve projects to $T^{r(a)+j} (I_j)$, by (\ref{flowdef}). In particular, 
 the curves project to the orbit $T^{r(a)+j} (I)$, for $j=0,\dots, J\!-\!1$. 


\paragraph{Properties of the mixing partitions.}
Given $\delta>0$, if $[b_1,b_2]$ is the base of the rectangle $R$, denote $ \chi$ the indicator of $[b_1+\delta,b_2-\delta]$. Choose $t_2 \geq t_1$ so that $2/(t_2 \log \log t_2) < \delta$ and the mesh of the partitions $\eta_{p}(t)$ for $t\geq t_2$ is bounded by $\delta$ by Property 2 in Proposition \ref{preliminarypartitionprop}. Denoting by $h(R)$ the height of the rectangle $R$, for $j=0,\dots, J\!-\!1$, let
\be \label{Ijhdef}
I_{j}^{h(R)} \doteqdot \{ x | \quad t - h(R)  \leq \BS{f}{r(a)+j}(x) \leq t  \}.
\ee
Points in $I_{j}^{h(R)}$ are the ones that reach the correct height to intersect $R$, i.e. if $x\in I_{j}^{h(R)}$, then $\varphi_t(x,0)$ is contained in the horizontal strip $I^{(0)}\times h(R)$ (as shown in the proof of Lemma \ref{areaboundlemma}).
%
For $j=0,\dots, J-1$, 
denote by 
\bes
\Delta f^j = \Delta f^j([a,b[,t) \doteqdot \BS{f}{r(a)+j}(a) -  \BS{f}{r(a)+j}(b) ; \label{Deltafjdef} \qquad \Delta f 
\doteqdot \Delta f^0. \ees
Remark that $\Delta f^j\geq 0$. The quantity $\Delta f$ express the the delay accumulated between the endpoints in time $t$. 
Also, the quantity $\Delta f^j$ gives the vertical stretch of the graph of $t - \BS{f}{r(a)+j}|_{[a,b[}$.

The last refinement of the partitions is needed to guarantee the following properties.
\begin{prop}[Mixing partitions]\label{mainpropertiespartitionprop}
Let $T\in \M^+$, $C^+>C^-$. Given $\epsilon>0$ and $0<\delta<1$, there exists $M(\epsilon)$, $t_3>t_2$ and refined partial partitions $\eta_m (t) \subset \eta_{s}(t)$, where $\eta_{s}(t)$ are the partitions given by Proposition \ref{growthBSf'partition}, such that $Leb(\eta_m(t))> Leb(\eta_{s}(t))-2\delta$ and for $t\geq t_3$ and each $I=[a,b[ \, \in \eta_m(t)$, $J(I,t)\rightarrow + \infty$ as $t\rightarrow +\infty$ and the following properties hold: 
\begin{enumerate}

\item \emph{Uniform vertical distribution}: 
\be\label{uniformstretch}
\left|\frac{ h(R)(b-a) }{\Delta f^j(I,t)} - Leb(I^{h(R)}_j) \right| \leq \frac{ h(R)(b-a) }{\Delta f^j(I,t)} \epsilon , \qquad j=1,\dots, J(I,t)-2; 
\ee

\item \emph{Variation of slopes}: 
\be \label{variationslopes}
\left| \frac{\Delta f(I,t)}{\Delta f^j(I,t)} -1 \right|     \leq  \epsilon, \quad j=0,\dots, J(I,t)\!-\!2; 
\ee

\item \emph{Asymptotic number of curves}:
\be \label{asymptoticscurves}
\left| \frac {J([a,b[,t)-1}{\Delta f(I,t)} -1 \right| \leq \epsilon ; 
\ee

\item \emph{Equidistribution on the base}: for some $\overline{x}\in I$,
\be \label{equidistributionbase}
\left| \frac{1}{J([a,b[,t)-1} \left( \sum_{j=0}^{J([a,b[,t)-2} \chi
 \left( T^{r(a)+j}(\overline{x}) \right)\right) - (b_2-b_1 - 2\delta) \right|   \leq    \epsilon . 
\ee
\end{enumerate}
\end{prop}

\subsection{Area estimates.}\label{areaestimatesec}
Let us show that the properties in Proposition \ref{mainpropertiespartitionprop} are enough to deduce the estimate (\ref{mainestimate}) of the Mixing Criterium (Lemma \ref{mixingcriterium}) and hence conclude the proof of the Theorem \ref{mixing}. 
\begin{lemma}\label{areaboundlemma}
For each $I=[a,b[ \in \eta_{m}(t)$, $t\geq t_3$ and $x\in I$,
\be\label{areabound}
Leb ([a,b[ \cap \varphi^{-t} R) \geq 
\sum_{j=1}^{J([a,b[,t)-2} \chi
\left( T^{r(a)+j}x \right) Leb \left( I_j ^{h(R)} \right) . 
\ee
\end{lemma}
\begin{proof}
Assume that $j_0$ is such that $1\leq j_0 \leq J([a,b[,t)-2$ and $\chi (T^{r(a)+j_0}x) =1 $, i.e. $T^{r(a)+j_0}x \in [b_1+\delta,b_2-\delta]$. Recalling that $t_2$ was chosen so that $\sup_{I\in \eta_{m}(t)} Leb(I) < \delta$ for $t\geq t_2$, when  $t\geq t_3 \geq t_2$, we have $T^{r(a)+j_0} I \subset [b_1,b_2]$. 

It is enough to show that $I_{j_0} ^{h(R)} \subset [a,b[ \cap \varphi^{-t} R$ to conclude. 
If $y\in I_{j_0} ^{h(R)}$, by definition $t - h(R)  \leq \BS{f}{r(a)+j_0}(y) \leq t$. Also, recalling that $h(R)< m_f$,
\bes
t \leq t - h(R) + f(T^{r(a)+j_0}y)  < \BS{f}{r(a)+j_0+1}(y) ,
\ees
which shows that $r(y,t)=r(a)+j_0$ and also that $y\in [y_{j_0} ,  y_{j_0+1}] \subset [a,b[$, by monotonicity 
and definition of splitting points (\ref{splittingpointsdef}). 

It follows, by definition (\ref{Ijhdef}) of $I_{j_0} ^{h(R)}$ and (\ref{flowdef}) of the flow action, that
\bes
\varphi_t(y,0) = (T^{r(a)+j_0}y , t-   \BS{f}{r(a)+j_0 }(y)   ) \in [b_1,b_2] \times [0,h(R)]= R.
\ees
This shows that $I_{j_0} ^{h(R)} \subset [a,b[ \cap \varphi^{-t} R$.
\end{proof}
Let us estimate the RHS of (\ref{areabound}). For $t\geq t_3$,
\begin{eqnarray}
&& \sum_{j=1}^{J([a,b[)-2} \chi
\left( T^{r(a)+j}(x) \right) Leb \left( I_j ^{h(R)} \right) \geq \nonumber \\
&\geq  & (1-\epsilon) h(R) (b-a) \sum_{j=1}^{J([a,b[,t)-2} \chi
\left( T^{r(a)+j}(x) \right) 
\frac{1}{\Delta f^j} \geq   \label{useuniformstretch} \\
&\geq & (1-\epsilon)^2 h(R) (b-a) \sum_{j=1}^{J([a,b[)-2} \chi
\left( T^{r(a)+j}(x) \right) 
\frac{1}{\Delta f}
 \geq   \label{usevariationslopes}
 \\
&\geq & (1-\epsilon)^3 h(R) (b-a) \sum_{j=1}^{J([a,b[)-2} \frac{ \chi
\left( T^{r(a)+j}(x) \right) }{J([a,b[,t)-1}  
\geq  \label{useasymptoticscurves} \\
&\geq & (1-\epsilon)^3 h(R) (b-a) (b_2-b_1 - 2\delta - 2\epsilon ) \xrightarrow{\epsilon , \delta \rightarrow 0} \mu(R) (b-a)  \label{useequidistributionbase} .
\end{eqnarray}
We used, in order, Property $1$ 
to get (\ref{useuniformstretch}), 
Property $2$ 
to get (\ref{usevariationslopes}),
Property $3$ to get (\ref{useasymptoticscurves}) and, eventually, to get (\ref{useequidistributionbase}) we combined 
Property $4$ 
with $\chi ( T^{r(a)}x )/(J(I,t)-1)\leq \epsilon$ for $t\geq t_3$ if $t_3$ is enlarged if necessary, since $J(I,t)$ tends to infinity.

When $\epsilon$ and $\delta$ are chosen sufficiently small, 
together with the Lemma \ref{areaboundlemma}, this concludes the proof of (\ref{mainestimate}). From Lemma \ref{mixingcriterium}, we get Theorem \ref{mixing}.

\subsection{Preliminary partitions.}\label{preliminarypartitionsec}
Let us prove Proposition \ref{preliminarypartitionprop}.  

Consider a fixed continuous time $t$. The maximum number of discrete iterations of $T$ when flowing by $t$, i.e. $r_M(t) \doteqdot \sup_{x\in I^{(0)}} r(x,t)$, can be bounded from above for each $x$ by using that $f\geq m_f >0$ and the definition of $r(x,t)$. We get 
\bes r(x,t) m_f \leq \BS{f}{r(x,t)} (x) \leq t.
\ees
Recalling the definition (\ref{RMdef}), we get 
$r_M(t)+1\leq R_M(t)$.

\paragraph{Continuity intervals of controlled size.} \label{continuitysec}
It is easy to see that any iterate $T^n\ \doteqdot T \cdot \phantom{,}\! \dots \phantom{,}\! \cdot T$ is again an IET: denoting by $\beta_0 =0 < \beta_1 < \! \phantom{,}\dots \!\phantom{,} < \beta_{d-1} <1$ the discontinuities of $T$, the discontinuities of $T^N$ are 
\be \label{discontinuities}
\{ T^{-j} \beta_i |\quad i=0, \dots , d-1; \quad  0\leq j < N \}.
\ee
Remark that $T^N$ is an exchange of at most $N d +1$ intervals.

Let $\eta_0(t)$ be the \emph{partition} of $I^{(0)}$ \emph{into continuity intervals} for $T^{R_M(t)}$, i.e. the partition into semi-open intervals whose endpoints coincide with the set (\ref{discontinuities}) where $N=R_M(t)$. By construction, for each $0\leq j\leq R_M(t)$, $T^j$ restricted to any $[a,b[ \, \in \eta_0(t)$ is continuous.

Given $M>1$, consider the following set
\bes
U_1 \doteqdot \bigcup _{\begin{subarray}{c} 0\leq i \leq d \\ 0\leq j \leq R_M(t) \end{subarray}} \overline{Ball}(T^{-j}\beta_i , \frac{2M}{t \log \log t} ),
\ees
which consists of closed balls of radius $2M/t \log \log t$ centered at the endpoints of $\eta_0(t)$.
Let $\eta_1(t)$ be the partial partition obtained from $\eta_0(t)$ by throwing away all intervals completely contained in $U_1$. Since, using (\ref{RMdef}),
\be \label{measureS1tendstozero}
Leb(U_1) \leq \frac{4 M}{t \log \log t} d \left( \frac{t}{m_f} + 3 \right) \xrightarrow{t \rightarrow + \infty} 0,
\ee
it follows that 
$Leb(\eta_1(t))\geq  1- Leb(U_1)$ 
converges to one.
Moreover, by construction, each $I\in \eta_1(t)$ contains at least one $y\notin U_1$. Hence, since the endpoints of $I$ are centers of the balls in $U_1$, $Leb(I) \geq 4M/t \log \log t$.

\paragraph{Distance from singularities.}
Let 
\bes
U_2 \doteqdot \bigcup _{0\leq j \leq R_M(t)}
T^{-j}\left[0 , \frac{M}{t \log \log t}\right)  \cup 
  \bigcup _{0\leq j \leq R_M(t)}  T^{-j}\left[1- \frac{M}{t \log \log t},1 \right).
\ees
Let $\eta_2(t) = \eta_1(t) \backslash U_2$. By construction, if $x\in \eta_2(t)$, 
\bes dist (T^{s}x, 0 ) \geq \frac{M}{t \log \log t}, \quad  dist (T^{s}x, 1 ) \geq \frac{M}{t \log \log t}, \qquad 0\leq  s\leq R_M(t),
\ees
which is Property $2$ of Proposition \ref{preliminarypartitionprop}.
Similarly to (\ref{measureS1tendstozero}), also $Leb(U_2)\xrightarrow{t\rightarrow \infty} 0$. 
Given $\delta >0$, choose $t_0$ so that $Leb(\eta_2(t))\geq Leb(\eta_1(t))- Leb(U_2)> 1-\delta/2$ for $t\geq t_0$. 
Intervals $I \in \eta_2(t)$ are either intervals of $\eta_1(t)$ or are obtained by some $I'\in \eta_1(t)$ by cutting an interval of length at most $M/(t \log \log t)$ on one or both sides of $I'$. 
Hence, 
$Leb(I') \geq 2M/t \log \log t$.

Let $\widetilde{\eta_2}(t)$ be the union of intervals of the form $[a,b'[\subset [a,b[$ associated to each $[a,b[\in \eta_2(t)$. Choosing each $b'$ close enough to each $b$, one still has $Leb([a,b'[)> 2M/t \log \log t$ and $Leb(\widetilde{\eta_2}(t)) > 1-\delta/2$. Since $T^j$, for $0\leq j \leq  R_M(t)$, is continuous on $[a,b']$, Property $1$ holds for $[a,b'[\, \in \widetilde{\eta_2}(t)$.

Construct $\eta_3(t)$ from $\widetilde{\eta}_2(t)$ by cutting each of the intervals $I \in \widetilde{\eta}_2(t)$ in pieces which satisfy the Lengths Control Property $2$ of Proposition \ref{preliminarypartitionprop}. 
For example, cut first $[Leb(I)/ (1/t\log \log t)] -1 $ intervals of length exactly $1/(t\log \log t)$ starting from the left, so that the last remaining interval has length at most $2/t\log \log t$.

Properties $1$ and $3$ still hold and 
$Leb(\eta_3(t)) = Leb(\widetilde{\eta}_2(t)) >1-\delta/2$ for $t\geq t_0$.


\paragraph{Control of the number of discrete iterations.}
Let us bound from below $r(x,t)$ when $x\in \eta_3(t)$. As a consequence of Property $3$,
\be \label{controlf}
f(T^j x) \leq const \log (t \log \log t) , \qquad 0\leq j \leq R_M(t). 
\ee
Hence $\BS{f}{r(x,t)+1}(x) \leq (r(x,t)+1)  const \log ( t \log \log t)$ and since by definition of $r(x,t)$, we have $\BS{f}{r(x,t)+1}(x)>t$, 
\be \label{roughboundbelowforr}
r(x,t) \geq \frac{t}{const \log (t \log \log t)} -1 \xrightarrow{t \rightarrow + \infty} + \infty,
\ee
uniformly for all $x\in \eta_3(t)$.
Since $f\in L^1$ and $T$ is ergodic, by Birkhoff ergodic theorem, for each $\delta >0$ there exists a measurable set $E_{\delta}$ and $N_{\delta}>0$ such that $Leb (E_{\delta})< \delta/2$ and 
\bes
\left| \frac{1}{r} \BS{f}{r}(x) - \int f(s) \ud s \right|  < 1, \qquad \forall \, x \notin E_{\delta}, \quad r\geq N_{\delta}.
\ees
Define a refined partial partition
$ \eta_4(t) \doteqdot  \eta_3(t) \backslash \{ I  \in \eta_3(t) | \quad I \subset  E_{\delta} \}.$
For $t\geq t_0$, 
 $Leb( \eta_4(t) ) \geq Leb(\eta_3(t))-\delta/2\geq  1 - \delta$. By construction for each $I \in \eta_4(t)$ there is at least one $x_I \in I$ such that $| \frac{1}{r} \BS{f}{r}(x_I) - 1| < 1$ for all $r\geq N_{\delta }$.
Enlarging $t_0$ if necessary, by (\ref{roughboundbelowforr}) we can assure $r(x_I,t) > N_{\delta}$ for each $x_I$, $I \in \eta_4(t)$. Hence 
$\BS{f}{r(x_I,t)+1}(x_I)< 2(r(x_I,t)+1)$, which, together with $\BS{f}{r(x_I,t)+1}(x_I) >t $ gives
\be \label{rxI}
 r(x_I,t) > t/2-1 , \quad \forall x_I , \, I \in \eta_4(t).
\ee 

To control all other $r(x,t)$, $x \in \eta_4(t)$, let us estimate the variation $r(x,t)-r(x_I,t)$ when $x \in I\in \eta_4(t)$. 
Assume $r(x,t) < r(x_I,t)$, otherwise we already have the lower bound.
By Properties 1 and 3, 
 $\BS{f}{r}$ is continuous on $I$ and $(T^r)'=1$ for $0\leq r\leq R_M$, so $(\BS{f}{r}(x))'= \BS{f'}{r}(x)$. By mean value theorem there exists $z$ between $x_I$ and $x$ such that
\bes
\left| \BS{f}{r(x,t)}(x_I) - \BS{f}{r(x,t)}(x) \right| \leq 
| \BS{f'}{r(x,t) } (z)  |
|x-x_I|.
\ees
Apply the rough bound on $\BS{f'}{r}$ in Corollary \ref{BSr^2}, enlarging again $t_0$ by (\ref{roughboundbelowforr}) so that $r(x,t)\geq r_0$ for $t\geq t_0$. Combining it with Property $3$ already proved, which gives $1/x_m \leq t \log\log t /M$ and $1/(1-x_M) \leq t \log\log t /M$, we get that $|\BS{f'}{r(x,t)}(z)| \leq const \, t^2$ for $t\geq t_0$. Since $Leb(I)\leq 2/( t \log \log t)$,
\be \label{roughboundvariation}
\left| \BS{f}{r(x,t)}(x_I) - \BS{f}{r(x,t)}(x) \right| \leq \frac{const \, t}{\log \log t}.
\ee
Hence, using (\ref{roughboundvariation}) and $\BS{f}{r(x_I,t)}(x_I)\leq t$ 
 and then $ \BS{f}{r(x,t)}(x)> t - f(T^{r(x,t)} x ) $ and (\ref{controlf}),
\bes
\begin{split}
&(r(x_I,t) - r(x,t)) m_f  \leq  \BS{f}{r(x_I,t)}(x_I) -  \BS{f}{r(x,t)}(x_I)  \leq \\
&\, t - \BS{f}{r(x,t)}(x) + \frac{const \, t}{\log \log t} \leq 
const \log(t \log t \log t ) + \frac{const \, t}{\log \log t}   =  o(t).
\end{split}
\ees
Rearranging and using the control for $x_I$ given by (\ref{rxI}), $r(x,t) \geq r(x_I,t) -o(t) \geq t/2 -1 - o(t)$. 
Hence, recalling also $r(x,t)\leq R_M(t)\leq  t/m_f +2$, if $t_0$ is large enough
, when $t\geq t_0$,  for each $x \in I \in \eta_4(t)$,
$t/3 \leq r(x,t) \leq R_M(t) \leq 2t/m_f$,
which is Property $4$. Since 
also the other Properties still hold, setting $\eta_p(t)\doteqdot \eta_4(t)$ this proves Proposition \ref{preliminarypartitionprop}.



\subsection{Stretching partitions.}\label{stretchingpartitionsec}
Let us prove Proposition \ref{growthBSf'partition}. 

Let $T\in \M$. For each $t$, let $l(t)$ and $l'(t)$ be uniquely determined by
\be \label{tbetweenh}
h^{(n_{l(t)})} \leq  R_M(t) < h^{(n_{l(t) + 1})}; \qquad h_\I^{(n'_{l'(t)})} \leq  R_M(t) < h_\I^{(n'_{l'(t) + 1})},
\ee
where  $\{n_l\}_{l\in \mathbb{N}}$ and $\{n'_{l'}\}_{l'\in \mathbb{N}}$ are the sequences of balanced times given by Proposition \ref{existencebalancedtimes} for $T$ and $T^{\I}$ respectively.
\begin{lemma}\label{controlhbetweenrlemma}
There exists $L\in \mathbb{N}$ independent of $t$ such that if $\frac{t}{3}\leq r \leq R_M(t)$ then
\be\label{controlhbetweenr}
h^{(n_{l(t)-L})} \leq  r < h^{(n_{l(t) + 1})}; \qquad h_\I^{(n'_{l'(t)-L})} \leq  r < h_\I^{(n'_{l'(t) + 1})}
\ee
\end{lemma}
\begin{proof} 
Let $l\in \mathbb{N}$ be such that $d^{l}\geq \max \{ 6\kappa /m_f, 6\kappa' /m_f \} $\footnote{Recall that $\kappa$ and $\kappa'$are given by Proposition \ref{existencebalancedtimes} and Corollary \ref{BSv}.}. 
By Property (\ref{positive}) in Proposition \ref{existencebalancedtimes}, we can apply Lemma \ref{productpositive} considering products of positive matrices that appear every $\overline{l}$ balanced steps and get, recalling the choice of $l$, by balance of the induction steps and (\ref{tbetweenh}),
$h^{(n_{l(t)-l\overline{l}})}$
$\leq (\kappa/d^l) \min_j h_j^{(n_{l(t)})} \leq  (\kappa/d^l)(2t/m_f) \leq t/3 \leq r$.
Analogous expressions can be obtained also for $h_\I^{(n'_{l'(t)-l\overline{l}})}$ and show that setting $L\doteqdot l\overline{l}+1$ we get (\ref{controlhbetweenr}).
\end{proof}

Define the set  $\Sigma_t = \Sigma_t (T)$ as 
\be\label{defSt}
\Sigma_t \doteqdot \bigcup_{l=l(t)-L}^{l(t)-1} \Sigma_l^+(T) \cup \bigcup_{l=l'(t)-L}^{l'(t)-1} \Sigma^-_{l'}(T) \cup  \overline{\Sigma}^+_{l(t)}(T) \cup \overline{\Sigma}^-_{l'(t)}(T),
\ee
where the sets $\Sigma^+_l(T)$ and  $\Sigma^-_{l'}(T)$ were defined in (\ref{badsetforBSu}) and (\ref{DefSigma-}) and where
\bes
 \overline{\Sigma}^+_{l(t)}(T) \doteqdot  \bigcup_{i=0}^{\min\{ R_M(t), [\sigma_{l(t)} h^{(n_{l(t)+1})}]\}} T^{-i} [0, \sigma_{l(t)} \lambda^{(n_{l(t)})}].
\ees
and $\overline{\Sigma}^-_{l'(t)}$ is the analogous truncation of $\Sigma^-_{l'(t)}$. Remark that
\be\label{LebSigmat}
Leb(\Sigma_t) \xrightarrow{t\rightarrow +\infty} 0,
\ee
because $Leb(\Sigma^+_l)\xrightarrow{t\rightarrow +\infty} 0$ for each $l\geq l(t)-L$ from (\ref{LebSigma}), the same holds for each $Leb(\Sigma^-_{l'})$, $l'\geq l'(t)-L$ and $\Sigma_t$ is union of at most $2(L+1)$ such sets. 
\begin{proofof}{Proposition}{growthBSf'partition}
Fix $T$, $C^+>C^-$, $\delta>0$ and $M>1$ and let $\eta_p(t)$ be the preliminary partitions given by Proposition \ref{preliminarypartitionprop} for $t\geq t_0$. 
Consider the set $\Sigma_t (T)$ defined in (\ref{defSt}) and, by (\ref{LebSigmat}), choose $t_1\geq t_0$ so that $Leb(\Sigma_t)<\delta/2$ for $t\geq t_1$.  
Define $\eta_{s}(t)$ as the partition obtained from $\eta_p(t)$ throwing away all the intervals which intersect $\Sigma_t$. 
If $I\in \eta_p(t)$ and $I\cap \Sigma_t\neq \emptyset$, then, from Property $1$ in Proposition \ref{preliminarypartitionprop}, either $I\subset \Sigma_t$ or, for some $0\leq j\leq R_M(t)$, $T^j I$ contains either some points $\sigma_{l} \lambda^{(n_{l})}$ with $l(t)-L\leq l \leq l(t)$ or some $1-\sigma_{l'} \lambda_{\I}^{(n_{l'})}$ with $l'(t)-L\leq l' \leq l'(t)$. Hence, using (\ref{LebSigmat}), Property $2$ in Proposition \ref{preliminarypartitionprop} and bounding the number of such points, we get 
\bes
Leb(\eta_{s}(t)) \geq Leb(\eta_p(t))-Leb(\Sigma_t)- \frac{2}{t \log \log t} 2 (L+1)\frac{2t}{m_f} 
\ees
and, enlarging $t_1$ if necessary, both the last two terms are than $\delta/2$.

Let $t/3  \leq r \leq 2 / m_f $ and $x\in \eta_{s}(t)$. Let us show that the assumptions of Corollary \ref{BSf'cor} and Corollary \ref{growthBSf''} on the growth of $\BS{f'}{r}$ and $\BS{f''}{r}$ hold. 
By Lemma \ref{controlhbetweenrlemma}, there exists $l, l'$, with $l(t)-L \leq l \leq l(t)$ and $l'(t)-L\leq l' \leq l'(t)$ such that $h^{(n_{l})} \leq  r < h^{(n_{l + 1})}$ and $ h_\I^{(n'_{l'})} \leq  r < h_\I^{(n'_{l'+1})}$. Since by construction of $\eta_s(t)$, $x\notin \Sigma_t$, in particular, if $l<  l(t) $ and $l'<  l'(t)$, $x\notin \Sigma_l^+$, $x\notin \Sigma^-_{l'}$.
Hence in this case the assumptions (\ref{assumptionsBSf'}) of Corollaries \ref{BSf'cor} and \ref{growthBSf''} hold. In the case where $l= l(t) $ or $l'=l'(t)$, we only have $x\notin \overline{\Sigma}^+_l$ or $x\notin \overline{\Sigma}^-_{l'}$, but 
also in this case the Corollaries hold since the only property needed in their proof is that $T^i x \notin [0,\sigma_{l(t)} \lambda^{(n_{l(t)})}]$ or $[1-\sigma_{l'(t)} \lambda_\I^{(n'_{l'(t)})}]$ for $0\leq i\leq r$ and $r \leq R_M(t)$. 

Since $r\geq t/3 \geq t_1/3$, enlarging $t_1$, one can assure that $r\geq r_o$, for the $r_o$ given by the Corollaries \ref{BSf'cor} and \ref{growthBSf''} and since $x\in I \in \eta_p(t)$, $x$ is not a singularity of $\BS{f}{r}$ by Property $1$ and $3$ of Proposition \ref{preliminarypartitionprop}.  Hence one can apply the Corollaries \ref{BSf'cor} and \ref{growthBSf''}. Moreover, by Property $3$, 
\bes
x_m \doteqdot \min_{0\leq i<r} T^i x \geq \frac{M}{t\log \log t}  \quad \mathrm{and} \quad 1-x_M \doteqdot \min_{0\leq i<r}(1- T^i x)  \geq \frac{M}{t\log \log t}. 
\ees
Thus, we get respectively
\begin{eqnarray}
 \BS{f'}{r}(x)  &\leq &
- C_1 r\log r \left( 1 - \frac{C^+ (\kappa'+1) t\log \log t}{C_1  M r \log r}\right) ;\label{BSf'forpartition}
\\ 
| \BS{f'}{r}(x) | &\leq &
 C_2 r\log r  \left( 1 + \frac{(\kappa + \kappa' +2) t\log \log t}{ M r \log r}\right) ;\label{BSf'aboveforpartition}
\\
| \BS{f''}{r}(x) | & \leq & 
C_4 \frac{t \log \log t}{M} r\log r \left( 1 +     \frac{(\kappa' +\kappa+2) t\log \log t}{ M r \log r} \right) \label{BSf''forpartition} .
\end{eqnarray}
 Recalling that $t\leq 3r$ and again by enlarging $t_1$ if necessary, one can assure that the last terms in (\ref{BSf'forpartition}, \ref{BSf'aboveforpartition}, \ref{BSf''forpartition}), involving $t\log \log t / r \log r$, are less than $1/2$. Hence we get respectively (\ref{f'decreases}, \ref{SBf'boundabove}, \ref{BSf''}).
\end{proofof}


\begin{cor}\label{growthBSpartitioncor}
For each $I=[a,b[\in \eta_{s}(t)$, $x\in I$ and $r(a,t)\leq r \leq r(b,t)$,
\begin{eqnarray}
&& const ( t\log t )  \leq   |\BS{f'}{r}(x)|   \leq   const' ( t\log t ) ; \label{growthf'} \\
&& \Delta f (I , t)  \geq   const \left( \frac{\log t }{\log \log t} \right) \xrightarrow{t \rightarrow + \infty } + \infty. \label{growthDeltaf}
\\
&& | \Delta f (I , t)|  =   o(\log t); \qquad  | \Delta f^1 (I , t)|  =   o(\log t) .\label{deltafologt}
\end{eqnarray}
\end{cor}
\begin{proof} Equation (\ref{growthf'}) follows from (\ref{f'decreases}) and (\ref{SBf'boundabove}) of Proposition \ref{growthBSf'partition}, since $t/3 \leq r(a,t),r(b,t) \leq 2t/m_f$. 
Since $\BS{f}{r(a,t)}$ is continuous with its derivative on $[a,b]$ (Property 1 and 3 in Proposition \ref{preliminarypartitionprop}), by mean value theorem there exists $z\in I$ such that 
\bes
 \Delta f = \BS{f}{r(a,t)}(a) - \BS{f}{r(a,t)}(b) = - \BS{f'}{r(a,t)}(z)(b-a) \geq   \frac{const ( t\log t )}{2 t \log \log t} ,
\ees 
where we applied (\ref{f'decreases}) 
 and the control on the interval sizes (Property 2 in Proposition \ref{preliminarypartitionprop}).
The proof of (\ref{deltafologt}) are obtained similarly using (\ref{growthf'}).
\end{proof}

\subsection{Mixing partitions.}\label{mainpropertiespartitionsec}
Let us prove Proposition \ref{mainpropertiespartitionprop}. 
\paragraph{Outline.}
The uniform vertical distribution (\ref{uniformstretch}) is proved in Section \ref{uniformstretchsec}.
In order to show that $J(I,t)$
is tending to infinity and moreover that it is asymptotic to $\Delta f$ (i.e. (\ref{asymptoticscurves})), we need a further refinement of the partitions that guarantees that, exploiting ergodicity, the number of fibers covered in time $\Delta f$ is asymptotic to $\Delta f$, since $\int_{I^{(0)}}f(x)\ud x =1$.

The refinement is constructed in two steps. In Section \ref{upperlowerboundsec} we prove a first rough upper bound (Corollary \ref{Jupperboundlemma}), which allows us to refine the partitions in order to have a better control of the distance from singularities for $r(a)\leq r \leq r(a)+1$, i.e. in the curves range (Lemma \ref{curvesdistsinglemma}). Using this refinement one can prove the control on the variation of the curves slopes (Section \ref{variationslopessec}).

In Section \ref{asymptoticscurvesandequidistributionsec}, the same refinement is used to give a lower bound on the number of curves and hence to construct a second refinement of the partitions on which one can get by ergodicity both equidistribution along the base (\ref{equidistributionbase}) and the curve asymptotics (\ref{asymptoticscurves}).

\subsubsection{Uniform vertical distribution.}\label{uniformstretchsec}
Let us show that given $\epsilon >0$, choosing $M>M_0(\epsilon)$, each $I\in \eta_{s}(t)$ satisfies Property $1$. 
%
Let us recall the following definition used in \cite{Fa:ana} (see also \cite{Ko:mix,Ko:nonI,Ko:nonII}).
\begin{defn}
Given $\epsilon>0$ 
, the function $g$ on the interval $[a,b]$ is \emph{
$\epsilon$-uniformly distributed} if 
for any $c,d$ such that $\inf_{[a,b[}g \leq c \leq d \leq \sup_{[a,b[}g$,
the measure of the set 
$I_{c,d} = \{x \in [a,b[ \, | \quad c\leq g(x) \leq d \}$ 
satisfies
\be \label{measIcd}
(1-\epsilon) \frac{ d-c  }{\sup_{[a,b[}g - \inf_{[a,b[}g}  \leq \frac{Leb(I_{c,d})}{b-a} \leq (1+\epsilon) \frac{ d-c  }{\sup_{[a,b[}g - \inf_{[a,b[}g} .
\ee
\end{defn}

In \cite{Fa:ana} Fayad proves the following criterium for uniform distribution.
\begin{lemma}[Fayad]\label{Fayadunifstretchlemma}
If $g$ is monotonic and 
\be
\sup_{[a,b[}|g''(x)||b-a| \leq \epsilon \inf_{[a,b[}|g'(x)| ,
\ee
then $g$ is $\epsilon$-uniformly distributed on $[a,b]$.
\end{lemma}

For each $[a,b[ \, \in \eta_{s}(t)$, consider $\BS{f}{r(a)+j}$, for $j=1,\dots, J([a,b[,t)-2$. 
From Property 2 of Proposition \ref{preliminarypartitionprop}, (\ref{f'decreases}, \ref{BSf''}) of Proposition \ref{growthBSf'partition} and $r\leq 2t/m_f$,
\bes \begin{split}
\frac{\sup_{[a,b[}|\BS{f''}{r(a)+j}(x)||b-a|}{\inf_{[a,b[}|\BS{f'}{r(a)+j}(x)| } 
\leq \frac{  C'' r^2 (\log r) ( \log \log r )(b-a)  }{M  C_1' r \log r }  \leq \\
\leq \frac{C'' (\frac{2t}{m_f} \log \log \frac{2t}{m_f})  \frac{2}{  t \log \log t}}{M C_1'}   \xrightarrow{t\rightarrow + \infty}  \frac{4 C'' }{ M C_1' m_f} .
     \end{split}
\ees
Choosing  $M > M_0(\epsilon ) \doteqdot 8 C''/ C_1' m_f \epsilon$ and then $t_3\geq t_2$ large enough, the last expression is less than $\epsilon$ for $t\geq t_3$.
Hence, each $\BS{f}{r(a)+j}$, for $j=1,\dots, J\!-\! 2$, being also decreasing, is $\epsilon$-uniformly distributed on $[a,b[\, \in \eta_{s}(t)$ by Lemma \ref{Fayadunifstretchlemma}.
The set $I_{j}^{h(R)}$ defined in (\ref{Ijhdef}) is of the form $I_{c,d}$ for $g=\BS{f}{r(a)+j}$, $c= t -h(R)$ and $d=t$; by the definition (\ref{splittingpointsdef}) of splitting points, one can check that $d=g(y_{j})\leq \sup_{[a,b[} g $ and $c \geq g(y_{j+1})\geq \inf_{[a,b[} g$. 
Hence, by (\ref{measIcd}), its measure is bounded by
\be
(1-\epsilon ) \frac{h(R)}{\Delta f^j ([a,b[,t) } (b-a)\leq  Leb(I_{j}^{h(R)} ) \leq (1+\epsilon ) \frac{h(R)}{\Delta f^j ([a,b[,t) } (b-a) .\ee
This proves the uniform vertical distribution Property (\ref{uniformstretch}).

\subsubsection{Rough upper and lower bound on the number of curves.}\label{upperlowerboundsec}
Recall that the number of curves generated from each $\varphi_t(I)$  is given by $J(I,t)=r(b,t)-r(a,t)+1$. 
\begin{lemma}\label{equivalentJlemma} For each $I=[a,b[$,
\be\label{equivalentJ}
r\left( T^{r(a)}b,\Delta f\right) + 1 \leq J(I,t)\leq  r\left( T^{r(a)+1}b,\Delta f^1\right) + 3.
\ee
\end{lemma}
Lemma \ref{equivalentJlemma} shows that the number of strips is related to the number of fibers that the point $T^{r(a)}b$ still has to cover in time $\Delta f$  when $a$ stops, because of the delay accumulated through the stretching of $\BS{f}{r(a)}$.



\begin{proof}
Applying the relation $\BS{f}{r_1+r_2}(x)=\BS{f}{r_1}(x)+\BS{f}{r_2}(T^{r_1}x)$ and the definition (\ref{defr}) of $r(\cdot,\cdot)$,
\bes
\begin{split}
 \BS{f}{ r ( T^{r(a)} b  , \Delta f)  + r(a) }(b)& =  \BS{f}{r(a)} (b) + \BS{f}{r( T^{r(a)} b ) , \Delta f)} (T^{r(a)}b) \leq \\ & \leq \BS{f}{r(a)} (b) + \Delta f = \BS{f}{r(a)} (a) \leq t.  \end{split}
\ees
Hence, $r(b,t) \geq r ( T^{r(a)} b  , \Delta f  ) + r(a,t) $, which is the first inequality in (\ref{equivalentJ}).

For the second inequality, 
\bes
\begin{split}
 \BS{f}{ r ( T^{r(a)+1 } b  , \Delta f^1)  + r(a) +2 }(b)& =  \BS{f}{r(a)+1} (b) + \BS{f}{r( T^{r(a)+1} b  , \Delta f^1 )+1 } (T^{r(a)+1}b)  \\ & > \BS{f}{r(a)+1} (b) + \Delta f^1 = \BS{f}{r(a)+1} (a) > t, 
\end{split}
\ees
which implies that $r(b)<  r ( T^{r(a)+1 } b  , \Delta f^1)  + r(a) +2 $.
\end{proof}

\paragraph{Rough upper bound on the number of curves.}
\begin{cor} 
\label{Jupperboundlemma}
Let $J(t)\doteqdot \sup_{I\in \eta_{s}(t)} J(I,t)$. Then
\be\label{Jupperbound}
J(t) =  o(\log t).
\ee
\end{cor}
\begin{proof}
By Lemma \ref{equivalentJlemma},
$J(I,t) \leq   r ( T^{r(a)+1 } b  , \Delta f^1) +3\leq \frac{\Delta f^1(I,t)}{m_f} +3$. 
Recalling that $\Delta f^1(I,t) = o\left( \log t \right)$  by (\ref{deltafologt}) of Corollary \ref{growthBSpartitioncor} we get the bound.
\end{proof}

\paragraph{Pull-back using discrete times.}
Consider the map $R_t: I^{(0)}\rightarrow  I^{(0)}$ given by $R_t(x)= T^{r(x,t)}x $, which is the projection of $ \varphi_t(x,0)\in X_f$ to the base $I^{(0)}$.
Remark that $R_t$ in general is not one to one.
The following lemma is used by Kochergin in \cite{Ko:mix} (Lemma 1.3).
\begin{lemma}[Kochergin]
For any measurable set $S \subset I^{(0)}$,
\be \label{Kochlemma}
Leb(R_t^{-1} S) \leq \int_S \left( \frac{f(x)}{m_f} +1 \right) \ud x.
\ee
\end{lemma}
Since $f\in L^1$, by absolute continuity of the integral, for any $\delta>0$ it is possible to choose $\delta_1$ such that the RHS of (\ref{Kochlemma}) is bounded by $\delta$ as long as $Leb(S)< \delta_1$. Hence we get the following corollary.
\begin{cor}\label{Kochlemmacor}
For each $\delta >0$, there exists $\delta_1>0$ such that for any measurable $S \subset I^{(0)}$ , if $Leb(S)< \delta_1$, then
$Leb(R_t^{-1} S) < \delta $.
\end{cor}

\paragraph{Refinement of the partitions.}
We want to refine the partitions in order to throw away points $x$ such that when considering $r$ such that $r(a)\leq r \leq r(a)+ J$, i.e. in the range of the curves, the distance of $T^{r}x$ from singularities is bounded from below by $1/\log t$. 
\begin{lemma}\label{curvesdistsinglemma}
There exist partitions $\eta_{5}(t)\subset \eta_{s}(t)$ and $t_3\geq t_2$ such that, for $t\geq t_3 $, $Leb(\eta_5(t))\geq Leb(\eta_{s}(t))-\delta$ and for each $x\in I \in \eta_5(t)$, 
\be \label{curvesdistsing}
|T^r x| \geq \frac{1}{ (\log t)^2 };\quad |1- T^r x| \geq \frac{1}{(\log t)^2 },
\ee
for each $r(a,t)\leq r \leq r(a,t)+J(I,t)$.
\end{lemma}
\begin{proof}
Define 
\be
U_3(t) \doteqdot \bigcup_{i=-[J(t)]}^{[J(t)]} T^i \left( \left[ 0, \frac{1}{(\log t)^2} \right] \right)  \cup  \bigcup_{i=-[J(t)]}^{[J(t)]} T^i \left( \left[ 1- \frac{1}{(\log t)^2} ,1 \right] \right).
\ee
Since the continuity intervals for $T^{[J(t)]}$ and  $T^{-[J(t)]}$ are at most $d([J(t)]+1)$ (see Section \ref{continuitysec}),  
the set $U_3(t)$ consists of at most $O(J(t)^2)$ disjoint intervals. 
Consider the $2/t\log\log t$-neighbourhood of $U_3(t)$, i.e. let $U_4(t)\doteqdot \{ \, x \in I^{(0)}|\,\, d(x, U_3(t))\leq 2/t\log\log t \, \}$. 
Hence, using Corollary \ref{Jupperboundlemma}, 
\bes
Leb(U_4(t))\leq \frac{4J(t)+4}{(\log t)^2} + \frac{O(J(t)^2)}{t \log \log t}
\xrightarrow{t \rightarrow + \infty } 0.
\ees
Choosing $t_3>t_2$ so that for $t\geq t_3$,  $Leb(U_4) < \delta_1$ where $\delta_1$ is given by Corollary \ref{Kochlemmacor}, we get $Leb(R_t^{-1}(U_4) )< \delta$. Define a refined partition $\eta_5(t)\subset \eta_s(t)$ by
\bes
 \eta_5(t) \doteqdot  \eta_{s}(t) \, \backslash \, \{ I  \in \eta_{s}(t) | \quad I \subset  R_t^{-1} U_4(t) \}.
\ees
Clearly $Leb(\eta_5(t))\geq Leb(\eta_{s}(t)) - \delta$.
Let us show that for each $I\in  \eta_5(t)$, we get (\ref{curvesdistsing}). 
By construction there exists $x\in I$ such that $R^t (x) = T^{r(x,t) }x \notin  U_4(t)$. Hence, by Proposition \ref{preliminarypartitionprop}, $T^{r(x,t) }y \notin  U_3(t)$ for each $y\in I$.  
For each $r=r(a),\dots,r(a)+ J$, the point $T^r y$ satisfies the inequalities (\ref{curvesdistsing})
by definition of $U_3(t)$, because $T^{r(x,t)}y \notin U_3(t)$, as shown above, and  $|r(x,t)-r|\leq J$.
\end{proof}

\begin{lemma}[Rough lower bound on $J(t)$]\label{roughlowerboundJlemma} 
Let $J(t)\doteqdot \sup_{I\in \eta_{5}(t)} J(I,t)$.
\be\label{roughlowerboundJ} 
J(t)\geq O\left( \frac{\log t}{ (\log \log t)^2}\right) \xrightarrow{t\rightarrow + \infty} + \infty.
\ee
\end{lemma}
\begin{proof}
For each $I=[a,b[\, \in \eta_5(t)$, by Lemma \ref{curvesdistsinglemma}, $f(T^{r}b) \leq O(  \log (\log t))$ for each $r(a,t) \leq r \leq r(a,t) + J(I,t)$. Hence, since by Lemma \ref{equivalentJlemma}, $J(I,t)\geq r ( T^{r(a)} b  , \Delta f )+1$, 
\bes 
J(I,t)  \geq  \frac{\BS{f}{ r ( T^{r(a)}b  , \Delta f)+1 }(T^{r(a)}b )}{\max_{0\leq i <  r ( T^{r(a)}b  , \Delta f)+1 }  f(T^{r(a)+i}b)} \geq \frac{\Delta f}{O( \log (\log t))},
\ees 
which gives (\ref{roughlowerboundJ}) by using the bound (\ref{growthDeltaf}) on $\Delta f$. 
\end{proof}

\subsubsection{Variation of slopes.}\label{variationslopessec}
Given $I=[a,b[ \in \eta_5(t)$ the variation of the average slope of the curves $(t- \BS{f}{r(a)+j})|_{I}$, for $0\leq j\leq J(t)$ can be written as
\bes
\begin{split}
\left| \Delta f^j - \Delta f \right| & =
\left| \BS{f}{r(a)+j }(a) -  \BS{f}{r(a)+j }(b) -  \BS{f}{r(a)}(a) + \BS{f}{r(a)}(b) \right| \\
&= \left| \sum_{i=0}^{j-1} f(T^{r(a)+i} a)  -  \sum_{i=0}^{j-1} f(T^{r(a)+i} b ) \right| \leq   \sum_{i=0}^{j-1} \left| f'(T^{r(a)+i}c)\right|  (b-a) ,
\end{split}
\ees
where in the last estimate $a\leq c \leq b$ by mean value theorem. 
Using Lemma \ref{curvesdistsinglemma}, $|f'(T^{r(a)+i}c)| \leq O \left( (\log t)^2 \right) $ for each $0\leq i<j \leq  J(t)$. 
Hence, applying also the bound on $J(t)$ given by Corollary \ref{Jupperboundlemma}, the growth estimate (\ref{growthDeltaf}) for $ \Delta f$ and the size control of $(b-a)$ (Property 2 of Proposition \ref{preliminarypartitionprop}),
\bes
\left| \frac{\Delta f^j - \Delta f}{\Delta f} \right| \leq \frac{ J(t) \sup_{i=0}^{J-1} |f'(T^{r(a)+i}c)| (b-a)}{ \Delta f  } \leq  \frac{ o(\log t ) O \left( (\log t)^2 \right) }{ O \left( \frac{\log t }{\log \log t} \right) t \log \log t  } ,
\ees
which converges to zero as $t \rightarrow + \infty$. Enlarge $t_3 > 0 $ so that the RHS is less than $\epsilon$ for $t\geq t_3$ to get (\ref{variationslopes}).

\subsubsection{Equidistribution on the base and asymptotic of curves.}\label{asymptoticscurvesandequidistributionsec}
Both the equidistribution on the base and the exact asymptotic for the number of curves follow by proving uniform convergence on a large set for the Birkhoff sums of $\chi$ 
and $f$ respectively. More precisely, one seeks uniform control for the points of the form $T^{r(a)}b$ where $b$ are the endpoints of the partition intervals $[a,b[$. 


\paragraph{Equidistribution on the base.}
Let $\widetilde{\eta}_5(t)$ be a narrowing of $\eta_5(t)$ obtained keeping only the central third of each interval:
\bes
\widetilde{\eta}_5(t) \doteqdot  \left\{ \left. \left[ a + (b-a)/3 
 , b - (b-a)/3  
 \right[ \, \right|\quad  [a,b[ \in \eta_5(t)   \right\} .
\ees

For each $\varepsilon>0$ and $\delta_1>0$, by ergodicity of $T$ and $T^{-1}$ one can find $U_5$ and $N>0$ such that $Leb(U_5)< \delta$ and for each $x\notin U_5$ and $n\geq N$, 
\be \label{unifconvBS}
\left| \frac{\BS{f, T^i}{n}(x)}{n} -1  \right| < \varepsilon ; \qquad  \left| \frac{\BS{\chi, T^{i}}{n}(x)}{n} -(b_2-b_1-2\delta)  \right| < \varepsilon; \qquad i=0,-1. 
\ee
If $\delta_1$ is given by Corollary \ref{Kochlemmacor} in correspondence of $\delta/3>0$, we get $Leb(R^{-t}(U_5))<\delta/3$.  
Define $\eta_6(t)\subset \eta_5(t)$ by throwing away all intervals $I \in \eta_5(t) $ such that the corresponding $\widetilde{I}$ is completely contained in $R^{-t}(U_5)$. Hence,
$Leb(\eta_6(t))\geq Leb(\eta_5(t))-\delta  \geq Leb(\eta_{s}(t))-2\delta $ by Lemma \ref{curvesdistsinglemma}.

By construction, for each $I\in \eta_6(t)$, there exists $\overline{x}$ such that $|\overline{x}-a|$, $|\overline{x}-b|> (b-a)/3$ and $T^{r(\overline{x},t)}\overline{x} \notin U_5$ and hence satisfies (\ref{unifconvBS}).
Arguing as in Corollary \ref{growthBSpartitioncor} to prove (\ref{growthDeltaf}), both 
$\Delta f ([a,\overline{x}],t )$ and $\Delta f ([\overline{x},b],t )$, as $t\rightarrow \infty$, are bounded from below by $O ( \log t / \log \log t)$. 
As in Lemma \ref{equivalentJlemma}, 
\bes
r(\overline{x})-r(a) \geq r(T^{r(a)}\overline{x}, \Delta f ([a,\overline{x}],t ); \qquad r(b)-r(\overline{x}) \geq  r(T^{r(\overline{x})}b ,\Delta f ([\overline{x},b],t )).
\ees
Hence, by the same proof in Lemma \ref{roughlowerboundJlemma}, both $r(\overline{x})-r(a)$ and $r(b)-r(\overline{x})$ tend to infinity uniformely as $t$ increases.
Enlarge $t_3$ so that for $t\geq t_3$, both
$r(\overline{x})-r(a) > N$ and $r(b)-r(\overline{x}) > N$.
Hence, the estimates in (\ref{unifconvBS}) hold when $n=r(\overline{x})-r(a)$ or $r(b)-r(\overline{x})$ and $x=T^{r(\overline{x},t)} \overline{x}$. 
Moreover they also hold for $i=-1$ and $x=T^{r(\overline{x},t)-1} \overline{x}$. To see it, in the case of $f$, use that
\bes
 \BS{f,T^{-1}}{r(\overline{x})-r(a)}(T^{r(\overline{x},t)-1}\overline{x}) = \BS{f ,T^{-1}}{r(\overline{x})-r(a)+1}(T^{r(\overline{x},t)} \overline{x}) - f(T^{r(\overline{x},t)} \overline{x}) 
\ees
and from Lemma \ref{curvesdistsinglemma} and the analogous of Lemma \ref{roughlowerboundJlemma} for $ \Delta f ([a,\overline{x}],t )$,
\bes \frac{f(T^{r(\overline{x},t)}\overline{x}) }{r(\overline{x})-r(a) } \leq O\left( \log (\log t)^2 \frac{(\log \log t)^2}{\log t}\right) 
\ees
which can be made arbitrarly small enlarging $t_3$ if necessary.
In the case of $\chi$, just use that $\chi\leq 1$ and $r(\overline{x})-r(a)$ tends to infinity.

Let us combine these estimates decomposing the Birkhoff sums as
\bes 
\BS{\chi,T}{r(b)-r(a)} (T^{r(a)} \overline{x} ) 
= \BS{\chi,T^{-1}}{r(\overline{x})-r(a)} (T^{r(\overline{x})-1} \overline{x} ) + \BS{\chi,T}{r(b)-r(\overline{x})} (T^{r(\overline{x})} \overline{x} )
\ees
and using that $\frac{r(b)-r(\overline{x})}{r(b)-r(a)} + \frac{r(\overline{x})-r(a)}{r(b)-r(a)}  =1 $. We get
\be\label{BSdecompcons} 
 \left|  \frac{\BS{\chi,T}{r(b)-r(a)}(T^{r(a)}\overline{x})}{r(b)-r(a)} -(b_2-b_1-2 \delta)  \right|\leq 2\epsilon ,
\ee
which proves equidistribution on the base (\ref{equidistributionbase}) for $\overline{x}\in I$.

\paragraph{Asymptotic number of curves.}\label{asymptoticscurvessec}
\begin{lemma}\label{convergenceBSrab} Enlarging $t_3$ if necessary, for each $[a,b[ \, \in \eta_6(t)$, $t\geq t_3$,
\bes \label{BSforb}
\left|  \frac{1}{r(b)-r(a)}  \BS{f}{r(b)-r(a)} (T^{r(a)} b)  -1  \right| \leq  2 \varepsilon.
\ees
\end{lemma}
\begin{proof}
By mean value, there exists $z \in [\overline{x},b]$ such that
\bes
\begin{split}
\left| \BS{f}{r(b)-r(\overline{x})} (T^{r(\overline{x},t)} \overline{x}) - \BS{f}{r(b)-r(\overline{x})} (T^{r(\overline{x},t)} b) \right| \leq  \left| \BS{f'}{r(b)-r(\overline{x})}(T^{r(\overline{x},t)}z) \right|\cdot \\ \cdot (b-\overline{x}) \leq J(t) \sup_{r(\overline{x})  \leq i\leq r(b)} |f'(T^{i}z) | (b-\overline{x}) \leq 
  \frac{ o(\log t) O((\log t)^2) }{t \log \log t}
\end{split}
\ees
where we used Corollary \ref{Jupperboundlemma} to bound $J(t)$ , Lemma \ref{curvesdistsinglemma} to bound $|f'(T^{i}z)|$ and Property 2 in Proposition \ref{preliminarypartitionprop} to control the size $(b-a)$. 
Hence, enlarging $t_3$, from the analogous estimate for $T^{r(\overline{x},t)}\overline{x}$, we get for $t \geq t_3$, 
\be \label{BSTrxbpast}
\left| \frac{1}{r(b)-r(\overline{x})}\BS{f}{r(b)-r(\overline{x})}(T^{r(\overline{x},t)}b) -1  \right| \leq 2 \varepsilon .
\ee
In a similar way, from the analogous estimate for  $T^{r(\overline{x},t)-1}\overline{x}$, we get
\be \label{BSTrxbfuture} 
\left| \frac{1}{r(a)-r(\overline{x})}\BS{f,T^{-1}}{r(a)-r(\overline{x})}(T^{r(\overline{x},t)-1}b) -1  \right|  < 2 \varepsilon .
\ee
Combining (\ref{BSTrxbpast}) and (\ref{BSTrxbfuture}) and decomposing the Birkhoff sums as 
\bes
\BS{f,T}{r(b)-r(a)} (T^{r(a)} b )=\BS{f,T^{-1}}{r(\overline{x})-r(a)} (T^{r(\overline{x})-1} b ) + \BS{f, T}{r(b)-r(\overline{x})} (T^{r(\overline{x})}b)
\ees
we get the Lemma.
\end{proof}

\begin{lemma}\label{reductionasymptoticslemma}
Enlarging $t_3$ if necessary, for each $[a,b[\, \in \eta_6(t)$, if $t \geq t_3$,
\bes
\left| \frac{ \BS{f}{r(b)-r(a)} (T^{r(a)}b)}{\Delta f} - 1\right|  \leq 
\varepsilon .
\ees
\end{lemma}
\begin{proof}
Since we can rewrite
\bes
\begin{split}
 \BS{f}{r(b)-r(a)} (T^{r(a)}b) &= \BS{f}{r(b)} (b) -  \BS{f}{r(a)} (b) \pm \BS{f}{r(a)} (a) = 
\\ &= \Delta f + \BS{f}{r(b)} (b) -  \BS{f}{r(a)} (a),
\end{split}
\ees
from $t - f( T^{r(a)} a ) <  \BS{f}{r(a)} (a) \leq t $ and $t-   f( T^{r(b)} b) <  \BS{f}{r(b)} (b)\leq t $, we get
\bes
\left| \frac{ \BS{f}{r(b)-r(a)} (T^{r(a)}b)}{\Delta f} - 1\right| 
\leq \frac{\max \{ f( T^{r(b)} b) , f( T^{r(a)}a) \}}{\Delta f} \leq 
\frac{ O((\log  \log t)^2)}{\log t}\rightarrow 0,
\ees
where we used by Lemma \ref{curvesdistsinglemma} and (\ref{growthDeltaf}).
\end{proof}

Lemma  \ref{convergenceBSrab} and Lemma \ref{reductionasymptoticslemma} give
\bes
\left| \frac{ r(b,t)-r(a,t)}{ \BS{f}{r(b)-r(a)} (T^{r(a)}b)} \frac{ \BS{f}{r(b)-r(a)} (T^{r(a)}b)}{\Delta f} - 1\right| \leq (1+\epsilon)
\ees
for an appropriate choice of $\varepsilon$ and $t\geq t_3$. Recalling that $J([a,b[,t)-1 =  r(b,t)-r(a,t) $, this concludes the proof of 
 the asymptotic of curves (\ref{asymptoticscurves}).

Setting $\eta_{m}(t)= \eta_6(t)$, this completes the verification that the partitions $\eta_{m}(t)$, for an appropriate choice of $\delta$ and $t\geq t_3$, satisfy  all the Properties listed in Proposition \ref{mainpropertiespartitionprop}. 

\subsection*{Acknowledgments.}
I would like to thank my advisor Prof. Ya. G. Sinai for proposing the problem and constantly guiding me with suggestions, patience and encouragement.
I also would like to thank both K. Khanin and A. Avila for useful discussions. 
Thanks also to Pavel Batchourine for listening to some parts of this work.


\bibliography{biblioflows}
\bibliographystyle{alpha}

\end{document}